\documentclass{article}

\newtheorem{fed}{\textbf{Definition}}[section]
\newtheorem{thm}[fed]{\textbf{Theorem}}
\newtheorem{lemma}[fed]{\textbf{Lemma}}

\newtheorem{prop}[fed]{\textbf{Proposition}}
\newtheorem{cor}[fed]{\textbf{Corollary}}
\usepackage{amssymb,bbold, bbm,graphicx,epsfig,psfrag,epic,eepic,latexsym}
\usepackage{amsmath}
\usepackage{mathrsfs}
\usepackage{wrapfig, graphicx}
\usepackage[dvips]{color}

\begin{document}

\title{The moduli space of gradient flow lines and Morse homology}
\author{Urs Frauenfelder, Robert Nicholls}
\maketitle

\tableofcontents

\newpage

\section{Introduction}

Given a Morse function on a closed manifold one constructs a chain complex generated by the critical points of the Morse function and graded by their Morse index. The boundary operator is obtained by counting gradient flow lines between critical points. A priori it is far from obvious that with this recipe one obtains actually a boundary operator, i.e., that the boundary of the boundary vanishes. That this is actually the case rests on a careful analysis of the moduli space of gradient flow lines. 

While the complex depends on the choice of the Morse function the resulting homology is independent of it. Critical points can be born and die but homology never dies nor is born. In particular, one has a stable lower bound on the number of critical points. An immediate consequence of this are the Morse inequalities.

For people just interested in the Morse inequalities on a closed manifold the construction of Morse homology which involves 
quite heavy tools from functional analysis might be a bit too much and the wonderful book by Milnor \cite{milnor} gives a quicker way to the Morse inequalities. On the other hand the amazing work of Floer \cite{floer1,floer2,floer3} shows that there is a semi-infinite dimensional analogon of Morse homology referred to as Floer homology. For example periodic orbits of a Hamiltonian system can be obtained variationally as critical points of the action functional of classical mechanics and Floer showed how they can be used as generators of a chain complex. The Morse inequalities for Floer's chain complex than lead to topological lower bounds on the number of periodic orbits conjectured before by Arnold. 

Although we do not discuss Floer homology in these notes they are written having Floer homology in mind. Therefore we use concepts which can be generalized to the semi-infinite dimensional case. In the finite dimensional case they might not always be the most efficient ones but on the other hand the finite dimensional case has the advantage that these concepts can be visualized much more easily than in the semi-infinite dimensional one giving the reader the right intuition how Floer homology works. 

These notes grew out of lectures the first named author gave at Ludwig-Maximilian University in Munich, Seoul National University, and the University of Augsburg. He would like to thank all the participants of his lectures which gave him all a lot of inspiration. Of course there are many excellent textbooks and articles on these topics from which these notes profited a lot, like \cite{audin-damian, banyaga-hurtubise, nicolaescu, salamon, schwarz, weber}.

\section{Morse functions and Morse indices}

\subsection{Critical points}

Assume that $M$ is a manifold and $f \colon M \to \mathbb{R}$ is a smooth function. A \emph{critical point} of the function
$f$ is  a point $u \in M$ such that the differential of $f$ at this point vanishes, i.e., $df(u)=0$. Alternatively suppose that we are in a local coordinate chart around $u$, then we can rephrase this condition that all partial derivatives at $u$ vanish
$$\frac{\partial f}{\partial u_i}(u)=0, \quad 1 \leq i \leq n$$
where $n$ is the dimension of the manifold. Note that the concept of a critical point does not depend on the choice of coordinates. Indeed, if $u=u(v)$ is a coordinate change we obtain by the chain rule
$$\frac{\partial f}{\partial v_j}(v)=\frac{\partial f}{\partial u_i}(u)\frac{\partial u_i}{\partial v_j}(v)=0, \quad
1 \leq j \leq n.$$

\subsection{The Hessian}

The Hessian at a point $u$ in a coordinate chart is the symmetric $n \times n$-matrix given by
$$H_u=\bigg[\frac{\partial^2 f}{\partial u_i \partial u_j}(u)\bigg]_{1 \leq i,j \leq n}.$$
Let us see how the Hessian transforms under coordinate change. By the chain rule we obtain
$$\frac{\partial^2 f}{\partial v_i \partial v_j}(v)=\frac{\partial^2 f}{\partial u_k \partial u_\ell}(u)\frac{\partial u_k}{\partial v_j}(v) \frac{\partial u_\ell}{\partial v_i}(v)+\frac{\partial f}{\partial u_k}(u)\frac{\partial^2 u_k}{\partial v_j \partial v_i}(v).$$
This looks a bit scary. However, note that if $u$ is a critical point, the last term just vanishes so that the above formula simplifies to 
$$\frac{\partial^2 f}{\partial v_i \partial v_j}(v)=\frac{\partial^2 f}{\partial u_k \partial u_\ell}(u)\frac{\partial u_k}{\partial v_j}(v) \frac{\partial u_\ell}{\partial v_i}(v).$$
If we abbreviate by
$$\Phi=\bigg[\frac{\partial u_j}{\partial v_i}(v)\bigg]_{1 \leq i,j \leq n}$$
the Jacobian of the coordinate change we obtain the following transformation behaviour for the Hessian at a critical point
$$H_u=\Phi^T H_v \Phi$$
where $\Phi^T$ denotes the transpose of the Jacobian. That means that the Hessian transforms as a symmetric form under coordinate change.

\subsection{Invariants of a symmetric form and the Theorem of Sylvester}

We recall a classical theorem due to Sylvester \cite{sylvester} about invariants of a quadratic form. 
Suppose that $H$ is a $n \times n$-matrix. We can associate to $H$ the following three nonnegative integers
\begin{eqnarray*}
n_0&:=&\mathrm{dim} \ker H,\\
n_+&:=&\max\Big\{\mathrm{dim}V: V \subset \mathbb{R}^n\,\,\mathrm{subvectorspace},\,\,v^T H v>0,\,\,\forall \,\,v 
\in V \setminus \{0\}\Big\},\\
n_-&:=&\max\Big\{\mathrm{dim}V: V \subset \mathbb{R}^n\,\,\mathrm{subvectorspace},\,\,v^T H v<0,\,\,\forall \,\,v 
\in V \setminus \{0\}\Big\}.
\end{eqnarray*}
If $H$ is symmetric, i.e.,
$$H=H^T$$
we can find an orthogonal matrix $\Phi$ which diagonalizes $H$, i.e.,
$$\Phi^{-1}H \Phi=D$$
where $D$ is diagonal. Because $\Phi$ is orthogonal we can write this equation equivalently as
$$\Phi^T H \Phi=D.$$
By scaling we can arrange that all entries of the diagonal matrix $D$ are either $1$, $-1$ or $0$. In this last step we usually loose the property that $\Phi$ is orthogonal, however, $\Phi$ is still invertible. In view of this we obtain the following relation between the three invariants
\begin{equation}\label{sylvester}
n_-+n_0+n_+=n.
\end{equation}
Moreover, we see that there are no other invariants of a quadratic form. 

\subsection{Morse critical points and its index}

Suppose that $p \in M$ is a critical point of a smooth function $f \colon M \to \mathbb{R}$ and let $H_f(p)$ be the Hessian
of $f$ at $p$. We say that the critical point is \emph{Morse} if 
$$\ker H_f(p)=\{0\}$$
or equivalently if
$$n_0(H_f(p))=0.$$
A smooth function is called Morse, if all its critical points are Morse. 
Suppose that $p$ is a Morse critical point of $f$. Then we define its Morse index 
$$\mu(p):=n_-\big(H_f(p)\big)$$
namely the number of negative eigenvalues of the Hessian counted with multiplicity. If $n$ is the dimension of the manifold
$M$, then we infer from (\ref{sylvester}) that
$$n_+\big(H_f(p)\big)=n-\mu(p).$$
\begin{figure}
\begin{center}
\includegraphics{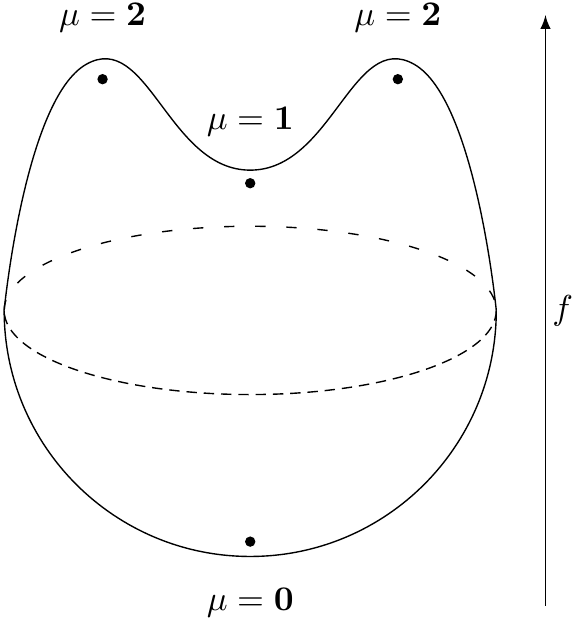}
\caption{Distorted sphere with Morse critical points and Morse indices.}
\end{center}
\end{figure}
In fact apart from the value the function $f$ attains at a Morse critical point, its Morse index is the only local invariant of the function as the following Lemma tells
\begin{lemma}[Morse]\label{morselemma}
Suppose that $p$ is a Morse critical point of a function $f \colon M \to \mathbb{R}$. Then there exist local coordinates 
$u$ around $p$ such that $p$ lies in the origin and the function $f$ in this coordinates becomes
$$f(u)=f(p)-\sum_{i=1}^{\mu(p)}u_i^2+\sum_{i=\mu(p)+1}^n u_i^2.$$
\end{lemma}
\textbf{Proof: } A proof of the Morse Lemma can for example be found in the classical book by Milnor on Morse theory, 
\cite[Chapter 2]{milnor}. \hfill $\square$

\section{Examples of computations of Morse homology}

\subsection{A cooking recipe}

Suppose that $M$ is a closed manifold and $f \colon M \to \mathbb{R}$ is a Morse function. Abbreviate
$$\mathrm{crit}(f):=\big\{x \in M: df(x)=0\big\}$$
the set of critical points of $f$. For $k \in \mathbb{N}$ we set
$$\mathrm{crit}_k(f):=\big\{x \in \mathrm{crit}(f): \mu(x)=k\big\}$$
the subset of critical points having Morse index $k$. We declare the set $\mathrm{crit}_k(f)$ to be the basis of
a $\mathbb{Z}_2$-vector space, namely we set
$$CM_k(f):=\mathrm{crit}_k(f) \otimes \mathbb{Z}_2.$$
Vectors $\xi \in CM_k(f)$ are formal sums
$$\xi=\sum_{x \in \mathrm{crit}_k(f)} a_x x$$
where the coefficients $a_x$ belong to the field $\mathbb{Z}_2$. 
We further define a boundary operator
$$\partial_k \colon CM_k(f) \to CM_{k-1}(f).$$
We require that $\partial_k$ is linear. Therefore it suffices to define it on basis vectors. On a basis vector
$x \in \mathrm{crit}_k(f)$ it is given by
$$\partial_k x=\sum_{y \in \mathrm{crit}_{k-1}(f)}\#_2\big\{\textrm{gradient flow lines from}\,\,x\,\,\textrm{to}\,\,y\big\}y$$
where $\#_2$ denotes the cardinality modulo two. 

We remark that instead of counting gradient flow lines modulo two one could also count them with integer coefficients. However, in this case one has to do a signed count. One way to define a sign for gradient flow lines is to choose a an orientation of stable and unstable manifolds. This is for example explained in \cite{salamon}. Another way is to define coherent orientations in the sense of Floer and Hofer \cite{floer-hofer}. This is for example explained in \cite{schwarz}.

In the following we illustrate the cooking recipe in a couple of examples. 

\subsection{The heart}

\begin{figure}[h]
\begin{center}
	\includegraphics{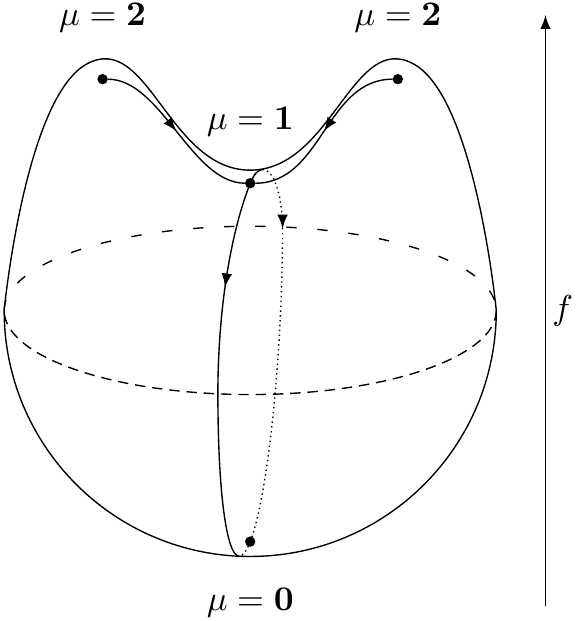}
	\caption{Contributing gradient flow lines and critical points of deformed sphere.}
\end{center}
\end{figure}
The heart will appear in these notes again and again since one can illustrate basically the whole intricate behaviour of gradient flow lines with this example. 
Topologically it is a sphere but it is distorted so that its height function has two maxima, one saddle point and one minimum. Let us denote the two maxima by $x_1$ and $x_2$, the saddle point by $y$ and 
the minimum by $z$.
Then we have
$$CM_2(f)=\langle x_1,x_2\rangle \cong \mathbb{Z}_2^2,$$
i.e., the $\mathbb{Z}_2$-vector space has the basis $x_1$ and $x_2$. Similarly,
$$CM_1(f)=\langle y \rangle \cong \mathbb{Z}_2 \qquad CM_0(f)=\langle z \rangle \cong \mathbb{Z}_2.$$
There is one gradient flow line from $x_1$ to $y$ and one gradient flow line from $x_2$ to $y$. Therefore we have
\begin{equation}\label{heartbound}
\partial_2 x_1=y=\partial_2 x_2.
\end{equation}
There are two gradient
 flow lines from $y$ to $z$ but because we count modulo two we have
$$\partial_1 y=0.$$
In particular, we have 
$$\partial_1 \circ \partial_2=0,$$
i.e., the boundary of the boundary vanishes. There are obviously no critical point of negative Morse index so that we set
$$CM_{-1}=\{0\}$$
so that $\partial_0=0$ as well. Since there are no critical points of Morse index three either we also have
$\partial_3=0$. In view of the fact that the boundary of the boundary vanishes we can define Morse homology as the quotient vector space
$$HM_k(f)=\frac{\mathrm{ker}\partial_k}{\mathrm{im}\partial_{k+1}}.$$
From (\ref{heartbound}) we see that neither $x_1$ nor $x_2$ are in the kernel of $\partial_2$, but since we are working with $\mathbb{Z}_2$-coeffients, there sum is in the kernel, i.e.,
$$\partial_2(x_1+x_2)=0.$$
Since $\partial_3$ vanishes its image is trivial so that we have
$$HM_2(f)=\langle [x_1+x_2]\rangle \cong \mathbb{Z}_2.$$
Since $\partial_1$ vanishes the saddle point $y$ is actually in the kernel of $\partial_1$ but by (\ref{heartbound}) it is
as well in the image of $\partial_2$. Therefore
$$HM_1(f)=\{0\}.$$
Finally $z$ is in the kernel of $\partial_0$ but not in the image of $\partial_1$ so that consequently
$$HM_0(f)=\langle[z]\rangle \cong \mathbb{Z}_2.$$
Summarizing our computations we obtained
$$HM_*(f)=\left\{\begin{array}{cc}
\mathbb{Z}_2 & *=0,2\\
\{0\} & \mathrm{else}.
\end{array}\right.$$
But this precisely corresponds to the singular homology of the sphere.

\subsection{The round sphere}

\begin{figure}[h]
\begin{center}
	\includegraphics{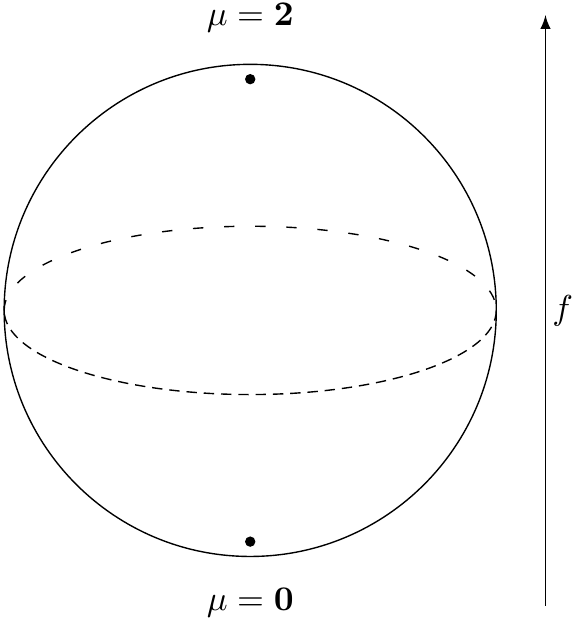}
	\caption{Gradient flow lines and critical points of round sphere.}
\end{center}
\end{figure}
The heart is definitely not the most efficient way to compute the Morse homology of the sphere. If one looks at the round sphere one obtains a much simpler Morse function having just one maximum, one minimum, and no saddle point. In particular, 
one has
$$CM_*(f)=\left\{\begin{array}{cc}
\mathbb{Z}_2 & *=0,2\\
\{0\} & \mathrm{else}.
\end{array}\right.$$
For degree reasons the boundary has to vanish, so that one obtains
$$HM_*(f)=CM_*(f)=\left\{\begin{array}{cc}
\mathbb{Z}_2 & *=0,2\\
\{0\} & \mathrm{else}.
\end{array}\right.$$
This is again the singular homology of a sphere and coincides with the computation for the heart. This two examples already illustrate the nontrivial fact that the Morse homology does not depend on the choice of the Morse function we have chosen on our closed manifold. 

\subsection{The torus}

\begin{figure}[h]
\begin{center}
	\includegraphics{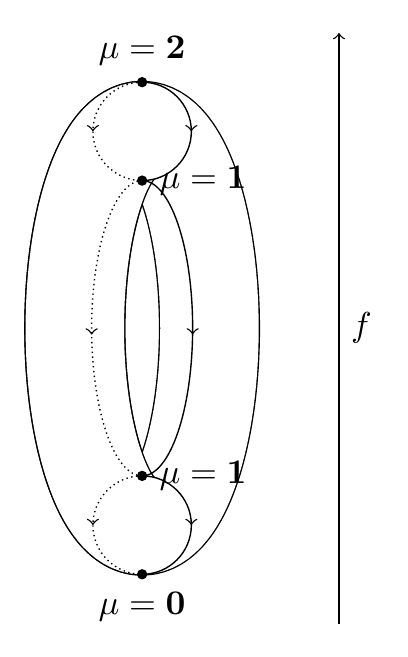}
	\quad
	\includegraphics{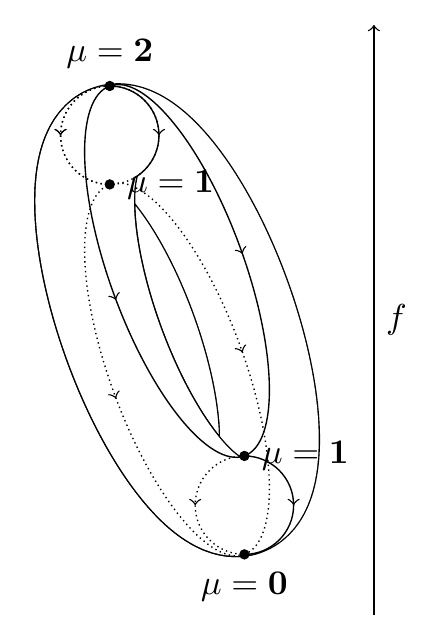}
	\caption{Upright torus with degenerate gradient flow lines (left) and tilted torus(right).}
\end{center}
\end{figure}

In our next example we consider a different topological type, namely the torus. We first put our torus upright. We see four critical points, one maximum, two saddle points, and a minimum. But then we spot something a bit strange. Namely there are as well two gradient flow lines from the higher saddle point to the lower saddle point. This is a new phenomenon. We have seen in the example of the heart gradient flow lines between critical points of index difference one but never gradient flow lines between critical points of the same index. The reader might point out that since there are two gradient flow lines and we anyway count modulo two we can just discard these. In fact this works in this example. However, in general if there occur dubious gradient flow lines one should be careful. Therefore we tilt the torus a little bit. We see that the two dubious gradient flow lines immediately disappear. We now see two gradient flow lines from the maximum to each of the saddle points
and two gradient flow lines from each saddle point to the minimum. Since we count modulo two the boundary vanishes and we have
$$CM_*(f)=HM_*(f)=\left\{\begin{array}{cc}
\mathbb{Z}_2 & *=0,2\\
\mathbb{Z}_2^2 & *=1.
\end{array}\right.$$
Again we recognise the singular homology of the torus.

The tilting of the torus can as well be interpreted as a change of the Riemannian metric on the torus. We will learn that for a generic choice of the Riemannian metric on any closed manifold we get a well-defined count of gradient flow lines leading to a boundary operator. However, there can be some ungeneric metrics with degenerate gradient flow lines which one first has to perturb a little.

\subsection{Surfaces of higher genus}

\begin{figure}
\begin{center}
	\includegraphics{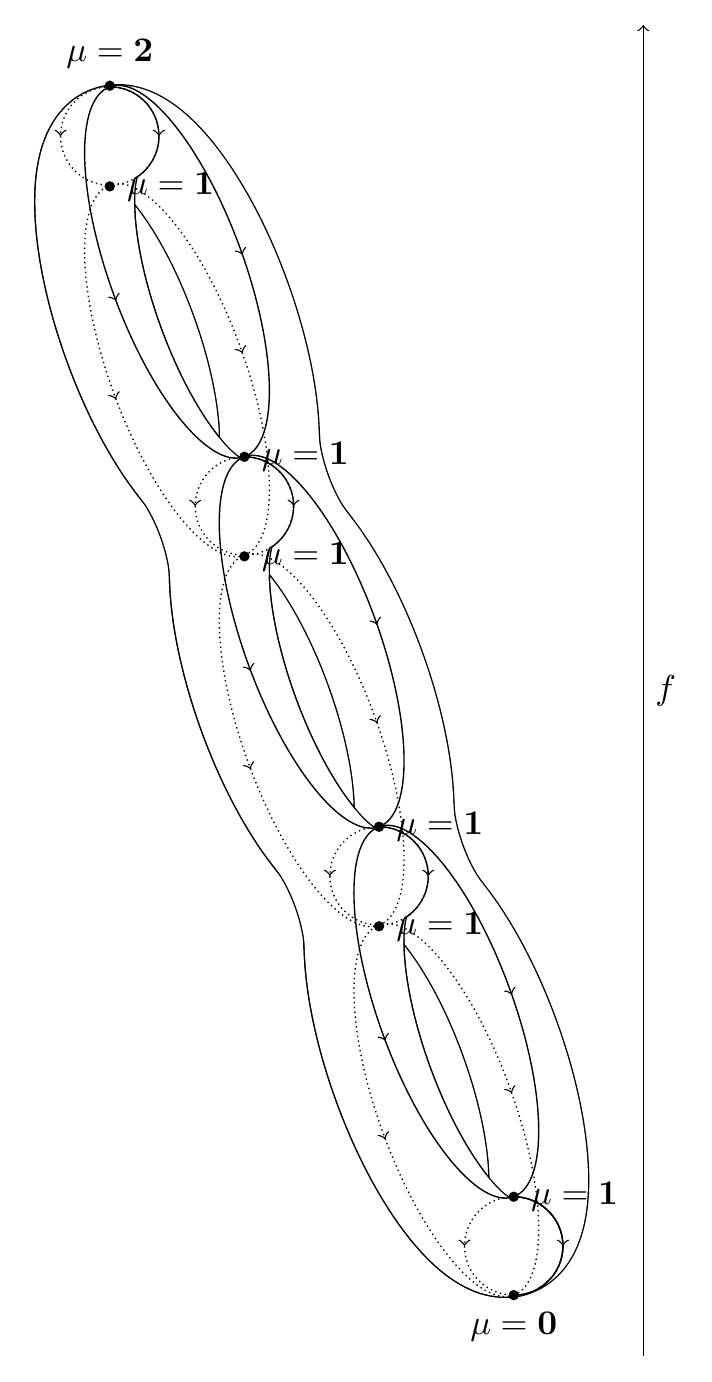}
	\caption{Gradient flow lines and critical points for surface of genus 3.}
\end{center}
\end{figure}
The example of the torus immediately generates to surfaces of higher genus. As for the torus we tilt the surface a little to avoid degenerate gradient flow lines. For a surface with $g$ holes we see one maximum, $2g$ saddle points and one minimum. 
There are two gradient flow lines from the maximum to each of the saddle points and two gradient flow lines from each saddle point to the minimum. Therefore the boundary again vanishes and we obtain
$$CM_*(f)=HM_*(f)=\left\{\begin{array}{cc}
\mathbb{Z}_2 & *=0,2\\
\mathbb{Z}_2^{2g} & *=1,
\end{array}\right.$$
namely the singular homology of a surface with $g$ holes. 

\subsection{The real projective plane}

\begin{figure}[h]
\begin{center}
	\includegraphics{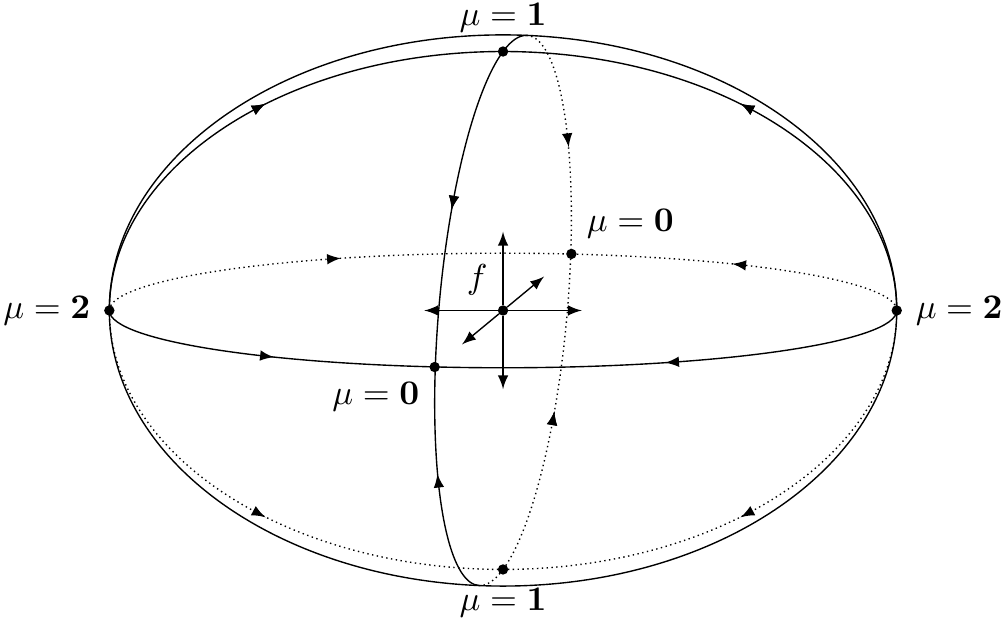}
	\caption{Critical points and gradient flow lines for ellipsoid with radial Morse function.}
\end{center}
\end{figure}
The real projective plane is the quotient space
$$\mathbb{R}P^2=S^2/\mathbb{Z}_2$$
where the group $\mathbb{Z}_2$ acts on the sphere $S^2=\{x \in \mathbb{R}^3: ||x||=1\}$ by antipodal involution
$x \mapsto -x$. A $\mathbb{Z}_2$-invariant Morse function on $S^2$ induces a Morse function on the quotient
$\mathbb{R}P^2$. For $a_1<a_2<a_3$ we consider the $\mathbb{Z}_2$-invariant Morse function
$$f \colon S^2 \to \mathbb{R}, \quad x \mapsto a_1x_1^2+a_2x_2^2+a_3x_3^2.$$
This Morse function has two maxima at $(0,0,\pm 1)$, two saddle points at $(0,\pm 1,0)$ and two minima at
$(\pm 1,0,0)$. 
We abbreviate by $u_1,u_2$ the two maxima, by $v_1,v_2$ the two saddle points and by $w_1,w_2$ the two minima. We see a gradient flow line from each of the maxima to each saddle point and a gradient flow line from each saddle point to each of the minima. Therefore on $S^2$ the chain complex of this Morse function is
$$\partial u_1=v_1+v_2=\partial u_2, \qquad \partial v_1=w_1+w_2=\partial v_2.$$
Observe that 
$$\partial^2 u_1=\partial v_1+\partial v_2=2w_1+2w_2=0$$
since we are counting modulo two and similarly $\partial^2 u_2=0$. 
Although neither $u_1$ nor $u_2$ lies in the kernel of the boundary operator their sum satisfies
$$\partial(u_1+u_2)=2v_1+2v_2=0.$$
Therefore we have 
$$HM_2(f)=\langle [u_1+u_2]\rangle \cong \mathbb{Z}^2.$$
Similarly we have $\partial(v_1+v_2)=0$ but now $v_1+v_2$ lies as well in the image of the boundary operator so that in degree one the kernel agrees with the image implying for the quotient space that
$$HM_1(f)=\{0\}.$$
In degree zero the two minima $w_1$ and $w_2$ both lie in the kernel of the boundary operator but their sum lies as well in the image so that in the quotient space we have
$$[w_1]=[w_2].$$
Consequently
$$HM_0(f)=\langle [w_1]\rangle \cong \mathbb{Z}_2.$$
Summarizing we have 
$$HM_*(f)=\left\{\begin{array}{cc}
\mathbb{Z}_2 & *=0,2\\
\{0\} & \mathrm{else}.
\end{array}\right.$$
which is again the singular homology of the sphere. 

We now consider the from $f$ induced function on the quotient
$$\bar{f} \colon \mathbb{R}P^2 \to \mathbb{R}.$$
On the quotient two two maxima are identified so that we have just one maximum left
$$u=[u_1]=[u_2] \in \mathbb{R}P^2.$$
Similarly we just have one saddle point and one minimum
$$v=[v_1]=[v_2], \qquad w=[w_1]=[w_2].$$
Gradient flow lines as well get identified in pairs. On $S^2$ we had a total of four gradient flow lines from the two maxima to the two saddle points. After identification two are left. From the two saddle points to the two minima we had as well in total four gradient flow lines on the sphere so that on the projective plane two are left as well. In particular, the boundary vanishes. Therefore
$$CM_*(\bar{f})=HM_*(\bar{f})=\left\{\begin{array}{cc}
\mathbb{Z}_2 & *=0,1,2\\
\{0\} & \mathrm{else},
\end{array}\right.$$
which corresponds to the singular homology of $\mathbb{R}P^2$.

\subsection{Noncompact examples}

For noncompact manifolds in general Morse homology cannot be defined. As an example we steal a point of the heart, namely precisely a point on one of the gradient flow lines from the saddle point to the minimum. On the punctured heart this gradient flow line does not exist anymore so that we have just one gradient flow line from the saddle point to the minimum.
Therefore we have 
$$\partial y=z$$
so that combined with (\ref{heartbound}) we get
$$\partial^2 x_1=z.$$
In particular, the boundary of the boundary does not vanish and we cannot define homology. 

Nevertheless in many examples we still can define Morse homology but the Morse homology might depend on the Morse function and not just the manifold. We illustrate this on the 1-dimensional manifold $\mathbb{R}$.

\begin{figure}[h]
\begin{center}
	\includegraphics{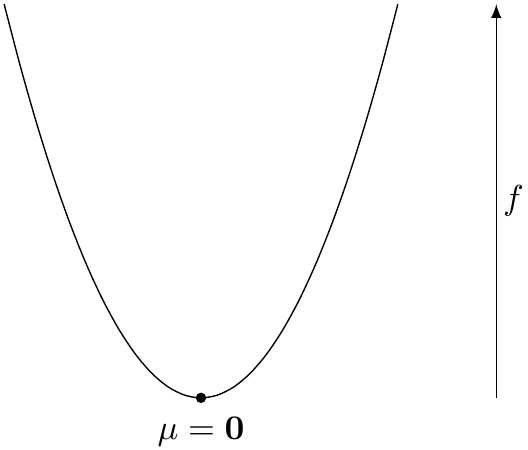} \quad
	\includegraphics{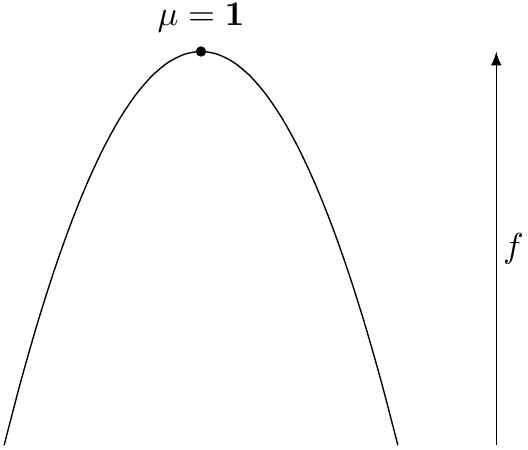}
	\caption{Morse function on the real line with only one minimum (left) and only one maximum (right).}
\end{center}
\end{figure}
Consider a parabola with one minimum. For such a Morse function we obtain
$$HM_*(f)=\left\{\begin{array}{cc}
\mathbb{Z}_2 & *=0\\
\{0\} & \mathrm{else},
\end{array}\right.$$
Now consider a parabola with one maximum. In this case we get
$$HM_*(f)=\left\{\begin{array}{cc}
\mathbb{Z}_2 & *=1\\
\{0\} & \mathrm{else},
\end{array}\right.$$

\begin{figure}[h]
\begin{center}
	\includegraphics{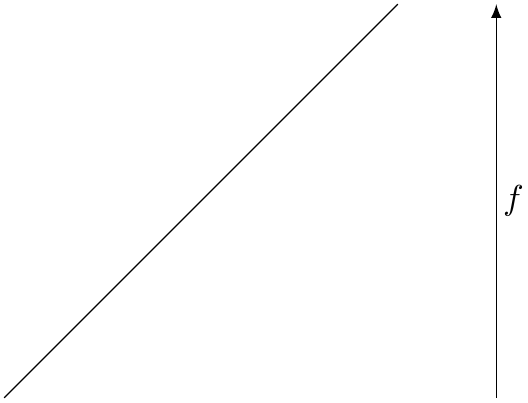}
	\caption{Morse function on the real line without any critical points.}
\end{center}
\end{figure}
If we look at a straight line with nontrivial slope there are no critical points at all and the Morse homology vanishes completely. 

Although on the real line the Morse homology depends on the Morse function as these three examples show it is nevertheless
invariant under homotopies as long as the asymptotic behaviour during the homotopy is fixed.  As an example we look at a Morse function with $n$ maxima and $n+1$ minima which goes asymptotically at both ends to plus infinity. Such a Morse function is homotopic to the parabola with one minimum through a homotopy fixing the asymptotic behaviour. 
\begin{figure}[h]
\begin{center}
	\includegraphics{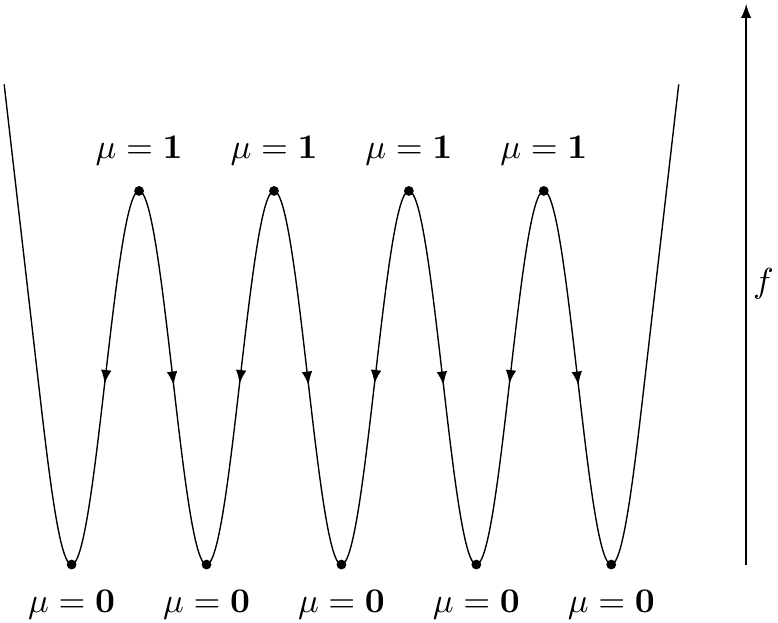}
	\caption{Morse function on the real line with many critical points.}
\end{center}
\end{figure}
We denote the maxima in increasing order by $x_1, \ldots x_n$ and the minima by $y_1, \ldots, y_{n+1}$. We observe that from each maximum two gradient flow lines go to neighbouring minima so that we have
$$\partial x_k=y_k+y_{k+1}, \quad 1 \leq k \leq n.$$
We see from this that the map
$$\partial \colon CM_1(f) \to CM_0(f)$$
is injective. In particular, the homology in degree $1$ vanishes. Since $CM_1(f)$ is $n$-dimensional it follows that
$\partial CM_1(f)$ is an $n$-dimensional subspace in the $n+1$-dimensional vector space $CM_0(f)$. We conclude that the
homology in degree zero is one-dimensional so that we have
$$HM_*(f)=\left\{\begin{array}{cc}
\mathbb{Z}_2 & *=0\\
\{0\} & \mathrm{else},
\end{array}\right.$$
as in the case of the parabola with one minimum. 

\subsection{The Morse inequalities}

Suppose that $M$ is a closed manifold and $f \colon M \to \mathbb{R}$ is a Morse function. We define the $k$-th Betti number of $M$ as the dimension of the vector space $HM_k(f)$, i.e.,
$$b_k(M):=\mathrm{dim}HM_k(f).$$
We show in these notes that the Betti numbers do not depend on the chosen Morse function on $M$. In fact they correspond to the Betti numbers defined via singular homology. This stronger fact was for example proved by Schwarz in \cite{schwarz} or with an alternative proof in \cite{schwarz2}. In particular, it follows from this, that the Betti numbers are independent of the smooth structure of our manifold and just depend on its topology, altough we need the smooth structure to actually make sense of the notion of a Morse function. An immediate consequence of this are the Morse inequalities, namely
\begin{thm}
Suppose that $f$ is a Morse function on a closed manifold $M$. Then the number of its critical points can be estimated from below by the sum of the Betti numbers of the manifold, i.e.,
$$\#\mathrm{crit}(f) \geq \sum_{k=0}^{\mathrm{dim}(M)} b_k(M).$$
\end{thm}
\textbf{Proof:} We have
\begin{eqnarray*}
\# \mathrm{crit}(f)=\sum_{k=0}^{\mathrm{dim}(M)} CM_k(f) \geq \sum_{k=0}^{\mathrm{dim}(M)} HM_k(f)
=\sum_{k=0}^{\mathrm{dim}(M)} b_k(M).
\end{eqnarray*}
This finishes the proof of the theorem. \hfill $\square$
\\ \\
For example for the sphere the Betti numbers are
$$b_k(S^2)=\left\{\begin{array}{cc}
1 & *=0,2\\
0 & \mathrm{else},
\end{array}\right.$$
so that their sum becomes two. This is not too impressive. In fact since $S^2$ is compact any continuous function on it assumes its maximum and minimum so that we already know that there have to be at least two critical points. However, the situation changes dramatically for the torus. In this case the Betti numbers are
$$b_k(T^2)=\left\{\begin{array}{cc}
1 & *=0,2\\
2 & *=1,
\end{array}\right.$$
so that their sum becomes four. Without the assumption that the function is Morse there do not need to exist four critical points. Namely identifying the torus with $\mathbb{R}^2/\mathbb{Z}^2$ the function
$$f \colon T^2 \to \mathbb{R}, \quad (x,y) \mapsto \sin(\pi x)\sin(\pi y)\sin(\pi(x+y))$$
has only three critical points, a maximum, a minimum and a degenerate saddle point also called monkey saddle. 

For a surface with $g$ holes $\Sigma_g$ the Morse inequalities get even more impressive. In this case the Betti numbers are
$$b_k(\Sigma_g)=\left\{\begin{array}{cc}
1 & *=0,2\\
2g & *=1,
\end{array}\right.$$
so that their sum becomes $2+2g$.

\section{Analysis of gradient flow lines}

In this chapter we analyse the gradient flow equation. The gradient flow equation is a first order ODE. It is invariant under time shift. When we count gradient flow line we always count unparametrized gradient flow lines, namely solutions of the gradient flow equation modulo the $\mathbb{R}$-action given by time shift. Since the group $\mathbb{R}$ is noncompact this leads to interesting analytical limit behaviour. Namely a sequence of gradient flow lines can break in the limit. 
On a closed manifold one has the following compactness result. A sequence of gradient flow lines has always a subsequence with converges to a broken gradient flow line. This compactness result is the clue why one obtains a boundary operator by counting gradient flow lines.   

\subsection{Parametrised and unparametrised gradient flow lines}

Suppose that $M$ is a smooth manifold and $f\colon M \to \mathbb{R}$ is a smooth function. Suppose that $g$ is a Riemannian metric on $M$, i.e., for every $x \in M$ we have a scalar product $g_x$ on the tangent space $T_xM$ and
$g_x$ depends smoothly on $x$. We define that the gradient of $f$ at a point $x \in M$ with respect to the metric $g$
implicitly by the condition
$$df(x)v=g_x(\nabla_g f(x),v),\quad \forall\,\,v \in T_x M.$$
If $g$ is fixed we just write for the gradient
$$\nabla f=\nabla_g f.$$
Critical points correspond to the points where the gradient vanishes
$$\mathrm{crit}f:=\{x \in M: df(x)=0\}=\{x \in M: \nabla f(x)=0\}.$$
\begin{figure}[h]
\begin{center}
	\includegraphics{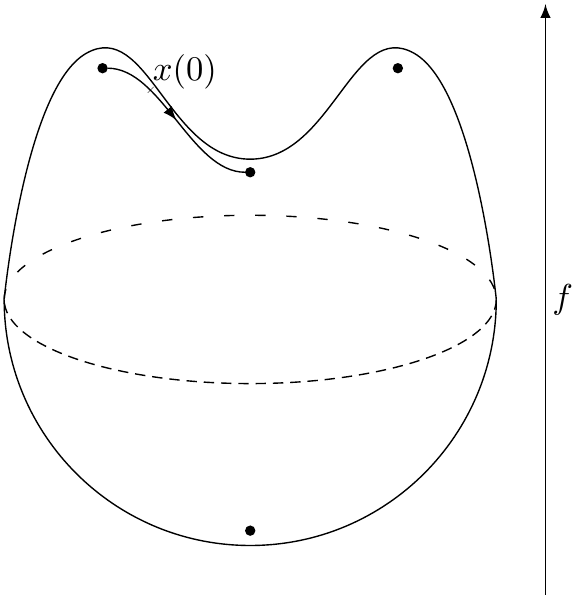}
	\caption{Parametrised gradient flow line on deformed sphere.}
\end{center}
\end{figure}
However note, that while the gradient depends on the choice of the Riemannian metric, the critical points are independent
of the Riemannian metric. 
A \emph{(parametrised) gradient flow line} $x \in C^\infty(\mathbb{R},M)$ is a solution of the ODE
\begin{equation}\label{para}
\partial_s x(s)+\nabla f(x(s))=0,\quad s \in \mathbb{R}.
\end{equation}
The group $\mathbb{R}$ acts on parametrised gradient flow lines as follows. If $r \in \mathbb{R}$ and $x$ is a solution of (\ref{para}) we define $r_* x \in C^\infty(\mathbb{R},M)$ by timeshift
$$r_* x(s)=x(s+r),\quad s \in \mathbb{R}.$$
Note that if $x$ is a solution of (\ref{para}) it follows that $r_* x$ for every $r \in \mathbb{R}$ is a solution as well.
We refer to a solution of (\ref{para}) modulo timeshift as an \emph{unparametrised gradient flow line}, i.e., an unparametrised gradient flow line is an equivalence class $[x]$ where $x$ is a solution of (\ref{para}) and the equivalence relation is given by
$$x \sim y \quad :\Longleftrightarrow\quad \exists\,\,r \in \mathbb{R},\,\,y=r_* x.$$ 
We point out that the boundary operator is defined by counting unparametrised gradient flow lines. Indeed, parametrised flow lines cannot be counted, because unless a gradient flow line is constant timeshift gives immediately rise to uncountably many others. 

\subsection{Gradient flow lines flow downhill}

A crucial property of gradient flow lines is that they flow downhill as the following lemma explains. In particular, unless they are constant, they can never come back to the same point. 
\begin{lemma}\label{downlem}
Suppose that $x \in C^\infty(\mathbb{R},M)$ is a solution of (\ref{para}) and assume that $s_1<s_2$. Then it holds that
\begin{equation}\label{downhill}
f(x(s_2)) \leq f(x(s_1))
\end{equation}
and equality holds if and only if $x$ is constant. In this case $x$ is a critical point of $f$. 
\end{lemma}
\textbf{Proof: } We differentiate and take advantage of the gradient flow equation (\ref{para}) as well as of the definition of the gradient
\begin{eqnarray*}
\frac{d}{ds}f(x(s))&=&df(x(s))\partial_s x(s)\\
&=&-df(x(s)) \nabla f(x(s))\\
&=&-|| \nabla f (x(s))||^2_{g_{x(s)}}\\
&\leq & 0.
\end{eqnarray*}
Here $||\cdot||_g$ denotes the norm induced from the metric $g$. Integrating the above inequality we immediately obtain (\ref{downhill}). Moreover, we see that equality can only hold if $\nabla f(x(s))=0$ for every $s \in [s_1,s_2]$ and therefore $x(s)\in \mathrm{crit}f$. In view of the uniqueness of the initial value problem for the ODE (\ref{para}) we
conclude that in this case $x(s)$ has to be constant for all $s \in \mathbb{R}$. This proves the lemma. \hfill $\square$

\subsection{Energy}

Suppose that $x\colon \mathbb{R}$ is a smooth map not necessarily a gradient flow line. Then we define its energy with respect to the metric $g$ as
$$E(x):=E_g(x):=\int_{-\infty}^\infty ||\partial_s x||^2_g ds \in [0,\infty].$$
For gradient flow lines the following Proposition holds.
\begin{prop}\label{enerprop}
Assume that $x \in C^\infty(\mathbb{R},M)$ is a gradient flow line. Then
$$E(x)=\lim_{s \to -\infty} f(x(s))-\lim_{s \to \infty} f(x(s))=\sup_{s \in \mathbb{R}}f(x(s))-\inf_{s \in \mathbb{R}}f(x(s)).$$
\end{prop}
\textbf{Proof:} Using the gradient flow equation (\ref{para}) we compute
\begin{eqnarray*}
E(x)&=&\int_{-\infty}^\infty ||\partial_s x||^2 ds\\
&=&-\int_{-\infty}^\infty\langle \nabla f(x),\partial_s x\rangle ds\\
&=&-\int_{-\infty}^\infty df(x(s))\partial_s x ds\\
&=&-\int_{-\infty}^\infty \frac{d}{ds}f(x(s))ds\\
&=&\lim_{s \to -\infty} f(x(s))-\lim_{s \to \infty}f(x(s)).
\end{eqnarray*}
The second equation follows from the fact that gradient flow lines flow downhill as explained in Lemma~\ref{downlem}. \hfill$\square$
\\ \\
As an immediate Corollary of the Proposition we obtain
\begin{cor}\label{enbound}
Assume that $x \in C^\infty(\mathbb{R},M)$ is a gradient flow line whose asymptotics satisfy $\lim_{s \to \pm \infty} x(s)=x^\pm \in \mathrm{crit}f$. Then its energy is given by
$$E(x)=f(x^-)-f(x^+).$$
\end{cor}
From the corollary we learn that the energy of a gradient flow line only depends on its asymptotics. 

\subsection{Local convergence}

In this subsection we prove that maybe after transition to a subsequence a sequence of gradient flow lines locally converges to another gradient flow line. The convergence is not global. Even if all gradient flow lines in the sequence have the same asymptotics the limit gradient flow lines can have different asymptotics. For local convergence it is not important that the smooth function $f$ is Morse. What we however take advantage of is that the manifold $M$ is compact. 
\begin{prop}\label{loccom}
Assume that $x_\nu \in C^\infty(\mathbb{R},M)$ for $\nu \in \mathbb{N}$ is a sequence of gradient flow lines. Then there
exists a subsequence $\nu_j$ and a gradient flow line $x$ such that $x_{\nu_j}$ converges in the $C^\infty_{\mathrm{loc}}$-topology to $x$, i.e., for every real number $R>0$ the restriction $x_{\nu_j}|_{[-R,R]}$ converges in the $C^\infty$-topology to the restriction $x|_{[-R,R]}$.
\end{prop}
\textbf{Proof:} The proof relies on the Theorem of Arzel\`a-Ascoli. We carry out the proof in three steps.
\\ \\
\textbf{Step\,1:} \emph{The sequence $x_\nu$ is equicontinuous.}
\\ \\
Because $M$ is compact there exists a finite constant $c$ such that
$$||\nabla f(x)|| \leq c, \quad x  \in M.$$
Hence for $s_1<s_2$ we estimate using  the gradient flow equation (\ref{para})
$$d\big(x_\nu(s_1),x_\nu(s_2)\big)\leq \int_{s_1}^{s_2}||\partial_s x(s)||ds=\int_{s_1}^{s_2}||\nabla f(x(s))||ds \leq
c(s_2-s_1).$$
In particular, the right hand side is independent of $\nu$. This proves equicontinuity. 
\\ \\
\textbf{Step\,2:} \emph{There exists a subsequence $\nu_j$ and a continuous map $x \colon \mathbb{R} \to M$ such that
$x_{\nu_j}$ converges in the $C^0_{\mathrm{loc}}$-topology to $x$.}
\\ \\
By Step\,1 the sequence $x_\nu$ is equicontinuous. Moreover, the target space $M$ is compact. Therefore the assertion of Step\,2 follows from the Theorem of Arzel\`a-Ascoli.
\\ \\
\textbf{Step\,3:} \emph{Bootstrapping.}
\\ \\
The convergence we obtained in Step\,2 is just $C^0_{\mathrm{loc}}$. Using the gradient flow equation inductively we improve
the convergence to $C^\infty_{\mathrm{loc}}$. First note that by (\ref{para}) we have
\begin{equation}\label{gradeq}
\partial_s x_{\nu_j}=-\nabla f(x_{\nu_j})
\end{equation}
and by Step\,2 the right hand side converges in the $C^0_{\mathrm{loc}}$-topology to $-\nabla f(x)$. This implies that we can assume that $x \in C^1(\mathbb{R},M$ and
$$\partial_s x=-\nabla f(x).$$
Assume that $s_0 \in \mathbb{R}$. Choose local coordinates around $x(s_0)$. Since $x_{\nu_j}$ converges in the 
$C^0_{\mathrm{loc}}$-topology to $x$, there exists $\epsilon>0$ and $j_0 \in \mathbb{N}$ such that 
$x_{\nu_j}(s)$ for $s \in (s_0-\epsilon,s_0+\epsilon)$ and $j \geq j_0$ are contained in these local coordinates. Therefore
we can assume for the following local discussion without loss of generality that $M$ is an open subset of $\mathbb{R}^n$.
In particular, we can talk about higher partial derivatives $\partial^k_s x_{\nu_j}$. Applying the chain rule to the equation 
(\ref{gradeq}) inductively we conclude that
\begin{equation}\label{highder}
\partial_s^k x_{\nu_j}=F_k(x_{\nu_j},\partial_s x_{\nu_j},\ldots,\partial_s^{k-1} x_{\nu_j})
\end{equation}
for continuous functions $F_k$ involving $f$ and its derivatives. Hence if we assume by induction that
$x_{\nu_j}$ converges to $x$ in the $C^{k-1}_{\mathrm{loc}}$-topology it follows from (\ref{highder}) that 
$x_{\nu_j}$ converges to $x$ in the $C^k_{\mathrm{loc}}$-topology. This proves that $x_{\nu_j}$ converges to $x$
in the $C^\infty_{\mathrm{loc}}$-topology. \hfill $\square$

\subsection{Broken gradient flow lines}

Even if all gradient flow lines in a sequence have the same asymptotics, the limit gradient flow line one obtains by the local convergence theorem does not need to converge asymptotically to the same critical points. Even worse, changing the
parametrization of the gradient flow lines might result in a different limit gradient flow line. To deal with this situation we consider all nonconstant limit gradient flow lines one can obtain by a change of parametrization. Such a limit then makes as well sense as a limit for unparametrized gradient flow lines. Here is the relevant definition.

\begin{fed}
Assume that $x^\pm \in \mathrm{crit}f$. A \emph{(parametrised) broken gradient flow line} from $x^-$ to $x^+$ is a tupel
$$y=\{x^k\}_{1 \leq k \leq n},\quad n \in \mathbb{N}$$
such that the following conditions are met.
\begin{description}
 \item[(i)] For every $k \in \{1,\ldots,k\}$ it holds that $x^k$ is a nonconstant gradient flow line.
 \item[(ii)] The asymptotics of the gradient flow lines satisfy
\begin{eqnarray*}
\lim_{s \to -\infty} x^1(s)&=&x^-,\\
\lim_{s\to \infty} x^k(s)&=&\lim_{s\to -\infty}x^{k+1}(s),\quad 1 \leq k \leq n-1,\\
 \lim_{s \to \infty} x^n(s)&=&x^+.
\end{eqnarray*} 
\end{description}
\end{fed}

\begin{figure}[h]
\begin{center}
\includegraphics[scale=1]{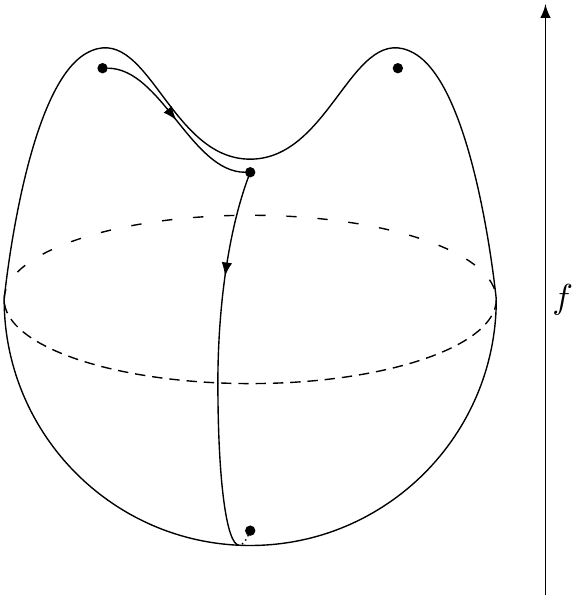}
\caption{Broken gradient flow line on deformed sphere.}
\end{center}
\end{figure}

\subsection{Floer-Gromov convergence}

On solutions of the gradient flow equation (\ref{para}) the noncompact group $\mathbb{R}$ acts by timeshift. It is a common phenomenon that if a noncompact group acts on a moduli space of solutions of a PDE, then the moduli space of unparametrized solutions has to be compactified with objects showing borderline analytic behaviour like breaking in the case of gradient flow lines, or bubbling in the case of holomorphic curves in Gromov-Witten theory. For example on holomorphic spheres the noncompact group $PSL(2,\mathbb{C})$ acts by reparametrisation of the domain via M\"obius transformations. In this case the combinatorics of the limiting curve gets even much more involved. While in the case of breaking we end up just with a string of gradient flow lines in the case of bubbling we can get a whole bubble tree, see for example \cite{mcduff-salamon}. For gradient flow lines the global notion of convergence is the following.
\begin{figure}[h]
\begin{center}
\includegraphics[scale=1]{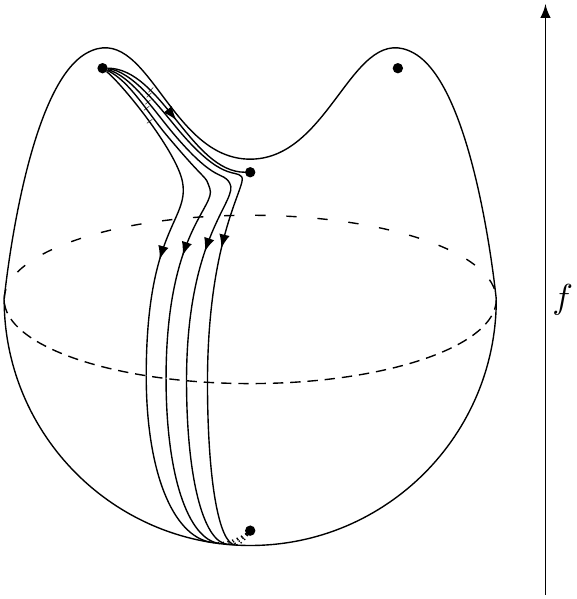}\quad
\includegraphics[scale=1]{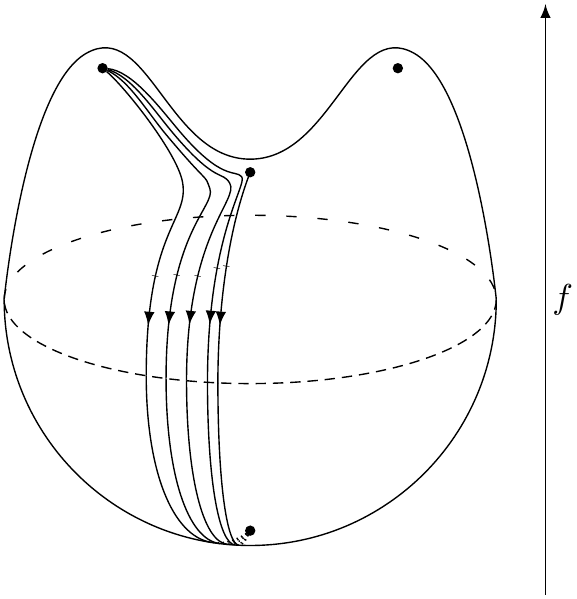}
\medskip

\includegraphics[scale=1]{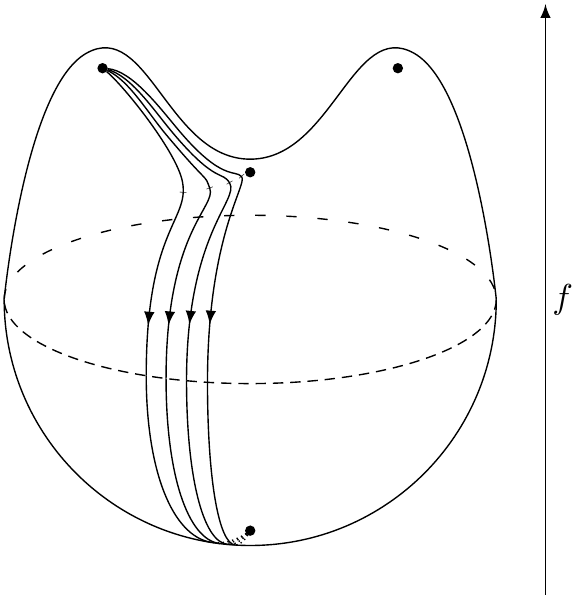}
\caption{Parametrised gradient flow lines converging to two components of a broken gradient flow line (above) and to a single critical point (below).}
\end{center}
\end{figure}
\begin{fed}
Suppose that $x_\nu \in C^\infty(\mathbb{R},M)$ is a sequence of gradient flow lines having fixed asymptotics
$\lim_{s \to \pm \infty} x_\nu(s)=x^\pm \in \mathrm{crit}f$ and $y=\{x^k\}_{1 \leq k \leq n}$ is a broken gradient flow line
from $x^-$ to $x^+$. We say that the sequence $x_\nu$ \emph{Floer-Gromov converges} to $y$ if the following holds true. Assume that for $1 \leq k \leq n$ there exists a sequence $r_\nu^k \in \mathbb{R}$ such that
$$(r_\nu^k)_*x_\nu \xrightarrow
 {
  \substack{C^\infty_{\mathrm{loc}}}
 }
 x^k.
$$
\end{fed}

\subsection{Floer-Gromov compactness}

Different from local compactness for the global Floer-Gromov compactness we assume now that $f$ is a Morse function. Moreover, the manifold $M$ is closed as usual. The main result of this section is the following.
\begin{thm}\label{floergromov}
Assume that $x_\nu$ is a sequence of gradient flow lines with fixed asymptotics $\lim_{s \to \pm \infty}x_\nu(s)=x^\pm
\in \mathrm{crit}f$ satisfying $x^-\neq x^+$. Then there exists a subsequence $\nu_j$ and a broken gradient flow line
$y=\{x^k\}_{1 \leq k \leq n}$ from $x^-$ to $x^+$ such that
$$x_{\nu_j} \xrightarrow
 {
  \substack{\mathrm{Floer-Gromov}}
 }y.$$
\end{thm}
Note that the assumption that the asymptotics $x^+$ and $x^-$ are different rules out the trivial case that the sequence just consists of constant gradient flow lines each of them the same critical point of the Morse function $f$.
In order to prove the Theorem we need two Lemmas. The first of them is the following.
\begin{lemma}\label{finite}
The Morse function $f$ has only finitely many critical points on the closed manifold $M$. 
\end{lemma}
\textbf{Proof: } Because $M$ is compact it suffices to show that the critical points of $f$ are isolated. This follows immediately from the Morse Lemma, namely Lemma~\ref{morselemma}. More elementary, not relying on the Morse Lemma, we can argue as follows. We consider the smooth function $e \colon M \to \mathbb{R}$ given for $x \in M$ by
$$e(x)=||\nabla f(x)||^2=df(x) \nabla f(x).$$
Suppose that $x_0$ is a critical point of $f$. We examine the Taylorexpansion of $e$ in local coordinates $U$ around $x_0$. 
We can assume without loss of generality that $x_0=0$ in the local coordinates chosen. We introduce the metric Hessian. If $x \in U$, then the metric Hessian is a linear map
$$H_f^g \colon T_x U \to T_x U$$
where the tangent space $T_x U$ in local coordinates can of course be canonically identified with $\mathbb{R}^n$. For $v,w \in T_x U$ it is defined implicitly by the equation 
$$d^2 f(x)(v,w)=g_x(v,H_f^g(x)w).$$
The metric Hessian is selfadjoint with respect to the inner product $g_x$ on $T_x U$. While at a critical point the usual Hessian transforms as a bilinear map the metric Hessian transforms as a linear map under coordinate change. Different from the usual Hessian it depends on the choice of the Riemannian metric. Using the metric Hessian we can write the differential of $e$ as follows
\begin{eqnarray}\label{de}
de(x)v&=&2d^2f(x)\big(\nabla f(x),v\big)\\ \nonumber
&=&2g_x\big(\nabla f(x),H^g_f(x)v\big)\\ \nonumber
&=&2 df(x) H^g_f(x)v.
\end{eqnarray}
In particular, 
\begin{equation}\label{de0}
de(0)=0
\end{equation}
since $df(0)=0$. Given two tangent vectors $v,w \in T_0 U$ we obtain from (\ref{de}) for the Hessian of $e$ at the critical point $0$ 
\begin{eqnarray}\label{d2e}
d^2 e(0)(v,w)&=&2 d^2 f(0)\big(H^g_f(0)v,w\big)\\ \nonumber
&=&2g_0\big(H^g_f(0)v,H^g_f(0)w\big).
\end{eqnarray}
Using (\ref{de0}), (\ref{d2e}) together with $e(0)=0$ we see by the Theorem of Taylor that there exists a finite constant $c$ such that
\begin{equation}\label{taylor}
\Big| e(x)-g_0\big(H^g_f(0)x,H^g_f(0)x\big)\Big| \leq c||x||^3.
\end{equation}
Because $f$ is Morse $H^g_f(0)$ is nondegenerate and therefore there exists a positive constant $c'$ such that
\begin{equation}\label{he}
g_0\big(H^g_f(0)x,H^g_f(0)x\big) \geq c'||x||^2.
\end{equation}
Combining inequalities (\ref{taylor}) and (\ref{he}) we obtain the inequality
$$e(x) \geq c'||x||^2-c||x||^3.$$
Note that the function $r \mapsto c'r^2-cr^3$ is positive for $r \in \big(0,\tfrac{c'}{c}\big)$ we conclude that
$0$ is an isolated zero of $e$ and therefore an isolated critical point of $f$. This finishes the proof of the Lemma. \hfill $\square$
\\ \\
Our next Lemma tells us that gradient flow lines asymptotically converge. This Lemma as well is not true without the assumption that $f$ is Morse. For functions which are not Morse it might happen that a gradient flow line spirals into
a set of critical points. An example of such a phenomenon can be found in \cite[page 13-14]{palis-demelo}
\begin{lemma}\label{asymp}
Suppose that $x$ is a gradient flow line. Then there exist $x^\pm \in \mathrm{crit}f$ such that
$$\lim_{s \to \pm \infty} x(s)=x^\pm.$$
\end{lemma}
\textbf{Proof: } By Lemma~\ref{finite} we know that the cardinality of critical points of $f$ is finite. Hence for some
$N \in \mathbb{N}$ we can write
$$\mathrm{crit}f=\{x_1,\ldots,x_N\}.$$
Choose open neighbourhoods $U_i$ of $x_i$ for $i \in \{1,\ldots,N\}$ which are pairwise disjoint, i.e.,
$$U_i \cap U_j=\emptyset, \quad i \neq j.$$
We first show the following claim which tells us that the gradient of $f$ is uniformly bounded away from zero on the complement of these open neighbourhoods of the critical points.
\\ \\
\textbf{Claim\,1: } \emph{There exists $\epsilon_0>0$ such that} 
$$||\nabla f(x)|| \geq \epsilon_0,\quad \forall\,\,x \in M \setminus \bigcup_{i=1}^N U_i.$$
\textbf{Proof of Claim\,1: } Since $M \subset \bigcup_{i=1}^N U_i$ is a closed subset of $M$ and $M$ is compact we conclude that $M \setminus \bigcup_{i=1}^N U_i$ is compact. We assume by contradiction that there exists a sequence of points
$y_n \in M \setminus \bigcup_{i=1}^N U_i$ for $n \in \mathbb{N}$ satisfying 
\begin{equation}\label{klein}
||\nabla f(y_n)||\leq \frac{1}{n}.
\end{equation}
Since $M\setminus \bigcup_{i=1}^NU_i$ is compact we conclude that $y_n$ has a convergent subsequence $y_{n_j}$ such that
$$y:=\lim_{j \to \infty}y_{n_j} \in M \setminus \bigcup_{i=1}^N U_i.$$
From (\ref{klein}) we deduce that $||\nabla f(y)||=0$. This implies that $y$ is a critical point of $f$ which contradicts
the fact that $y$ lies in $M \setminus \bigcup_{i=1}^N U_i$, where by construction no critical points of $f$ are. 
This contradiction establishes Claim\,1. \hfill $\square$
\\ \\
We now choose $\epsilon_0$ as in Claim\,1 and for $0<\epsilon \leq \epsilon_0$ we introduce the following open neighbourhood
of the critical point $x_i$
$$V_i^\epsilon:=\big\{x \in U_1:||\nabla f(x)||<\epsilon\big\} \subset U_i.$$
Our next claim is the following.
\\ \\
\textbf{Claim\,2: } \emph{For every $\epsilon>0$ there exists a sequence $s_\nu^\epsilon \to \infty$ such that}
$$x(s_\nu^\epsilon) \in \bigcup_{i=1}^N V_i^\epsilon.$$
\textbf{Proof of Claim\,2: } We argue by contradiction. If Claim\,2 fails to be true then there exists $\sigma \in \mathbb{R}$ with the property that
$$x(s) \notin \bigcup_{i=1}^N V_i^\epsilon, \quad \forall\,\,s \geq \sigma.$$
By construction of the sets $V_i^\epsilon$ this implies that
$$||\nabla f(x(s))|| \geq \epsilon,\quad \forall\,\,s \geq \sigma.$$
Using Proposition~\ref{enerprop} and the gradient flow equation (\ref{para}) we conclude from that
\begin{eqnarray*}
\max f-\min f&\geq& E(x)\\
&=&\int_{-\infty}^\infty||\partial_s x||^2 ds\\
&=&\int_{-\infty}^\infty ||\nabla f(x)||^2 ds\\
&\geq & \int_{\sigma}^\infty ||\nabla f(x)||^2ds\\
&\geq & \int_{\sigma}^\infty \epsilon^2 ds\\
&=& \infty.
\end{eqnarray*}
This contradiction shows that Claim\,2 is true. \hfill $\square$
\\ \\
Note that $\overline{V_i^{\epsilon/2}}$, the closure of the set $V_i^{\epsilon/2}$, and $M \setminus V_i^\epsilon$ are compact subsets of $M$, which are disjoint, i.e.,
$$\overline{V_i^{\epsilon/2}}\cap \big(M\setminus V_i^\epsilon\big)=\emptyset.$$
In particular, its distance satisfies
$$d\Big(\overline{V_i^{\epsilon/2}},M\setminus V_i^\epsilon\Big)>0.$$
We abbreviate
$$\kappa_\epsilon:=\min_{i \in \{1,\ldots,N\}}d\Big(\overline{V_i^{\epsilon/2}},M\setminus V_i^\epsilon\Big)>0.$$
Because
$$E=\int_{-\infty}^\infty||\partial_s x||^2 ds <\infty$$
there exists $\sigma_\epsilon \in \mathbb{R}$ satisfying
$$\int_{\sigma_\epsilon}^\infty||\partial_s x||^2ds \leq \frac{\epsilon \kappa_\epsilon}{4}.$$
Our third claim is the following.
\\ \\
\textbf{Claim\,3: } \emph{There exists $i \in \{1,\ldots,N\}$ such that $x(s) \in V_i^\epsilon$ for $s \geq \sigma_\epsilon$.}
\\ \\
\textbf{Proof of Claim\,3: } By Claim\,2 there exists $s_0 \geq \sigma_\epsilon$  and $i \in \{1,\ldots,N\}$ such that
$$x(s_0) \in V^{\epsilon/2}_i.$$
We assume by contradiction that there exists $s_1 \geq \sigma_\epsilon$ with the property that
$$x(s_1) \notin V_i^\epsilon.$$
This implies that there exist times $\sigma_\epsilon \leq t_0<t_1$ such that
$$x(t) \in V_i^\epsilon \setminus V_i^{\epsilon/2}, \quad t \in [t_0,t_1]$$
and either
$$x(t_0) \in \partial V_i^\epsilon, \quad x(t_1) \in \partial V_i^{\epsilon/2}$$
or
$$x(t_1) \in \partial V_i^\epsilon, \quad x(t_0) \in \partial V_i^{\epsilon/2}.$$
Note that for $t \in [t_0,t_1]$ we have
$$||\partial_s x(t)||=||\nabla f(x(t))|| \geq \frac{\epsilon}{2}.$$
Hence we estimate
\begin{eqnarray*}
\kappa_\epsilon &\leq& d\big(x(t_0),x(t_1)\big)\\
&\leq&\int_{t_0}^{t_1} ||\partial_s x||ds\\
&\leq&\frac{2}{\epsilon}\int_{t_0}^{t_1}||\partial_s x||^2 ds\\
&\leq&\frac{2}{\epsilon}\int_{\sigma_\epsilon}^{\infty}||\partial_s x||^2 ds\\
&\leq&\frac{2}{\epsilon}\frac{\epsilon \kappa_\epsilon}{4}\\
&=&\frac{\kappa_\epsilon}{2}.
\end{eqnarray*}
This contradiction shows that Claim\,3 has to be true. \hfill $\square$
\\ \\
Using Claim\,3 it is now straightforward to prove the Lemma. Assume that $V$ is an arbitrary open neighbourhood of $x_i$.
Because $||\nabla f||$ is a continuous function on $M$ there exists $\epsilon>0$ such that $V_i^\epsilon \subset V$.
Hence by Claim\,3 it holds that
$$x(s) \in V, \quad \forall\,\,s \geq \sigma_\epsilon.$$
Since $V$ was an arbitrary neighbourhood of $x_i$ we deduce that
$$\lim_{s \to \infty} x(s)=x_i.$$
A similar argument shows that $\lim_{s \to -\infty}x(s)$ exists as well. This finishes the proof of the Lemma. \hfill $\square$
\\ \\
We are now ready to embark on the proof of the main result of this section. 
\\ \\
\textbf{Proof of Theorem~\ref{floergromov}: } We prove the theorem by induction. For $m \in \mathbb{N}$ we establish the following assertion
\begin{description}
 \item[($A_m$)] There exists a subsequence $\nu_j$, a broken gradient flow line $y_m=\{x^k\}_{1 \leq k \leq \ell}$ for
 $\ell \leq m$, and sequences $r^k_j \in \mathbb{R}$ for $1 \leq k \leq \ell$ meeting the following requirements
 \begin{description}
 \item[(i)]$(r_j^k)_*x_{\nu_j} \xrightarrow
 {
  \substack{C^\infty_{\mathrm{loc}}}
 }
 x^k, \quad 1 \leq k \leq \ell,$
 \item[(ii)] $\lim_{s \to -\infty}x^1(s)=x^-$,
 \item[(iii)] If $\ell<m$ it holds that $\lim_{s \to \infty}x^\ell(s)=x^+$.
 \end{description}
\end{description}
We first do the base case and show that assertion $(A_1)$ is true. For that purpose we choose an open neighbourhood $V$
of $x^-$. In the proof of Lemma~\ref{finite} we learned that critical points of a Morse function are isolated. Therefore
we can choose $V$ having the additional property that
$$\overline{V} \cap \mathrm{crit}f=\{x^-\}$$
where $\overline{V}$ is the closure of the open set $V$. We define
$$r_\nu^1:=\inf\big\{s \in \mathbb{R}: x_\nu(s) \notin V\big\}$$
\begin{figure}[h]
\begin{center}
 \includegraphics[scale=1]{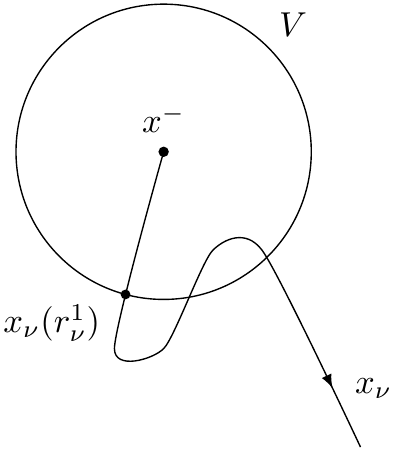}
 \caption{First exit time of gradient flow line out of $V$.}
\end{center}
\end{figure}
the first exit time of $V$. Because $\lim_{s \to \infty}x_\nu(s)=x^+ \notin V$ we conclude that
$r^1_\nu$ is finite. Using Proposition~\ref{loccom} we conclude that there exists a subsequence $\nu_j$ and a gradient
flow line $x^1$ such that
\begin{equation}\label{conv}
(r_{\nu_j}^1)_*x_{\nu_j} \xrightarrow
 {
  \substack{C^\infty_{\mathrm{loc}}}
 }
 x^1.
\end{equation}
The gradient flow line $x^1$ is not constant since
$$x^1(0)=\lim_{j \to \infty}(r_{\nu_j}^1)_*x_{\nu_j}(0)=\lim_{j \to \infty}x_{\nu_j}(r_{\nu_j}^1) \in \partial V$$
and by the choice of $V$
$$\partial V \cap \mathrm{crit}f=\emptyset.$$
Define
$$y^1:=\{x^1\}$$
and
$$r^1_j:=r^1_{\nu_j}.$$
We need to check that with these choices the three requirements of $(A_1)$ are satisfied. Assertion (i) is just (\ref{conv}). Assertion (ii) requires some argument. By Lemma~\ref{asymp} we know that 
$\lim_{s\to -\infty}x^1(s)$ exists in $\mathrm{crit}f$. The time $r^1_\nu$ is defined as the first exit time, therefore
$$(r^1_{\nu_j})_* x_{\nu_j}(s) \in \overline{V}, \quad s \leq 0.$$
Hence
$$x^1(s) \in \overline{V}, \quad s \leq 0.$$
Therefore
$$\lim_{s \to -\infty}x^1(s) \in \mathrm{crit}f \cap \overline{V}=\{x^-\}.$$
This proves assertion (ii) and assertion (iii) is empty. We have established the truth of $(A_1)$.
\\ \\
We next carry out the induction step. We assume that $(A_m)$ is true and show that this implies that $(A_{m+1})$ has to be true as well. Let $y_m=\{x^k\}_{1\leq k \leq \ell}$ be the broken gradient flow line provided by $(A_m)$. We distinguish two cases. The first trivial case is the following.
\\ \\
\textbf{Case\,1:} $\lim_{s \to \infty}x^\ell(s)=x^+$.
\\ \\
In this case we set $y_{m+1}=y_m.$
Then $y_{m+1}$ satisfies the requirements of $(A_{m+1})$. 
\\ \\
The interesting nontrivial case is the following.
\\ \\
\textbf{Case\,2:} $\lim_{s \to \infty}x^\ell(s)\neq x^+$.
\\ \\
In this case by (iii) we necessarily have $\ell=m$. By Lemma~\ref{asymp} we know that
$$(x^m)^+:=\lim_{s\to \infty}x^m(s) \in \mathrm{crit}f$$
exists. Choose an open neighbourhood $W$ of $(x^m)^+$ such that
$$\overline{W} \cap \mathrm{crit}f=\big\{(x^m)^+\big\}.$$
Using that $\lim_{s\to \infty}x^m(s)=(x^m)^+$ we conclude that there exists $s_0 \in \mathbb{R}$ such that
$$x^m(s) \in W, \quad s \geq s_0.$$
Taking further advantage of the fact that
$$
(r_j^m)_*x_{\nu_j} \xrightarrow
 {
  \substack{C^\infty_{\mathrm{loc}}}
 }
 x^m$$
we infer that there exists $j_0 \in \mathbb{N}$ such that for every $j \geq j_0$ it holds that
$$(r_j^m)_*x_{\nu_j}(s_0) \in W.$$
For $j \geq j_0$ we define
$$R_j:=\inf\big\{r \geq 0: (r_j^m)_* x_{\nu_j}(s_0+r) \notin W\big\}$$
\begin{figure}[h]
\begin{center}
 \includegraphics[scale=1]{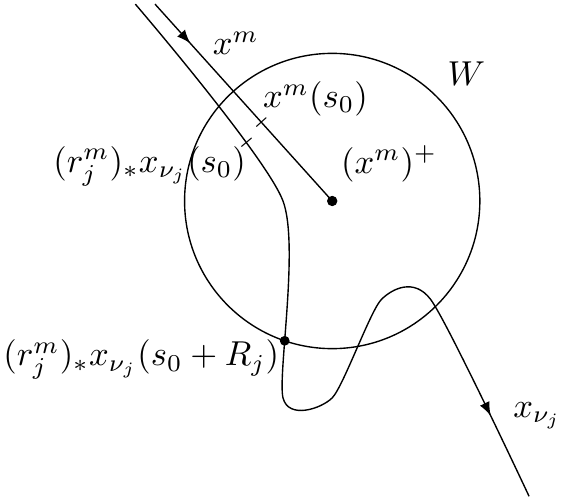}
 \caption{First exit time of gradient flow line out of $W$.}
\end{center}
\end{figure}
the first exit time of $(r_j^m)_*x$ from $W$ after $s_0$. Since
$$\lim_{s \to \infty}(r_j^m)_* x_{\nu_j}(s)=x^+ \notin \overline{W}$$
we conclude that
$$R_j<\infty, \quad j \geq j_0.$$
On the other hand because $x^m(s) \in W$ for $s \geq s_0$ and $(r_j^m)_*x_{\nu_j}$ converges in the 
$C^\infty_{\mathrm{loc}}$-topology to $x^m$ we deduce that
\begin{equation}\label{diver}
\lim_{j \to \infty} R_j=\infty.
\end{equation}
We define
$$r_j^{m+1}:=r_j^m+s_0+R_j.$$
After maybe transition to a further subsequence we can assume by Proposition~\ref{loccom} that there exists a
gradient flow line $x^{m+1}$ such that 
$$
(r_j^{m+1})_*x_{\nu_j} \xrightarrow
 {
  \substack{C^\infty_{\mathrm{loc}}}
 }
 x^{m+1}.$$
We first claim that $x^{m+1}$ is not constant. To see that we observe 
$$(r_j^{m+1})_*x_{\nu_j}(0)=x_{\nu_j}(r_j^{m+1})=x_{\nu_j}(r_j^m+s_0+R_j)=(r_j^m)_*x_{\nu_j}(s_0+R_j) \in \partial W.$$
Therefore
$$x^{m+1}(0)\in \partial W.$$
Because $\partial W \cap \mathrm{crit}f=\emptyset$ we conclude that $x^{m+1} \notin \mathrm{crit}f$
and therefore $x^{m+1}$ is not constant. We next show that
\begin{equation}\label{limi}
\lim_{s \to -\infty} x^{m+1}(s)=(x^m)^+.
\end{equation}
To understand that first recall that by Lemma~\ref{asymp} 
$\lim_{s \to -\infty}x^{m+1}(s) \in \mathrm{crit}f$ exists. We assume by contradiction that there exists a critical point $x'$ of $f$ such that
$$x'=\lim_{s \to -\infty}x^{m+1}(s) \neq (x^m)^+.$$
In this case we can choose an open neighbourhood $U$ of $x'$ such that
\begin{equation}\label{disjoint}
U \cap W=\emptyset.
\end{equation}
Because $\lim_{s\to-\infty} x^{m+1}(s)=x'$ there exists $s_1<0$ such that
$$x^{m+1}(s) \in U, \quad s \leq s_1.$$
Since $(r_j^{m+1})_* x_{\nu_j}$ converges to $x^{m+1}$ in the $C^\infty_{\mathrm{loc}}$-topology, we conclude  that there
exists a positive integer $j_0'$ such that for every $j \geq j_0'$ we have
\begin{equation}\label{inu}
(r_j^{m+1})_* x_{\nu_j}(s_1) \in U.
\end{equation}
On the other hand we know for $j \geq j_0$ that 
\begin{equation}\label{inw}
(r_j^m)_*x_{\nu_j}(s) \in W, \quad s \in (s_0,s_0+R_j).
\end{equation}
Shifting time back and forth we obtain
$$(r_j^m)_* x_{\nu_j}(s)=x_{\nu_j}(s+r_j^m)=x_{\nu_j}(s+r_j^{m+1}-s_0-R_j)=(r_j^{m+1})_*x_{\nu_j}(x-s_0-R_j).$$
Using (\ref{inw}) we infer from this formula that
\begin{equation}\label{inw2}
(r_j^{m+1})_* x_{\nu_j}(s) \in W, \quad s \in (-R_j,0).
\end{equation}
Since the sequence $R_j$ diverges as we noted in (\ref{diver}) above there exists an integer $j_1 \geq \max\{j_0,j_0'\}$ such that 
$$-R_{j_1} \leq s_1.$$
Because $s_1$ is chosen negative we obtain by combining (\ref{inu}) and (\ref{inw2}) that
$$(r_{j_1}^{m+1})_* x_{\nu_j}(s_1) \in U \cap W$$
which contradicts (\ref{disjoint}), which tells us that the sets $U$ and $W$ are disjoint. This contradiction proves that (\ref{limi}) is correct.  
\\ \\
We now set 
$$y_ {m+1}=\{x^k\}_{1 \leq k  \leq m+1}.$$
With this choice $y_{m+1}$ is a broken gradient flow line satisfying assertion $(A_{m+1})$.
This finishes the proof of the induction step.
\\ \\
To finish the proof of the Theorem we need to show that the induction terminates in the sense that there exists $m  \in  \mathbb{N}$ with the property that the positive asymptotics of the last gradient flow line $(x^m)^+=\lim_{s \to \infty} x^m(s)$ satisfies 
$$(x^m)^+=x^+.$$
To see that note that due to the fact that nonconstant gradient flow lines flow downhill
by Lemma~\ref{downlem} we have that for every $j \in \{1,\ldots,m-1\}$ it holds that
$$f((x^j)^+)>f((x^{j+1})^+).$$
This shows that all positive asymptotics are different. Since by Lemma~\ref{finite}
there are only finitely many critical points we conclude that the induction has to terminate and the proof of Theorem~\ref{floergromov} is complete. \hfill $\square$

\subsection{Uniqueness of limits}

We call two broken gradient flow lines $y_0=\{x_0^k\}_{1 \leq k \leq m_0}$ and $y_1=\{x_1^k\}_{1 \leq k \leq m_1}$
\emph{equivalent}, when $m_0=m=1=:m$ and for every $k \in \{1,\ldots,m\}$ there exists $r^k \in \mathbb{R}$ such that
$$x^k_0=r^k_* x^k_1.$$
We refer to an equivalence class $[y]$ of broken gradient flow lines as an \emph{unparametrized broken gradient flow line}.
We prove that the Floer-Gromov limit of a sequence of gradient flow lines is unique up to equivalence. 
\begin{thm}\label{unilim}
Suppose that $x_\nu$ is a sequence of gradient flow lines with fixed positive and negative asymptotics $x^\pm$ which converges in the sense of Floer-Gromov to a broken gradient flow line $y_0$ from $x^-$ to $x^+$ as well as to a broken gradient flow line $y_1$ from $x^-$ to $x^+$. Then $y_0$ and $y_1$ are equivalent. 
\end{thm} 
To prove the theorem we need the following Proposition.
\begin{prop}\label{repro}
Suppose that $x_\nu$ is a sequence of gradient flow lines which converges in the $C^\infty_{\mathrm{loc}}$-topology
to a nonconstant gradient flow line $x$ and $r_\nu$ is a sequence of real numbers such that the sequence of reparametrized
gradient flow lines $(r_\nu)_*x_\nu$ converges in the $C^\infty_{\mathrm{loc}}$-topology to a nonconstant gradient flow line
$\widetilde{x}$. Assume that the two gradient flow lines $x$ and $\widetilde{x}$ have at least one common asymptotic, i.e., $\lim_{s\to-\infty}
x(s)=\lim_{s\to -\infty}\widetilde{x}(s)$ or $\lim_{s\to\infty}
x(s)=\lim_{s\to \infty}\widetilde{x}(s)$. Then the sequence $r_\nu$ converges to some $r \in \mathbb{R}$ and $\widetilde{x}=r_*x$.
\end{prop}
\textbf{Proof: } We treat the case where the negative asymptotic agree, i.e., we assume that
$$\lim_{s\to-\infty}x(s)=\lim_{s\to -\infty}\widetilde{x}(s)=:x^-.$$
The case where the positive asymptotics agree can be treated similarly and follows as well from the case of the negative asymptotics by replacing the Morse function $f$ by the Morse function $-f$ and by observing that in this case the gradient flow lines flip their orientation such that the positive asymptotic becomes the negative one and vica versa. We prove the Proposition in four steps. 
\\ \\
\textbf{Step\,1: } \emph{The sequence $r_\nu$ is bounded from above.}
\\ \\
\textbf{Proof of Step\,1: } We argue by contradiction and assume that there exists a subsequence $\nu_j$ such that
\begin{equation}\label{divrho2}
\lim_{j \to \infty} r_{\nu_j}=\infty.
\end{equation}
We choose an open neighbourhood $V$ of $x^-$ with the property that
\begin{equation}\label{einzig}
\overline{V} \cap \mathrm{crit}f=\{x^-\}.
\end{equation}
We need the following Claim. 
\\ \\
\textbf{Claim: } \emph{There exists an open neighbourhood $W$ of $x^-$ contained in $V$ with the property that if
$x$ is a gradient flow line for which there exist times $s_0<s_1$ such that $x(s_0) \in W$ and $x(s_1) \notin V$
then for every $s \geq s_1$ it holds that $x(s) \notin W$.}
\\ \\
\textbf{Proof of Claim: }Since $\overline{V}$ is compact we can thanks to (\ref{einzig}) find an $\epsilon>0$ and
an open neigbourhood $W_1$ of $x^-$ satisfying
$$\overline{W}_1 \subset V$$
such that
$$||\nabla f(x)|| \geq \epsilon, \quad x \in \overline{V} \setminus W_1.$$
Abbreviate
$$\kappa:=d\big(\overline{W}_1,M\setminus V\big)>0$$
the distance between the compact disjoint sets $\overline{W}_1$ and $M \setminus V$. 
We set
$$W:=W_1 \cap f^{-1}\Big(f(x^-)-\tfrac{\epsilon \kappa}{2},f(x^-)+\tfrac{\epsilon \kappa}{2}\Big).$$
Observe that $W$ is an open neighbourhoof of $x^-$ which is contained in the open neighbourhood $W_1$.
In order to see that $W$ meets the requirements of the claim we first note that there exists times
$$s_0<t_0<t_1<s_1$$
satisfying
$$x(t_0) \in \partial W_1,\qquad x(t_1) \in \partial V, \qquad x(t) \in V \setminus W_1,\,\,t \in (t_0,t_1).$$
We estimate using the gradient flow equation (\ref{para})
\begin{eqnarray*}
f(x(s_0))-f(x(s_1))&=&\int_{s_0}^{s_1}\frac{d}{ds}f(x(s))ds\\
&=&\int_{s_0}^{s_1} df(x(s))\partial_s x(s) ds\\
&=&\int_{s_0}^{s_1}||\nabla f(x(s))||\cdot||\partial_s x(s)||ds\\
&\geq&\int_{t_0}^{t_1}||\nabla f(x(s))||\cdot||\partial_s x(s)||ds\\
&\geq&\epsilon \int_{t_0}^{t_1}||\partial_s x(s)||ds\\
&\geq& \epsilon \kappa.
\end{eqnarray*}
From this inequality we infer using that $x(s_0) \in W$
$$f(x(s_1)) \leq f(x(s_0))-\epsilon \kappa <f(x^-)+\frac{\epsilon \kappa}{2}-\epsilon \kappa=f(x^-)-\frac{\epsilon \kappa}{2}.$$
Using that gradient flow lines flow downhill by Lemma~\ref{downlem} we conclude that
$$f(x(s))<f(x^-)-\frac{\epsilon \kappa}{2}, \quad s\geq s_1.$$
By construction of $W$ this implies that
$$f(x(s)) \notin W, \quad s \geq s_1.$$
This establishes the truth of the Claim. 
\\ \\
We choose an open neighbourhood $W$ of $x^-$ contained in $V$ as in the Claim. Since the negative asymptotic of $x$ is 
$x^-$, there exists $s_0 \in \mathbb{R}$
such that
$$x(s_0) \in W.$$
Because $x_\nu$ converges in the $C^\infty_{\mathrm{loc}}$-topology to $x$ we conclude that
there exists $\nu_0 \in \mathbb{N}$ with the property that
\begin{equation}\label{schw1}
x_\nu(s_0) \in W, \quad \nu \geq \nu_0.
\end{equation}
Taking advantage of the fact that $x$ is a nonconstant gradient flow line, there exists
$$s_1>s_0$$
with the property that
$$x(s_1) \notin \overline{V}.$$
Again using that $x_\nu$ converges in the $C^\infty_{\mathrm{loc}}$-topology to $x$ we obtain
$\nu_1 \geq \nu_0$ such that 
\begin{equation}\label{schw2}
x_\nu(s_1) \notin \overline{V}, \quad \nu \geq \nu_1.
\end{equation}
Note that
\begin{equation}\label{schw3}
x_\nu(s_1)=(r_\nu)_*x_\nu(s_1-r_\nu).
\end{equation}
Applying (\ref{schw3}) to (\ref{schw1}) and (\ref{schw2}) we obtain
\begin{equation}\label{schw4}
(r_\nu)_*x_\nu(s_0-r_\nu) \in W,\qquad (r_\nu)_*x_\nu(s_1-r_\nu) \notin \overline{V}, \quad \nu \geq \nu_1.
\end{equation}
Since the negative asymptotic of $\widetilde{x}$ is $x^-$ as well, there exists $s_2 \in \mathbb{R}$ such that 
$$\widetilde{x}(s_2) \in W$$
from which we infer, using that $(r_\nu)_*x_\nu$ converges in the $C^\infty_{\mathrm{loc}}$-topology to
$\widetilde{x}$, that there exists $\nu_2 \geq \nu_1$ such that
\begin{equation}\label{schw5}
(r_\nu)_*x_\nu(s_2) \in W,\quad \nu \geq \nu_2.
\end{equation}
By (\ref{divrho2}) there exists $j_0 \in \mathbb{N}$ such that $\nu_{j_0} \geq \nu_2$ and 
\begin{equation}\label{schw6}
s_1-\rho^\ell_{\nu_{j_0}}<s_2.
\end{equation}
However (\ref{schw4}),(\ref{schw5}) combined with (\ref{schw6}) contradict the construction of $W$. This finishes the proof of Step\,1.
\\ \\
\textbf{Step\,2: } \emph{The sequence $r_\nu$ has a converging subsequence.}
\\ \\
By Step\,1 the sequence $r_\nu$ is bounded from above. Interchanging the roles of $x$ and $\widetilde{x}$ we see that
$-r_\nu$ as well is bounded from above. This implies that the sequence $r_\nu$ is bounded. Hence Step\,2 follows from the Theorem of Bolzano-Weierstrass. 
\\ \\
\textbf{Step\,3: } \emph{The sequence $r_\nu$ converges to a real number $r$.}
\\ \\
By Step\,2 there exists $r \in \mathbb{R}$ and a subsequence $\nu_j$ such that 
$$\lim_{j \to \infty}r_{\nu_j}=r.$$
We assume by contradiction that not the whole sequence $r_\nu$ converges to $r$. We have seen in the proof of
Step\,1b that the sequence $r_\nu$ is bounded. Hence in this case there exists
$$r' \neq r$$
and another subsequence $\nu'_j$ such that
$$\lim_{j \to \infty}r_{\nu'_j}=r'.$$
To derive a contradiction of this we first compute
\begin{eqnarray}\label{r1}
\widetilde{x}(0)&=&\lim_{j \to \infty}(r_{\nu_j})_*x_{\nu_j}(0)\\ \nonumber
&=&r_* x(0)\\ \nonumber
&=&x(r).
\end{eqnarray}
Redoing the same computation for the sequence $\nu'_j$ instead of $\nu_j$ leads to
$$\widetilde{x}(0)=x(r').$$
Combining these two equations we obtain
$$x(r)=x(r').$$
However, because $r \neq r'$ we obtain from Lemma~\ref{downlem} that $x$ is a constant gradient flow lines contradicting the assumption of the Proposition. This contradiction proves Step\,3. 
\\ \\
\textbf{Step\,4: } $\widetilde{x}=r_* x$.
\\ \\
The computation in (\ref{r1}) implies that
$$\widetilde{x}(0)=r_*x(0).$$
Because the gradient flow equation (\ref{para}) is an ODE and therefore a solution is uniquely determined by its initial
condition we deduce that
$$\widetilde{x}=r_* x.$$
This finishes the proof of Step\,4 and the proof of the Proposition is complete. \hfill $\square$
\\ \\
\textbf{Proof of Theorem~\ref{unilim}: } We prove the theorem by induction on the gradient flow lines appearing in the two broken ones. If
$y_0=\{x_0^k\}_{1 \leq k \leq m_0}$ and $y_1=\{x_1^k\}_{1 \leq k \leq m_1}$ we first show as the base case the following assertion.
\\ \\
\textbf{Step\,1: } \emph{There exists $r^1 \in \mathbb{R}$ such that $x^1_0=r^1_* x^1_1$.}
\\ \\
\textbf{Proof of Step\,1: } By definition of Floer-Gromov convergence there exist sequences $r^1_{0,\nu} \in \mathbb{R}$
and $r^1_{1,\nu} \in \mathbb{R}$ such that
$$
(r_{0,\nu}^1)_*x_\nu \xrightarrow
 {
  \substack{C^\infty_{\mathrm{loc}}}
 }
 x^1_0, \qquad (r_{1,\nu}^1)_*x_\nu \xrightarrow
 {
  \substack{C^\infty_{\mathrm{loc}}}
 }
 x^1_1.$$
Since $x^1_0$ and $x^1_1$ have the common negative asymptotic $x^-$ Proposition~\ref{repro} tells us that
$$r^1:=\lim_{\nu \to \infty}(r^1_{0,\nu}-r^1_{1,\nu})$$
exists and satisfies
$$x^1_0=r^1_*x^1_1.$$
This finishes the proof of the base case Step\,1.
\\ \\
Step\,2 is the induction step. We assume that for $\ell \leq \min\{m_0,m_1\}$ we have establishes the existence of real numbers $r^k$ for $k \in \{1,\ldots,\ell-1\}$ such that
$$x^k_0=r^k_* x^k_1, \quad 1 \leq k \leq \ell-1.$$ 
Under this induction hypothesis the induction step is now
\\ \\
\textbf{Step\,2: } \emph{There exists $r^\ell \in \mathbb{R}$ such that $x^\ell_0=r^\ell_* x^\ell_1$.}
\\ \\
\textbf{Proof of Step\,2: } Again by definition of Floer-Gromov convergence there exist sequences $r^\ell_{0,\nu} \in \mathbb{R}$
and $r^\ell_{1,\nu} \in \mathbb{R}$ such that
$$
(r_{0,\nu}^\ell)_*x_\nu \xrightarrow
 {
  \substack{C^\infty_{\mathrm{loc}}}
 }
 x^\ell_0, \qquad (r_{1,\nu}^1)_*x_\nu \xrightarrow
 {
  \substack{C^\infty_{\mathrm{loc}}}
 }
 x^\ell_1.$$
By induction hypothesis we know that
$$\lim_{s \to \infty}x^{\ell-1}_0(s)=\lim_{s \to \infty}x^{\ell-1}_1(s)$$
and consequently
$$\lim_{s \to -\infty}x^\ell_0(s)=\lim_{s \to -\infty}x^\ell_1(s).$$
Therefore we infer from Proposition~\ref{repro} that
$$r^\ell:=\lim_{\nu \to \infty}(r^\ell_{0,\nu}-r^\ell_{1,\nu})$$
exists and satisfies
$$x^\ell_0=r^\ell_*x^\ell_1.$$
This finishes the proof of the induction step and the Theorem follows. \hfill $\square$

\subsection{Floer-Gromov compactness for broken gradient flow lines}

We have explained so far the notion of Floer-Gromov convergence for gradient flow lines. More generally we can define the notion of Floer-Gromov convergence for a sequence of broken gradient flow lines. 
\begin{fed}
Assume that $y_\nu=\{x^k_\nu\}_{1\leq k \leq m_\nu}$ is a sequence of broken gradient flow lines from a critical point
$x^-$ to a critical point $x^+$ and $y=\{x^k\}$ is a broken gradient flow line from $x^-$ to $x^+$. We say that
$y_\nu$ \emph{Floer-Gromov converges} to $y$ if for every $k \in \{1,\ldots,m\}$ there exists a sequence of
real numbers $r^k_\nu$ as well as a sequence $j^k_\nu \in \{1,\ldots,m_\nu\}$ such that
$$
(r_\nu^k)_*x_\nu^{j^k_\nu} \xrightarrow
 {
  \substack{C^\infty_{\mathrm{loc}}}
 }
 x^k.$$
\end{fed}
In order to generalize results about Floer-Gromov convergence for sequences of gradient flow lines to sequences of broken
gradient flow lines the following notion is helpful. 
\begin{fed}
A sequence $y_\nu$ of broken gradient flow lines from $x^-$ to $x^+$ is called \emph{tame} if the number of gradient flow lines in each broken gradient flow
line $y_\nu$ is fixed independent of $\nu$ and moreover the asymptotics are fixed as well. That means that there exists
$n \in \mathbb{N}$ and critical points $c^k$ for $0 \leq k \leq n$ with $c^0=x^-$ and $c^{n+1}=x^+$ such that
$$y_\nu=\{x^k_\nu\}_{1\leq k \leq n}$$
and 
$$\lim_{s \to -\infty}x^k_\nu(s)=c^{k-1},\quad \lim_{s \to \infty}x^k_\nu(s)=c^k, \qquad k \in \{1,\ldots,n\}.$$ 
\end{fed}
\begin{lemma}\label{tame}
Each sequence of broken gradient flow lines from $x^-$ to $x^+$ has a tame subsequence.
\end{lemma}
\textbf{Proof: } This follows from the fact that the number of critical points of the Morse function $f$ is finite by Lemma~\ref{finite}. \hfill $\square$
\\ \\
Theorem~\ref{floergromov} about Floer-Gromov convergence of gradient flow lines generalizes to broken gradient flow lines as follows.
\begin{thm}\label{brokcomp}
Assume that $y_\nu$ is a sequence of broken gradient flow lines from $x^-$ to $x^+$. Then there exists a subsequence
$\nu_j$ and a broken gradient flow line $y$ from $x^-$ to $x^+$ such that
$$y_{\nu_j} \xrightarrow
 {
  \substack{\mathrm{Floer-Gromov}}
 }y.$$
\end{thm}
\textbf{Proof: } By Lemma~\ref{tame} we can assume maybe after transition to a subsequence that the sequence $y_\nu$ is tame.
We now apply Theorem~\ref{floergromov} to each sequence $x^k_\nu$ for fixed $k \in \{1,\ldots,n\}$ individually. This finishes the proof of the Theorem. \hfill $\square$
\\ \\
Theorem~\ref{unilim} about the uniqueness of limits up to equivalence generalizes to sequences of broken gradient flow lines as well. Namely we have the following Theorem. 
\begin{thm}\label{brokunilim}
Assume that $y_\nu$ is a sequence of broken gradient flow lines from $x^-$ to $x^+$ which converges in the sense of Floer-Gromov to broken gradient flow lines $y$ and $y'$. Then $y$ and $y'$ are equivalent. 
\end{thm}
\textbf{Proof: } Again we can assume by Lemma~\ref{tame} maybe after transition to a subsequence that $y_\nu$ is tame. Now the result follows from Theorem~\ref{unilim}. 

\subsection{Moduli spaces}\label{moduli}

Given a Morse function $f$ and a Riemannian metric $g$ on a closed manifold $M$ and given critical points $x^-$ and
$x^+$ of $f$ we abbreviate by
$$\widetilde{\mathcal{M}}(f,g;x^-,x^+):=\big\{x\,\,\textrm{solution of}\,\,(\ref{para}): \lim_{s \to \pm \infty}x(s)=x^\pm\big\}$$
the moduli space of all parametrized gradient flow lines from $x^-$ to $x^+$. On this moduli space we have an action of
$\mathbb{R}$ by timeshift $(x,r) \mapsto r_* x$. If $x^- \neq x^+$, then this action is free. Otherwise if
$x^-=x^+$ then the moduli space consists just of the constant gradient flow line to $x^-=x^+$. Hence we assume in the following that
$$x^- \neq x^+$$
and abbreviate the quotient 
$$\mathcal{M}(f,g;x^-,x^+):=\widetilde{\mathcal{M}}(f,g;x^-,x^+)/\mathbb{R},$$
namely the moduli space of unparametrized gradient flow lines from $x^-$ to $x^+$.
The moduli space of parametrized gradient flow lines from $x^-$ to $x^+$ is a subset of the bigger moduli space consisting
of broken gradient flow lines from $x^-$ to $x^+$ 
$$\widetilde{\mathcal{M}}^b(f,g;x^-,x^+):=\big\{y\,\,\textrm{broken gradient flow line from}\,\,x^-\,\,\textrm{to}
\,\,x^+\big\}.$$
There is a natural stratification
$$\widetilde{\mathcal{M}}^b(f,g;x^-,x^+)=\bigsqcup_{k=1}^\infty
\widetilde{\mathcal{M}}^b_k(f,g;x^-,x^+)$$
where $\widetilde{\mathcal{M}}^b_k(f,g;x^-,x^+)$ is the moduli space of $(k-1)$-fold broken gradient flow lines. In particular, we have
$$\widetilde{\mathcal{M}}^b_1(f,g;x^-,x^+)=\widetilde{\mathcal{M}}(f,g;x^-,x^+).$$
On $\widetilde{\mathcal{M}}^b_k(f,g;x^-,x-+)$ we have a free action of $\mathbb{R}^k$ by componentwise timeshift, namely if
$r=(r^1,\ldots,r^k) \in \mathbb{R}^k$ and $y=(x^1,\ldots,x^k) \in \widetilde{\mathcal{M}}^b_k(f,g;x^-,x^+)$, then
$$r_*y=(r^1_*x^1,\ldots,r^k_* x^k).$$
We denote the quotient by
$$\mathcal{M}^b_k(f,g;x^-,x^+):=\widetilde{\mathcal{M}}^b_k(f,g;x^-,x^+)/\mathbb{R}^k,$$
namely the moduli space of unparametrized $(k-1)$-fold broken gradient flow lines. We finally set
$$\overline{\mathcal{M}}(f,g;x^-,x^+):=\bigsqcup_{k=1}^\infty \mathcal{M}^b_k(f,g;x^-,x^+)$$
the moduli space of unparametrized arbitrary times broken gradient flow lines. We think of this moduli space as a compactification of $\mathcal{M}(f,g;x^-,x^+)$ since by Theorem~\ref{brokcomp} every sequence in this space has a Floer-Gromov converging subsequence. Moreover, by Theorem~\ref{brokunilim} the limit is unique. 

\subsection{The Floer-Gromov topology}

We continue using the notation of the previous paragraph. What we achieved so far is that we have introduced the moduli space
of broken unparametrized gradient flow lines $\overline{\mathcal{M}}(f,g;x^-,x^+)$ and defined on this moduli space a notion of convergence, namely Floer-Gromov convergence, which has the property that each sequence has a converging sequence and moreover the limit of this sequence is unique. We finally want to endow this moduli space with a topology so that Floer-Gromov convergence can be explained as convergence in this topology. We refer to this topology as the Floer-Gromov topology.
The construction of this topology is analogous to the construction of the Gromov topology in Gromov-Witten theory explained in \cite[Section 5.6]{mcduff-salamon}. There the following axioms appeared, although the structure did not get a name. 
\begin{fed}
Suppose that $X$ is a set. A \emph{convergence structure} on $X$ is a subset 
$$\mathcal{C} \subset X \times X^\mathbb{N}$$
meeting the following conditions. 
\begin{description}
 \item[(Constant):] If $x_n=x_0$ for all $n \in \mathbb{N}$, then $(x_0,(x_n)_n) \in \mathcal{C}$. 
 \item[(Subsequence):] If $(x_0,(x_n)_n) \in \mathcal{C}$ and $g \colon \mathbb{N} \to \mathbb{N}$ is strictly increasing,
  then $(x_0,(x_{g(n)})_n) \in \mathcal{C}$. 
 \item[(Subsubsequence):] If for every strictly increasing function $g \colon \mathbb{N} \to \mathbb{N}$, there is a
  strictly increasing function $f \colon \mathbb{N} \to \mathbb{N}$ such that $(x_0,(x_{g \circ f(n)})_n) \in \mathcal{C}$,
  then $(x_0,(x_n)_n) \in \mathcal{C}$. 
 \item[(Diagonal):] If $(x_0,(x_k)_k) \in \mathcal{C}$ and $(x_k,(x_{k,n})_n) \in \mathcal{C}$ for every $k$, then there
  exist sequences $k_i, n_i \in \mathbb{N}$ such that $(x_0,(x_{k_i,n_i})_i) \in \mathcal{C}$.
 \item[(Uniqueness):] If $(x_0,(x_n)_n) \in \mathcal{C}$ and $(y_0,(x_n)_n) \in \mathcal{C}$, then $x_0=y_0$.  
\end{description}
A tupel $(X,\mathcal{C})$ consisting of a set $X$ and a convergence structure $\mathcal{C}$ on $X$ is referred to as a
\emph{convergence space}.
\end{fed} 
Intuitively one should think of an element $(x_0,(x_n)_n)$ in a convergence structure $\mathcal{C}$ as a sequence $(x_n)_n$ on $X$ converging to $x_0$. Following \cite[Section 5.6]{mcduff-salamon} we now define a topology on a convergence space which has the property that a sequence in the convergence space converges, if and only if it converges in the topology.
\begin{fed}
Assume that $(X,\mathcal{C})$ is a convergence space. We say that $U \subset X$ is open if and only if for every
$(x_0,(x_n)_n) \in \mathcal{C} \cap (U \times X^\mathbb{N})$ there exists $n_0$ such that $x_n \in U$ for every
$n \geq n_0$. 
\end{fed} 
For a convergence space $(X,\mathcal{C})$ abbreviate
$$\mathcal{U}(\mathcal{C}):=\{U \subset X: U\,\,\textrm{open}\} \subset 2^X.$$
\begin{lemma}\label{topo}
$\mathcal{U}(\mathcal{C})$ is a topology on $X$ which has the property that $(x_0,(x_n)_n) \in \mathcal{C}$ if and only if
the sequence $(x_n)_n$ converges to $x_0$ with respect to $\mathcal{U}(\mathcal{C})$.
\end{lemma}
\textbf{Proof: } That $\mathcal{U}(\mathcal{C})$ is a topology is immediate and is actually true for any subset
$\mathcal{C}$ of $X \times X^\mathbb{N}$. Moreover, the definition of the topology $\mathcal{U}(\mathcal{C})$ implies
that if $(x_0,(x_n)_n) \in \mathcal{C}$, then the sequence $(x_n)_n$ converges to $x_0$ in the topology $\mathcal{U}(\mathcal{C})$. Again this holds true for any subset of $X \times X^\mathbb{N}$. What is nontrivial is that if a sequence
$(x_n)_n$ converges to $x_0$ in the topology $\mathcal{U}(\mathcal{C})$, then $(x_0,(x_n)_n) \in \mathcal{C}$. This implication requires the axioms of a convergence structure and is carried out in \cite[Lemma 5.6.4]{mcduff-salamon}. \hfill $\square$ 
\\ \\
Given a Morse function $f$ on a closed Riemannian manifold $(M,g)$ and two different critical points $x^-$ and $x^+$ of $f$ we consider the moduli space
$$\overline{\mathcal{M}}=\overline{\mathcal{M}}(f,g;x^-,x^+)$$
of unparametrized broken gradient flow lines from $x^-$ to $x^+$. We introduce the subset
$$\mathcal{C}_{FG} \subset \overline{\mathcal{M}}\times \overline{\mathcal{M}}^\mathbb{N}$$
consisting of $(x_0,(x_n)_n)$ such that 
$$x_n \xrightarrow
 {
  \substack{\mathrm{Floer-Gromov}}
 }x_0.$$
\begin{lemma}
$\mathcal{C}_{FG}$ is a convergence structure on $\overline{\mathcal{M}}$.
\end{lemma}
\textbf{Proof: } The Constant, Subsequence and Subsubsequence axioms are obvious. The Diagonal axiom follows from the fact
that $C^\infty_{\mathrm{loc}}(\mathbb{R},M)$ has a countable neighbourhood base. Finally the Uniqueness Axiom follows from Theorem~\ref{brokunilim}. \hfill $\square$
\\ \\
In view of Lemma~\ref{topo} 
$$\mathcal{U}_{FG}:=\mathcal{U}(\mathcal{C}_{FG})$$
is a topology on $\overline{\mathcal{M}}$ which has the property that a sequence $(x_n)_n$ in $\overline{\mathcal{M}}$
is Floer-Gromov converging if and only if it converges in the topology $\mathcal{U}_{FG}$. We refer to $\mathcal{U}_{FG}$
as the \emph{Floer-Gromov topology} on $\overline{\mathcal{M}}$ and in the following we think of $\overline{\mathcal{M}}$
as a topological space endowed with the Floer-Gromov topology. With this convention we can rephrase Theorem~\ref{brokcomp}
in the following way.
\begin{cor}
$\overline{\mathcal{M}}$ is sequentially compact. 
\end{cor}
\subsection{Strategy to show that the boundary of the boundary vanishes}\label{strategy}

For general Riemannian metrics the moduli spaces of gradient flow lines might look rather nasty. However, we will show
that for generic metrics they look quite nice. What the term ``generic metric" precisely means, we discuss in Section~\ref{generic}. What suffices to know for the moment is that ``many" metrics are generic and in particular that such metrics exist. Intuitively one might think that given any metric one can wiggle a bit at this metric to make it generic. 
\\ \\
Suppose that $f \colon M \to \mathbb{R}$ is a Morse function and $c_1$ and $c_2$ are two different critical points of $f$. We will prove the following three theorems.
\begin{thm}\label{dd1}
For a generic metric $g$ the moduli space $\mathcal{M}(f,g;c_1,c_2)$ is a manifold of dimension
$$\mathrm{dim}\mathcal{M}(f,g;c_1,c_2)=\mu(c_1)-\mu(c_2)-1.$$
\end{thm}
\begin{thm}\label{dd2}
If $\mu(c_1)=\mu(c_2)+1$, then for generic Riemannian metric $g$ it holds that
$$\overline{\mathcal{M}}(f;g;c_1,c_2)=\mathcal{M}(f,g;c_1,c_2)$$
and therefore the moduli space $\mathcal{M}(f,g;c_1,c_2)$ is a compact, zero-dimensional manifold, i.e., a finite number of
points.
\end{thm}
\begin{thm}\label{dd3}
If $\mu(c_1)=\mu(c_2)+2$, then for generic Riemannian metric $g$ it holds that
$$\overline{\mathcal{M}}(f;g;c_1,c_2)=\mathcal{M}(f,g;c_1,c_2) \sqcup \mathcal{M}^b_2(f,g;c_1,c_2)$$
and $\overline{\mathcal{M}}(f,g;c_1,c_2)$ has the structure of a one-dimensional manifold with boundary such that
$$\partial \overline{\mathcal{M}}(f,g;c_1,c_2)=\mathcal{M}^b_2(f,g;c_1,c_2).$$
\end{thm}
The precise version of Theorem~\ref{dd1} is Theorem~\ref{dimpar}, the one of Theorem~\ref{dd2} is Corollary~\ref{dd2f}, and
the one of Theorem~\ref{dd3} is Theorem~\ref{dd3f}.
\\ \\
From Theorem~\ref{dd2} we get that in the case of Morse index difference two 
$$\#_2\mathcal{M}(f,g;c_1,c_2) \in \mathbb{Z}_2$$
is welldefined so that by setting for $c \in \mathrm{crit}(f)$
$$\partial c=\sum_{\substack{c' \in \mathrm{crit}(f)\\ \mu(c')=\mu(c)-1}}\#_2 \mathcal{M}(f,g;c,c')c'$$
we obtain a welldefined linear map
$$\partial \colon CM_*(f) \to CM_{*-1}(f).$$
Theorem~\ref{dd3} tells us that in case the Morse index difference is two the boundary of the compactified moduli space consists of once broken gradient flow lines. We can write that is 
$$\partial \overline{\mathcal{M}}(f,g;c_1,c_2)=\bigcup_{c \in \mathrm{crit}(f)} \mathcal{M}(f,g;c_1,c)\times
\mathcal{M}(f,g;c,c_2).$$
Since $\mu(c_1)=\mu(c_2)+2$ we infer from Theorem~\ref{dd1} that in order to contribute nontrivially to the union on the righthand side a critical point $c$ has to satisfy
$$\mu(c)=\mu(c_1)-1=\mu(c_2)+1$$
so that we can write the formula above equivalently as
$$\partial \overline{\mathcal{M}}(f,g;c_1,c_2)=\bigcup_{\substack{c \in \mathrm{crit}(f)\\\mu(c)=\mu(c_1)-1}} \mathcal{M}(f,g;c_1,c)\times
\mathcal{M}(f,g;c,c_2).$$
To see how this formula implies that $\partial^2=0$ we compute for a critical point $c$
\begin{eqnarray}\label{rand}
\partial^2 c&=&\partial\Bigg(\sum_{\substack{c' \in \mathrm{crit}(f)\\ \mu(c')=\mu(c)-1}}\#_2\mathcal{M}(f,g;c,c')c'\Bigg)
\\ \nonumber
&=&\sum_{\substack{c' \in \mathrm{crit}(f)\\ \mu(c')=\mu(c)-1}}\#_2\mathcal{M}(f,g;c,c') \partial c'\\ \nonumber
&=&\sum_{\substack{c' \in \mathrm{crit}(f)\\ \mu(c')=\mu(c)-1}}\#_2\mathcal{M}(f,g;c,c') 
\Bigg(\sum_{\substack{c'' \in \mathrm{crit}(f)\\ \mu(c'')=\mu(c')-1}}\#_2 \mathcal{M}(f,g;c',c'')c''\Bigg)\\ \nonumber
&=&\sum_{\substack{c'' \in \mathrm{crit}(f)\\ \mu(c'')=\mu(c)-2}}\Bigg(\sum_{\substack{c' \in \mathrm{crit}(f)\\
\mu(c')=\mu(c)-1}}\#_2 \mathcal{M}(f,g;c,c')\cdot \#_2 \mathcal{M}(f,g;c',c'')\Bigg)c''\\ \nonumber
&=&\sum_{\substack{c'' \in \mathrm{crit}(f)\\ \mu(c'')=\mu(c)-2}}\#_2 \partial \overline{\mathcal{M}}(f,g;c,c'')c''
\end{eqnarray}
\begin{figure}[h]
\begin{center}
 \includegraphics[scale=1]{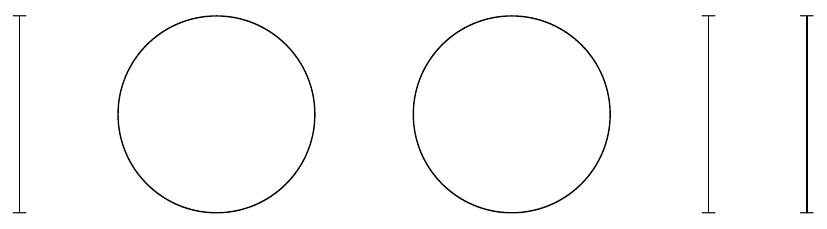}
 \caption{Conpact one-dimensional manifold with boundary.}
\end{center}
\end{figure}
By Theorem~\ref{dd3} the compactified moduli space $\overline{\mathcal{M}}(f,g;c,c'')$ is a compact one-dimensional manifold with boundary. But a compact one-dimensional manifold with boundary is diffeomorphic to a finite disjoint union of circles and intervals, see for instance \cite[Appendix]{milnor2}. In particular, the number of boundary points of a compact one-dimensional manifold is even, namely two times the number of intervals occuring in the disjoint union. Therefore (\ref{rand})
implies that
$$\partial^2=0,$$
i.e., a boundary has no boundary.

\section{Fredholm theory}

\subsection{Fredholm operators}

In this section we recall some basic facts about Fredholm operators. Suppose that $X$ and $Y$ are Banach spaces.
\begin{fed}
A bounded linear operator $D \colon X \to Y$ is called Fredholm if it meets the following conditions.
\begin{description}
 \item[(i)] The kernel of $D$ is finite dimensional.
 \item[(ii)] The image of $D$ in $Y$ is closed.
 \item[(iii)] The cokernel of $D$ is finite dimensional. 
\end{description}
\end{fed}
In our applications the Banach spaces $X$ and $Y$ are always Hilbert spaces, so that in the following we always assume that
$X$ and $Y$ are Hilbert spaces. For Hilbert spaces the last condition can be rephrased as
$$\mathrm{dim} D(X)^\perp<\infty$$
where $D(X)^\perp$ is the orthogonal complement of $D(X)$ in $Y$.
\begin{fed}
Assume that $D \colon X \to Y$ is a Fredholm operator. Then
$$\mathrm{ind}(D):=\mathrm{dim} (\mathrm{ker}D)-\mathrm{dim}D(X)^\perp \in \mathbb{Z}$$
is called the \emph{index} of $D$.
\end{fed}
It is interesting to examine this concept in the finite dimensional case. Hence suppose that $X$ and $Y$ are finite dimensional vector spaces. Then every linear map $D \colon X \to Y$ is Fredholm and the following formulas hold
$$\mathrm{dim}(X)=\mathrm{dim}(\mathrm{ker}D)+\mathrm{dim}D(X), \quad \mathrm{dim}(Y)=\mathrm{dim}D(X)+
\mathrm{dim}D(X)^\perp.$$
Hence the index of $D$ is simply the difference of the dimensions of $X$ and $Y$
$$\mathrm{ind}(D)=\mathrm{dim}(X)-\mathrm{dim}(Y).$$
In particular, the index is independent of $D$ and only depends on the vector spaces $X$ and $Y$. This is in sharp contrast to the dimensions of the kernel and cokernel which depend on the linear map $D$ and can jump even under small perturbations of $D$.
\\ \\
In infinite dimensions the index cannot be expressed anymore just in terms of $X$ and $Y$. However, both the Fredholm property as well as the index are stable under perturbations. Recall that the space $\mathcal{L}(X,Y)$ of bounded linear
operators from $X$ to $Y$ becomes itself a Banach space with norm the operator norm defined for $D \in \mathcal{L}(X,Y)$ by
$$||D||:=\sup_{x \neq 0 \in X}\frac{||Dx||_Y}{||x||_X}.$$
Using the operator norm on $\mathcal{L}(X,Y)$ we can now state the first stability theorem as follows.
\begin{thm}\label{stab1}
Suppose that $D \colon X \to Y$ is a Fredholm operator. Then there exists $\epsilon=\epsilon(D)>0$ such that every
$D' \in \mathcal{L}(X,Y)$ satisfying
$$||D'-D||<\epsilon$$
is itself a Fredholm operator and its index satisfies
$$\mathrm{ind}(D')=\mathrm{ind}(D).$$ 
\end{thm}
The theorem tells us that the Fredholm property is open on $\mathcal{L}(X,Y)$ and the index is constant on connected components of the subset of Fredholm operators. This is extremely useful in applications. In fact given a Fredholm operator 
it is often difficult to compute its index directly by computing its kernel and cokernel. However, what often works is that the given operator can be homotoped through Fredholm operators to a simpler operator for which kernel and cokernel can be computed. The stability theorem then tells us that the Fredholm index of the original operator did not change under the homotopy although the dimensions of the kernel and cokernel might well have. 
\\ \\
The second stability property is that the Fredholm property as well as the index remain unchanged under addition of compact operators. Recall that a linear operator
$$K \colon X \to Y$$
is called \emph{compact} if $\overline{K(B_1)}$ is compact in $Y$, where $B_1=\{x \in X: ||x||=1\}$ is the 1-ball around the origin in $X$ and $\overline{K(B_1)}$ is the closure of its image in $Y$. The second stability theorem is the following.
\begin{thm}\label{stab2}
Suppose that $D \colon X \to Y$ is a Fredholm operator and $K \colon X \to Y$ is a compact operator, then its sum
$$D+K \colon X \to Y$$
is a Fredholm operator as well and its index satisfies
$$\mathrm{ind}(D+K)=\mathrm{ind}(D).$$
\end{thm}
To check in practize that a given operator is Fredholm one usually has to prove an estimate. Namely one has the following lemma. 
\begin{lemma}\label{semifredi0}
Suppose that $X,Y,Z$ are Banach spaces and $D \colon X \to Y$ is a bounded operator. Suppose that there exists a compact linear operator $K \colon X \to Z$ and a constant $c>0$ such that 
\begin{equation}\label{semiest}
||x||_X \leq c\big(||Dx||_Y+||Kx||_Z\big), \quad \forall\,\,x \in X.
\end{equation}
Then $D(X)$ is closed in $Y$ and $\mathrm{ker}(D)$ is finite dimensional. 
\end{lemma}
An operator which has a finite dimensional kernel and a closed image but not necessarily a finite dimensional cokernel is referred to as a \emph{semi-Fredholm operator} and therefore an estimate of the form (\ref{semiest}) is called a 
\emph{semi-Fredholm estimate}. Although a semi-Fredholm estimate does not yet guarantee the third property of a Fredholm operator it can often used to derive this as well. The trick is to look at the adjoint of $D$. If this has a similar form as
$D$ itself the semi-Fredholm estimate can as well be applied to the adjoint and since the kernel of the adjoint coincides with the cokernel of $D$ this than establishes the Fredholm property of $D$.
\\ \\
Proofs of the results mentioned in this paragraph can for example be found in \cite[Appendix A.1]{mcduff-salamon}.

\subsection{The linearized gradient flow equation}

Suppose that $A \in C^0\big(\mathbb{R},\mathrm{End}(\mathbb{R}^n)\big)$ is a continuous family of endomorphisms from 
$\mathbb{R}^n$ to itself which we can write as a continuous family of $n\times n$-matrices. Assume that the limit operators
$$\lim_{s \to \pm \infty} A(s)=A^\pm$$
exist and are symmetric and nondegenerate in the sense that $\mathrm{ker}A^\pm=\{0\}$. We do neither require that $A(s)$ is symmetric for finite $s$ nor that it is nondegenerate. We abbreviate by
$$\mathfrak{A} \subset C^0\big(\mathbb{R},\mathrm{End}(\mathbb{R}^n)\big)$$
the space of all families of matrices satisfying the above described asymptotic behaviour. For $A \in \mathfrak{A}$ we abbreviate by
$$\mu(A^\pm) \in \{0,\ldots,n\}$$
the Morse index of the asymptotic operators $A^\pm$, namely the number of negative eigenvalues of the symmetric matrices
$A^\pm$ counted with multiplicity. Moreover, we introduce the following bounded linear operator
$$D_A \colon W^{1,2}(\mathbb{R},\mathbb{R}^n) \to L^2(\mathbb{R},\mathbb{R}^n)$$
which is given for $\xi \in W^{1,2}(\mathbb{R},\mathbb{R}^n)$ by
\begin{equation}\label{da}
D_A \xi(s)=\partial_s \xi(s)+A(s) \xi(s), \quad s \in \mathbb{R}.
\end{equation}
\begin{thm}\label{index}
Suppose that $A \in \mathfrak{A}$. Then $D_A$ is a Fredholm operator and its index is given by
$$\mathrm{ind}(D_A)=\mu(A^-)-\mu(A^+).$$
\end{thm}
We prove Theorem~\ref{index} in Section~\ref{indcomp}. Let us explain how linearizing a gradient flow line gives rise to an
operator $D_A$. Suppose that $x$ is a gradient flow line of $\nabla f$ converging asymptotically to critical points $x^\pm$
of $f$. We assume that the closure of $x$ is completely contained in a chart $U \subset \mathbb{R^n}$. Since gradient flow lines flow downhill such a chart always exists. In coordinates the gradient flow equation becomes
$$\partial_s x_i+\sum_{j=1}^n g_{ij}(x)\frac{\partial f}{\partial x_j}=0$$
where $g_{ij}$ are the coefficients of the inverse of the metric $g$.  Linearizing the gradient flow equation one obtains
the equation
$$\partial_s \xi_i+\sum_{j=1}^n\sum_{k=1}^n\frac{\partial g_{ij}(x)}{\partial x_k}\xi_k\frac{\partial f}{\partial x_j}
+\sum_{j=1}^n \sum_{k=1}^n g_{ij}(x)\frac{\partial^2 f}{\partial x_j \partial x_k}(x)\xi_k=0.$$
Defining
$$A_{ik}=\sum_{j=1}^n\Bigg(\frac{\partial g_{ij}(x)}{\partial x_k}\frac{\partial f}{\partial x_j}(x)+g_{ij}(x)
\frac{\partial^2 f}{\partial x_j\partial x_k}(x)\Bigg)$$
we can write this as
$$\partial_s \xi_i+\sum_{k=1}^n A_{ik}\xi_k.$$
Note that because $x$ converges asymptotically to critical points we have the asymptotics
$$\lim_{s \to \pm \infty} A_{ik}(s)=\sum_{j=1}^n g_{ij}(x^\pm)\frac{\partial^2 f}{\partial x_j \partial x_k}(x^\pm)$$
so that by choosing proper coordinates around $x^\pm$ we can assume that the asymptotic matrices are symmetric and nondegenerate. 
\\ \\
For later reference we as well introduce linear operators one obtains by linearizing partial gradient flow lines, i.e.,
gradient flow lines restricted to a subset of the real line. For that purpose we introduce subsets
$$\mathfrak{A}_- \subset C^0\big((-\infty,0],\mathrm{End}(\mathbb{R}^n)\big), \quad \mathfrak{A}_+\subset C^0\big([0,\infty),\mathrm{End}(\mathbb{R}^n)\big)$$
where $A \in \mathfrak{A}_-$ if there exists a symmetric, nondegenerate $n \times n$-matrix $A^-$ such that
$$\lim_{s \to -\infty}A(s)=A^-$$
and similarly, $A \in \mathfrak{A}_+$ if
$$\lim_{s \to \infty}A(s)=A^+$$
where $A^+$ as well is required to be symmetric and nondegenerate. Moreover, for $T>0$ we abbreviate
$$\mathfrak{A}_T=C^0\big([-T,T],\mathrm{End}(\mathbb{R}^n)\big).$$
For $A$ in $\mathfrak{A}_-$, $\mathfrak{A}_+$ or $\mathfrak{A}_T$ we introduce operators $D_A$ is in (\ref{da}) where
the domains of the Sobolev spaces $W^{1,2}$ and $L^2$ have to be restricted to the corresponding domain of the operator $A$. The analogon of Theorem~\ref{index} for the restricted operators is the following result. 
\begin{thm}\label{index2}
Suppose that $A$ belongs to $\mathfrak{A}_-$, $\mathfrak{A}_+$ or $\mathfrak{A}_T$. Then $D_A$ is a Fredholm operator and for the index the following formulas hold.
\begin{description}
 \item[(i)] If $A \in \mathfrak{A}_-$, then $\mathrm{ind}(D_A)=\mu(A^-)$.
 \item[(ii)] If $A \in \mathfrak{A}_+$, then $\mathrm{ind}(D_A)=n-\mu(A^+)$.
 \item[(iii)] If $A \in \mathfrak{A}_T$, then $\mathrm{ind}(D_A)=n$.
\end{description}
Moreover, $D_A$ is always surjective so that we have
$$\mathrm{dim}\ker D_A=\mathrm{ind}(D_A).$$
\end{thm}
The proof of Theorem~\ref{index2} is carried out in Section~\ref{indcomp}.

\subsection{The semi-Fredholm estimate}

The crucial ingredient to prove that the operators $D_A$ are Fredholm is a semi-Fredholm estimate (\ref{semiest}). The following lemma provides such an estimate. 

\begin{lemma}\label{semifredi}
Suppose that $A \in \mathfrak{A}$. There exists $T>0$ and $c>0$ such that for every $\xi \in W^{1,2}(\mathbb{R})$
we have the following estimate
$$||\xi||_{W^{1,2}(\mathbb{R})} \leq c\big(||\xi||_{L^2([-T,T])}+||D_A \xi||_{L^2(\mathbb{R})}\big).$$
\end{lemma}
The reader might wonder if it were not sufficient to prove an estimate of the form
$$||\xi||_{W^{1,2}(\mathbb{R})}\leq c\big(||\xi||_{L^2(\mathbb{R})}+||D_A\xi||_{L^2(\mathbb{R})}\big).$$
This were much easier and is indeed the first step of the proof of Lemma~\ref{semifredi}. Unfortunately, while the restriction operator
$$W^{1,2}(\mathbb{R}) \to L^2([-T,T])\quad \xi \mapsto \xi_{[-T,T]}$$
is compact, this is not true anymore for the inclusion operator from $W^{1,2}(\mathbb{R})$ into $L^2(\mathbb{R})$. To see that
pick a compactly supported bump function $\beta$ satisfying
$$||\beta||_{1,2}=1.$$
For $\nu \in \mathbb{N}$ consider the shifted bump function
$$\beta_\nu \in W^{1,2}(\mathbb{R}), \quad s \mapsto \beta(s-\nu).$$
Note that 
$$||\beta_\nu||_{1,2}=1, \quad \forall\,\,\nu \in \mathbb{N}$$
so that $\beta_\nu$ is a sequence in the unit ball of $W^{1,2}(\mathbb{R})$ which has no convergent subsequence. 
\\ \\
\textbf{Proof of Lemma~\ref{semifredi} } The proof of the Lemma has three steps. We follow the arguments of Robbin and Salamon in \cite{robbin-salamon}.
\\ \\
\textbf{Step\,1: } \emph{There exists a constant $c>0$ such that for every $\xi \in W^{1,2}(\mathbb{R})$ the following estimate holds}
$$||\xi||_{W^{1,2}(\mathbb{R})}\leq c\big(||\xi||_{L^2(\mathbb{R})}+||D_A\xi||_{L^2(\mathbb{R})}\big).$$
\textbf{Proof of Step\,1: } Since $\lim_{s \to \pm \infty}A=A^\pm$ there exists $c_0>0$ such that 
the operator norm of $A(s)$ is uniformly bounded by $c_0$, i.e.,
$$||A(s)||_{\mathcal{L}(\mathbb{R}^n)}:=\sup_{v \in \mathbb{R}^n \setminus \{0\}}\bigg\{\frac{||Av||}{||v||}\bigg\}\leq c_0.$$
Using the equation
$$\partial_s \xi=D_A \xi-A \xi$$
we estimate using the triangle inequality
\begin{eqnarray*}
||\partial_s \xi||_{L^2(\mathbb{R}} &\leq& ||D_A \xi||_{L^2(\mathbb{R}}+||A\xi||_{L^2(\mathbb{R})}\\
 &\leq& ||D_A \xi||_{L^2(\mathbb{R}}+c_0||\xi||_{L^2(\mathbb{R})}.
\end{eqnarray*}
From that we infer
\begin{eqnarray*}
||\xi||_{W^{1,2}(\mathbb{R})}&=&\sqrt{||\xi||^2_{L^2(\mathbb{R})}+||\partial_s \xi||^2_{L^2(\mathbb{R})}}\\
&\leq&||\xi||_{L^2(\mathbb{R})}+||\partial_s \xi||_{L^2(\mathbb{R})}\\
&\leq&||D_A \xi||_{L^2(\mathbb{R})}+(c_0+1)||\xi||_{L^2(\mathbb{R})}\\
&\leq&(c_0+1)\big(||\xi||_{L^2(\mathbb{R})}+||D_A\xi||_{L^2(\mathbb{R})}\big).
\end{eqnarray*}
This finishes the proof of Step\,1.
\\ \\
\textbf{Step\,2: } \emph{Assume that $A(s)=A_0$ is constant. Then there exists a constant $c>0$ such that}
$$||\xi||_{W^{1,2}(\mathbb{R})}\leq c||D_{A_0} \xi||_{L^2(\mathbb{R})}.$$
\textbf{Proof of Step\,2: } Because $A_0 \in \mathfrak{A}$ it follows that it is symmetric, i.e., $A_0=A_0^T$, and
nondegenerate meaning that its kernel is trivial. We decompose the proof of Step\,2 into three substeps.
\\ \\
\textbf{Step\,2a: } \emph{Assume that $A_0$ is positive definite.}
\\ \\
\textbf{Proof of Step\,2a: } For $\xi \in W^{1,2}(\mathbb{R})$ define
\begin{equation}\label{et}
\eta:=D_{A_0}\xi=\partial_s \xi+A_0\xi.
\end{equation}
Given $\eta$ this is an ODE for $\xi$ and its general solution is given for $c \in \mathbb{R}^n$ by
$$\xi_c(s)=\int_{-\infty}^s e^{A_0(t-s)}\eta(t)dt+ce^{-A_0 s}.$$
The only solution in $L^2(\mathbb{R},\mathbb{R}^n)$ is $\xi=\xi_0$. Therefore if we set
$$\Phi(s):=\left\{\begin{array}{cc}
e^{-A_0 s} & s \geq 0\\
0 & s<0
\end{array}\right.$$
we can write $\xi$ as the convolution of $\eta$ with $\Phi$, i.e.,
$$\xi(s)=\Phi  * \eta(s)=\int_{-\infty}^\infty \Phi(s-t)\eta(t)dt.$$
We recall Young's inequality. Namely if $f \in L^p$, $g \in L^q$ and $1 \leq p,q,r \leq \infty$ satisfy
$$\frac{1}{p}+\frac{1}{q}=\frac{1}{r}+1$$
then
$$||f * g||_r \leq ||f||_p \cdot ||g||_q.$$
Choosing $r=2$, $p=1$, and $q=2$ we obtain the inequality
\begin{equation}\label{young}
||\xi||_2 \leq ||\Phi||_1\cdot ||\eta||_2.
\end{equation}
Note that because $A_0$ is positive definite we have
$$||\Phi||_1<\infty.$$
By (\ref{et}) we can write
$$\partial_s \xi=\eta-A_0\xi$$
so that using (\ref{young}) we obtain the estimate
\begin{eqnarray}\label{young2}
||\partial_s \xi||_2 \leq ||\eta||_2+||A_0\xi||_2 \leq \big(||A_0||\cdot ||\Phi||_1+1\big)||\eta||_2
\end{eqnarray}
where $||A_0||=||A_0||_{\mathcal{L}(\mathbb{R}^n)}$ is the operator norm of $A_0$. Combining (\ref{young}) and
(\ref{young2}) and using that by definition $\eta=D_{A_0}\xi$ we get the inequality
$$||\xi||_{1,2} \leq ||\xi||_2+||\partial_s \xi||_2 \leq \big((||A_0||+1)||\Phi||_1+1\big)||D_{A_0}\xi||_2.$$
This finishes the proof of Step\,2a.
\\ \\
\textbf{Step\,2b: } \emph{Assume that $A_0$ is negative definite.}
\\ \\
For $\xi \in W^{1,2}$ we define $\xi^- \in W^{1,2}$ by the formula
$$\xi^-(s):=-\xi(-s), \quad s \in \mathbb{R}.$$
Note that
$$D_{-A_0}\xi^-(s)=\partial_s \xi(-s)+A_0(\xi(-s),\quad s \in \mathbb{R}.$$
Since $-A_0$ is positive definite we infer using Step\,2a 
$$||\xi||_{1,2}=||\xi^-||_{1,2} \leq c||D_{-A_0}\xi^-||_2=||D_{A_0}\xi||_2.$$
This proves Step\,2b. 
\\ \\
\textbf{Step\,2c: } \emph{Proof of Step\,2.}
\\ \\
We decompose
$$\mathbb{R}^n=V_- \oplus V_+$$
where $V_-$ is the direct sum of the eigenspaces of $A_0$ to negative eigenvalues and $V_+$ is the direct sum of
eigenspaces to positive eigenvalues. Note that because $A_0$ is symmetric we have
$$V_- \perp V_+$$
i.e., $V_-$ is orthogonal to $V_+$. If $\xi \in W^{1,2}(\mathbb{R},\mathbb{R}^n)$ we decompose $\xi$ 
$$\xi=\xi_-+\xi_+, \qquad \xi_- \in W^{1,2}(\mathbb{R},V_-),\quad \xi_+ \in W^{1,2}(\mathbb{R},V_+).$$
From Step\,2a and Step\,2b we infer that for a constant $c$ it holds that
\begin{eqnarray*}
||\xi_-||_{1,2} \leq c||D_{A_0}\xi_-||_2,\\
||\xi_+||_{1,2} \leq c||D_{A_0}\xi_+||_2.
\end{eqnarray*}
Because $V_-$ and $V_+$ are orthogonal we have for $s \in \mathbb{R}$ with respect to the standard norm in
$\mathbb{R}^n$
$$||D_{A_0}\xi(s)||^2=||D_{A_0}\xi_-(s)||^2+||D_{A_0}\xi_+(s)||^2$$
from which we obtain by integration
$$||D_{A_0}\xi||_2^2=||D_{A_0}\xi_-||_2^2+||D_{A_0}\xi_+||_2^2.$$
Therefore
$$||D_{A_0}\xi||_2=\sqrt{||D_{A_0}\xi_-||^2_2+||D_{A_0}\xi_+||^2_2}\geq \frac{1}{\sqrt{2}}\big(||D_{A_0}\xi_-||_2+
||D_{A_0}\xi_+||_2\big).$$
Using this we estimate
\begin{eqnarray*}
||\xi||_{1,2}&=&||\xi_-+\xi_+||_{1,2}\\
&\leq&||\xi_-||_{1,2}+||\xi_+||_{1,2}\\
&\leq&c \big(||D_{A_0}\xi_-||_2+||D_{A_0}\xi_+||_2\big)\\
&\leq&\sqrt{2}\cdot c||D_{A_0}\xi||_2
\end{eqnarray*}
This finishes the proof of Step\,2.
\\ \\
\textbf{Step\,3: } \emph{We prove the Lemma}
\\ \\
Let $c>0$ be a constant so that the assertion of Step\,1 is true. Moreover, by Step\,2 maybe after choosing $c$ bigger we may assume in addition that 
$$||\xi||_{1,2}\leq c||D_{A^\pm}\xi||_2, \quad \forall\,\,\xi \in W^{1,2}(\mathbb{R},\mathbb{R}^n).$$
Since $\lim_{s \to \pm \infty}A(s)=A^\pm$ there exists $T>1$ with the property that
$$||A(s)-A^\pm||_{\mathcal{L}(\mathbb{R}^n)}\leq \frac{1}{2c},\quad \pm s \geq T-1.$$
Suppose that $\xi \in W^{1,2}(\mathbb{R},\mathbb{R}^n)$ satisfies $\xi(s)=0$ for $s \leq T-1$. We estimate
\begin{eqnarray*}
||\xi||_{1,2} &\leq& c||D_{A^+}\xi||_2\\
&=&c||D_A\xi+(A^+-A)\xi||_2\\
&\leq & c\big(||D_A\xi||_2+||(A-A^+)\xi||_2\big)\\
&\leq& c\big(||D_A\xi||+\tfrac{1}{2c}||\xi||_2\big)\\
&\leq & c||D_A\xi||_2+\tfrac{1}{2}||\xi||_{1,2}.
\end{eqnarray*}
In the second last inequality we have used that $\xi(s)$ vanishes for $s \leq T-1$. Hence we obtain
\begin{equation}\label{ab}
||\xi||_{1,2} \leq 2c||D_A\xi||_2
\end{equation}
In the same way we obtain inequality (\ref{ab}) for $\xi \in W^{1,2}(\mathbb{R},\mathbb{R}^n)$ satisfying
$\xi(s)=0$ for $s \geq -T+1$. 
We choose a smooth cutoff function $\beta \in C^\infty(\mathbb{R},[0,1])$ satisfying
$$\beta(s)=0,\quad |s| \geq T,\qquad \beta(s)=1,\quad |s|\leq T-1.$$
With this preparations we can estimate now for a general $\xi \in W^{1,2}$ 
\begin{eqnarray*}
||\xi||_{1,2}&=&||\beta \xi+(1-\beta)\xi||_{1,2}\\
&\leq&||\beta \xi||_{1,2}+||(1-\beta)\xi||_{1,2}\\
&\leq& 2c\Big(||\beta \xi||_2+||D_A(\beta \xi)||_2+||D_A((1-\beta)\xi)||_2\Big)\\
&\leq& 2c\Big(||\beta \xi||_2+||\beta D_A\xi+(\partial_s\beta)\xi||_2+||(1-\beta)D_A\xi-(\partial_s \beta)\xi||_2\Big)\\
&\leq& 2c\Big(||\beta \xi||_2+||\beta D_A\xi||_2+||(1-\beta)D_A\xi||_2+2||(\partial_s \beta)\xi||_2\Big)\\
&\leq& 2c\Big(||\xi||_{L^2([-T,T])}+2||D_A \xi||_2+2||\partial_s \beta||_\infty \cdot ||\xi||_{L^2([-T,T])}\Big)\\
&\leq& 2c(2+||\partial_s \beta||_\infty)\big(||\xi||_{L^2([-T,T])}+||D_A\xi||_{L^2(\mathbb{R})}\big).
\end{eqnarray*}
Here we have used in the second inequality Step\,1 and (\ref{ab}). This finishes the proof of the Lemma. \hfill $\square$

\subsection{Adjoints}\label{adjoint}

If $A \in \mathfrak{A}$ then Lemma~\ref{semifredi0} and Lemma~\ref{semifredi} imply that the bounded linear operator
$D_A \colon W^{1,2}(\mathbb{R}) \to L^2(\mathbb{R})$ has a finite dimensional kernel and a closed image. To show that it is actually Fredholm we need to check that its cokernel is finite dimensional as well. Hence suppose that $\eta \in L^2(\mathbb{R})$ lies in the orthogonal complement of the image of $D_A$, i.e.,
$$\langle \eta, D_A \xi\rangle=0, \quad \forall\,\,\xi \in W^{1,2}(\mathbb{R}).$$
By definition of $D_A$ this means that for every $\xi \in W^{1,2}(\mathbb{R})$ 
$$\int_{-\infty}^\infty\langle \eta(s),\partial_s \xi(s)\rangle ds+\int_{-\infty}^\infty \langle \eta(s),A(s)\xi(s)\rangle ds=0$$
which we rearrange to
$$\int_{-\infty}^\infty \langle \eta(s),\partial_s \xi(s)\rangle ds=-\int_{-\infty}^\infty\langle A^T(s)\eta(s),\xi(s)\rangle ds.$$
But this means that $\eta$ has a weak derivative given by $A^T \eta$. In particular, $\eta \in W^{1,2}(\mathbb{R})$ 
and 
$$\partial_s \eta(s)=A^T(s)\eta(s), \quad s \in \mathbb{R}.$$
In other words
$$\eta \in \mathrm{ker}D_{-A^T}.$$
However $-A^T$ lies as well in $\mathfrak{A}$ so that $D_{-A^T}$ is as well a semi-Fredholm operator. In particular
$$\mathrm{coker}D_A=\mathrm{ker}D_{A^T}$$
is finite dimensional. We have proved the following Proposition. 
\begin{prop}\label{fre1}
If $A \in \mathfrak{A}$, then $D_A$ is a Fredholm operator. 
\end{prop}
We next consider the case where $A \in \mathfrak{A}_+$. We show that 
$$D_A \colon W^{1,2}\big([0,\infty)\big) \to L^2\big([0,\infty)\big)$$
is surjective. Suppose that $\eta$ is in the orthogonal complement of the image of $D_A$. As in the previous case we 
obtain the following condition on $\eta$
\begin{equation}\label{ad1}
\int_0^\infty\langle \eta(s),\partial_s \xi(s)\rangle ds=-\int_0^\infty\langle A^T(s)\eta(s),\xi(s)\rangle ds,\quad
\forall\,\,\xi \in W^{1,2}\big([0,\infty)\big).
\end{equation}
Again this implies that $\eta \in W^{1,2}\big([0,\infty)\big)$ and
\begin{equation}\label{ad2}
\partial_s \eta(s)=A^T(s)\eta(s), \quad s \geq 0.
\end{equation}
Given $\delta>0$ we can choose a smooth monotone decreasing cutoff function $\beta_\delta \in C^\infty\big([0,\infty),[0,1]\big)$ satisfying
$$\beta_\delta(0)=1,\qquad \beta_\delta(s)=0,\quad s \geq \delta.$$
We define
$$\xi_\delta \in W^{1,2}\big([0,\infty)\big)$$
by
$$\xi_\delta(s)=\beta_\delta(s)\eta(0), \quad s \in [0,\infty).$$
Using (\ref{ad1}) we compute
\begin{eqnarray*}
|\eta(0)|^2&=&-\int_0^\infty \big\langle \eta(0),\beta_\delta'(s)\eta(0)\big\rangle ds\\
&=&-\int_0^\infty \big\langle \eta(0),\partial_s \xi_\delta(s)\big\rangle ds\\
&=&-\int_0^\infty \big\langle \eta(s),\partial_s \xi_\delta(s)\big\rangle ds+\int_0^\infty \big\langle \eta(s)-\eta(0),\partial_s \xi_\delta(s)\big\rangle ds\\
&=&\int_0^\infty \big\langle A^T(s)\eta(s),\xi_\delta(s)\big\rangle ds+\int_0^\infty \big\langle \eta(s)-\eta(0),\partial_s \xi_\delta(s)\big\rangle ds\\
&=&\int_0^\delta \big\langle A^T(s)\eta(s),\beta_\delta(s) \eta(0) \big\rangle ds+\int_0^\delta \big\langle \eta(s)-\eta(0),\partial_s \beta_\delta(s) \eta(0)\big\rangle ds
\end{eqnarray*}
By (\ref{ad2}) we see that $\eta$ is continuous so that for every $\epsilon>0$ there exists $\delta>0$ with the property that
$$\int_0^\delta \big\langle A^T(s)\eta(s),\beta_\delta(s) \eta(0) \big\rangle ds+\int_0^\delta \big\langle \eta(s)-\eta(0),\partial_s \beta_\delta(s) \eta(0)\big\rangle ds<\epsilon$$
so that
$$|\eta(0)|^2<\epsilon.$$
This implies that
$$\eta(0)=0.$$
However, by (\ref{ad2}) we see that $\eta$ is a solution of a first order ODE which is uniquely determined by its initial condition so that
$$\eta(s)=0, \quad s \geq 0.$$
This proves that $D_A$ is surjective. The same arguments apply if $A \in \mathfrak{A}_-$ or $\mathfrak{A}_T$ for $T>0$. We have proved the following Proposition.
\begin{prop}\label{fre2}
If $A \in \mathfrak{A}_\pm$ or $A \in \mathfrak{A}_T$ for $T>0$, then $D_A$ is a surjective Fredholm operator.
\end{prop}

\subsection{Computation of the indices}\label{indcomp}

In this paragraph we proof Theorem~\ref{index} and Theorem~\ref{index2}.
\\ \\
\textbf{Proof of Theorem~\ref{index}: } By Proposition~\ref{fre1} we know that $D_A$ for $A \in \mathfrak{A}$ is a Fredholm operator and it remains to compute its index. Suppose that $A_r$ for $r \in [0,1]$ is a continuous family in $\mathfrak{A}$. Because the index is invariant under homotopies by the stability Theorem~\ref{stab1} it follows that the index
$\mathrm{ind}(D_{A_r})$ is independent of $r$. 
\\ \\
Since $A^\pm$ are symmetric and nondegenerate there exist $\Psi^\pm\in SO(n)$ such that
$$(\Psi^\pm)^TA^\pm \Psi^\pm=D^\pm$$
where $D^\pm$ are diagonal matrices of the form
$$D^\pm=\left(\begin{array}{ccc}
a^\pm_1 & \cdots & 0\\
\vdots & \ddots &\vdots\\
0 & \cdots & a^\pm_n
\end{array}\right).$$
such that
$$a_i^\pm \neq 0,\quad 1 \leq i \leq n,\qquad a_i^\pm \leq a_j^\pm,\quad i \leq j.$$
Since $SO(n)$ is connected we can choose continuous paths
$$\Phi^\pm \colon [0,1] \to SO(n)$$
satisfying
$$\Phi^\pm(0)=\mathrm{1}, \quad \Phi^\pm(1)=\Psi^\pm.$$
Choose further
$$B \in C^0\big(\mathbb{R} \times [0,1],\mathrm{End}(\mathbb{R}^n)\big)$$
such that
$$B(r,s)=\left\{\begin{array}{cc}
\Phi^+(r) & s \geq 1\\
\Phi^-(r) & s \leq -1
\end{array}\right.$$
as well as
$$B(0,s)=\mathbb{1}, \quad \forall\,\,s \in \mathbb{R}.$$
We set for $r \in [0,1]$ and $s \in \mathbb{R}$
$$A_r(s):=B(r,s)^T A(s) B(r,s).$$
Note that $A_r$ is a continuous family in $\mathfrak{A}$ such that 
$$A_0=A$$
and
$$\lim_{s \to \pm \infty}A_1(s)=(\Phi^\pm(1))^T A^\pm \Phi^\pm(1)=(\Psi^\pm)^T A^\pm \Psi^\pm=D^\pm.$$
That means to compute the index we can assume without loss of generality that
$$A^\pm=D^\pm,$$
i.e., the asymptotic matrices are diagonal with respect to the same basis, namely the standard one. After a further homotopy but now with fixed asymptotics we can even assume that for every $s \in \mathbb{R}$
$$A(s)=\left(\begin{array}{ccc}
a_1(s) & \cdots & 0\\
\vdots & \ddots &\vdots\\
0 & \cdots & a_n(s)
\end{array}\right).$$
where the diagonal entries are smooth functions $a_i \in C^\infty(\mathbb{R},\mathbb{R})$ satisfying
$$a_i(s)=a^\pm_i,\quad \pm s \geq 1.$$
We now compute the kernel of $D_A$, when $A$ is of this form. If $\xi \in \mathrm{ker}D_A$, then $\xi \in W^{1,2}(\mathbb{R},\mathbb{R}^n)$ is a solution of the ODE
$$\partial_s \xi=-A\xi.$$
Writing $\xi$ in components
$$\xi(s)=(\xi_1(s),\ldots \xi_n(s)) \in \mathbb{R}^n$$
each component individually is a solution of the ODE
$$\partial_s \xi(s)=-a_i(s)\xi_i(s).$$
Since $a_i(s)=a^+_i$ for $s \geq 1$ we obtain that
$$\xi_i(s)=\xi_i(1)e^{-a_i^+(s-1)}, \quad s \geq 1,$$ 
and using that $a_i(s)=a^-(s)$ for $s \leq -1$ we have that
$$\xi_i(s)=\xi_i(-1)e^{-a_i^+(s+1)}, \quad s \leq -1.$$ 
On the other hand because $\xi_i \in W^{1,2}(\mathbb{R},\mathbb{R})$, we see from this that $\xi_i$ has to vanish identically if $a_i^+<0$ or $a_i^->0$. Since the asymptotic eigenvalues are ordered increasingly it holds that
\begin{eqnarray*}
a_i^+>0\,\,&\Longleftrightarrow&\,\,i>\mu(A^+)\\
a_i^-<0\,\,&\Longleftrightarrow&\,\,i\leq \mu(A^+)
\end{eqnarray*}
so that
\begin{equation}\label{dike}
\mathrm{dim}(\mathrm{ker}D_A)=\max\big\{\mu(A^-)-\mu(A^+),0\big\}.
\end{equation}
To compute the cokernel we take advantage that the cokernel coincides with the kernel of the adjoint so that by the computation of the adjoint in Section~\ref{adjoint} we have
$$\mathrm{coker}D_A=\mathrm{ker}D_{-A^T}=\mathrm{ker}D_A.$$
Therefore we compute using (\ref{dike})
\begin{eqnarray}\label{dicoke}
\mathrm{dim}(\mathrm{coker} D_A)&=&\max\big\{\mu(-A^-)-\mu(-A^+),0\big\}\\ \nonumber
&=&\max\big\{n-\mu(A^-)-(n-\mu(A^+)),0\big\}\\ \nonumber
&=&\max\big\{\mu(A^+)-\mu(A^-),0\big\}.
\end{eqnarray}
Combining (\ref{dike}) and (\ref{dicoke}) we can now compute the index
\begin{eqnarray*}
\mathrm{ind}(D_A)&=&\mathrm{dim}(\mathrm{ker}D_A)-\mathrm{dim}(\mathrm{coker}D_A)\\
&=&\max\big\{\mu(A^-)-\mu(A^+),0\big\}-\max\big\{\mu(A^+)-\mu(A^-),0\big\}\\
&=&\mu(A^-)-\mu(A^+).
\end{eqnarray*}
This finishes the proof of the theorem. \hfill $\square$
\\ \\
\textbf{Proof of Theorem~\ref{index2}: } By Proposition~\ref{fre2} we know that $D_A$ is a surjective Fredholm operator. Therefore it suffices to compute the dimension of its kernel. This can be carried out as in the proof of Theorem~\ref{index}.
\hfill $\square$

\section{Exponential decay}

\subsection{Exponential convergence to critical points}

Exponential decay tells us that if a gradient flow line converges asymptotically to a critical point, then the convergence
is exponential for all derivatives. 

\begin{thm}\label{expo}
Suppose that $x \in C^\infty(\mathbb{R},M)$ is a gradient flow line which converges asymptotically to a Morse critical point, i.e,
$$\lim_{s \to \infty} x(s)=x^+ \in \mathrm{crit}f.$$
Choose local coordinates around the critical point $x^+$ on which the critical point lies at the origin. Then there exists a constant $c>0$ and for each $k \in \mathbb{N}_0$ there are constants $\mu_k$ such that there exists a time $T \in \mathbb{R}$ with the property that
$$||\partial_s^k x(s)|| \leq \mu_k e^{-cs}, \quad s \geq T,\,\,k \in \mathbb{N}_0.$$
An analogous statement holds at the negative asymptotics. 
\end{thm}
\textbf{Remark: } The precise value of the constants in Theorem~\ref{expo} depends on the choice of the local coordinates around the critical point. However, the fact that such constants exist is independent of the choice of local coordinates.

\subsection{Exponential decay in the Euclidean case}

We first discuss a special case of Theorem~\ref{expo}. By the Morse Lemma reviewed as Lemma~\ref{morselemma} there exist local coordinates $U \subset \mathbb{R}^n$ around the critical point $x^+$ such that in these coordinates $x^+$ is at the origin and the function $f$ can be written as
$$f(x)=f(0)+\frac{1}{2}\langle x, A x\rangle$$
where $\langle \cdot, \cdot \rangle$ is the standard Euclidean inner product on $\mathbb{R}^n$ and $A$ is a diagonal matrix
\begin{equation}\label{diag}
A=\left(\begin{array}{ccc}
a_1 & \cdots & 0\\
\vdots & \ddots &\vdots\\
0 & \cdots & a_n
\end{array}\right).
\end{equation}
In fact, the Morse Lemma allows us even to arrange that all eigenvalues $a_i$ of $A$ are either $1$ or $-1$, but we prefer not to assume that to give the reader a feeling how the constants in Theorem~\ref{expo} are related to the eigenvalues of the
metric Hessian. However, we assume that the eigenvalues are ordered, i.e.,
$$a_i \leq a_j, \quad i \leq j$$
so that if $\mu \in \mathbb{N}_0$ is the Morse index of $f$ at $0$ it holds that
\begin{equation}\label{order}
a_i<0, \quad i \leq \mu,\qquad a_i>0,\quad i>\mu.
\end{equation}
Our additional assumption is that the metric on $U$ is just the restriction of the standard Euclidean metric on $\mathbb{R}^n$ to $U$. Note that this assumption cannot be in general arranged by a change of coordinates. Indeed, if a metric has nonvanishing curvature such a transformation is clearly impossible. 
\\ \\
Under the Euclidean assumption the gradient of $f$ is simply given by
$$\nabla f(x)=Ax$$
and $A$ corresponds to the metric Hessian of $f$. The gradient flow equation is 
$$\partial_s x+Ax=0$$
and if we write $x=(x_1,\ldots,x_n)$ we obtain for each component the ODE
$$\partial_s x_i+a_i x=0, \quad 1 \leq i \leq n$$
whose solution is
$$x_i(s)=x_i(0)e^{-a_i s}.$$
Since $x$ converges to $0$ as $s$ goes to infinity we obtain from (\ref{order}) that
$$x_i \equiv 0, \quad 1 \leq i \leq \mu.$$
Note that 
$$\partial_s^k x_i(s)=(-1)^kx_i(0)a_i^k e^{-a_i s}$$
so that the assertion of Theorem~\ref{expo} holds with the constant $c$ given by the smallest positive eigenvalue $a_{\mu+1}$ of the metric Hessian at the critical point.

\subsection{The action-energy inequality}

If the metric in Morse coordinates is not Euclidean the gradient flow equation cannot be solved in general explicitly and the 
proof of exponential decay is much more involved. A way to prove exponential decay for general Riemannian metrics is taking advantage of an action-energy inequality. This procedure was for example used in \cite{albers-frauenfelder, gaio-salamon, ziltener}. In holomorphic curve theory the action-energy inequality is a consequence of the isoperimetric inequality. 
\\ \\
The action-energy inequality is local, therefore we assume that $U \subset \mathbb{R}^n$ is an open subset containing the origin on which we have a smooth function $f \colon U \to \mathbb{R}$ which at $0$ has a unique critical point, which is Morse. We further assume that the critical value of $f$ is zero, i.e.,
$$f(0)=0$$
which we can always achieve by adding a constant to the function which does not change the gradient. If further $g$ is a Riemannian metric on $U$ an action-energy inequality is an inequality of the form 
\begin{equation}\label{actionenergy}
|f(x)|\leq \kappa ||\nabla_g f(x)||_g^2, \quad x \in U
\end{equation}
where $\kappa>0$ is a constant. The reason for this terminology comes from the fact that for gradient flow lines
$||\partial_s x||^2=||\nabla f(x)||^2$ so that $||\nabla f(x)||^2$ can be interpreted as an energy density. 
\\ \\
We first explain that we always have an action-energy inequality in Morse coordinates for an Euclidean metric.
\begin{lemma}\label{acen1}
Suppose that the Morse function is given by the quadratic form
\begin{equation}\label{morseform}
f(x)=\langle x, Ax\rangle
\end{equation}
where $A$ is as in (\ref{diag}) and $g$ is the Euclidean metric. Then (\ref{actionenergy}) holds with the constant $\kappa$ given by
$$\kappa=\frac{1}{2\min\{|a_i|: 1 \leq i \leq n\}}=\frac{1}{2\min\{-a_\mu,a_{\mu+1}\}}.$$
\end{lemma}
\textbf{Proof: } We estimate
\begin{eqnarray*}
||\nabla f||^2&=&\sum_{i=1}^n a_i^2 x_i^2\geq \frac{1}{2\kappa}\sum_{i=1}^n |a_i| x_i^2\geq \frac{1}{\kappa}|f(x)|.
\end{eqnarray*}
This proves the lemma. \hfill $\square$
\\ \\
It is interesting to note that in the Euclidean case $\kappa$ can be expressed with the smallest absolute value of the eigenvalues of the metric Hessian at the critical point. Observe that eigenvalues of the metric Hessian depend on the metric.
\\ \\
We next estimate how the energy density changes under a change of the metric.
\begin{lemma}\label{acen2}
Suppose that $g_1$ and $g_2$ are two Riemannian metric on $U$ such that for a constant $c>0$ their norms are subject to the inequality
$$||\cdot||_{g_1} \leq c||\cdot||_{g_2}.$$
Then the norms of their corresponding gradients satisfy the following inequality
$$||\nabla_{g_2} f||_{g_2} \leq c||\nabla_{g_1}||_{g_1}.$$
\end{lemma} 
\textbf{Proof: } We estimate
\begin{eqnarray*}
0&\leq&g_1\Big(c^2 \nabla_{g_1} f-\nabla_{g_2}f,c^2\nabla_{g_1}f-\nabla_{g_2}f\Big)\\
&=&c^4 g_1\big(\nabla_{g_1}f,\nabla_{g_1}f\big)-2c g_1\big(\nabla_{g_1} f,\nabla_{g_2}f\big)+g_1\big(\nabla_{g_2}f,\nabla_{g_2}g\big)
\end{eqnarray*}
The second term we estimate using the definition of the gradient two times for both metrics $g_1$ and $g_2$
$$g_1\big(\nabla_{g_1}f,\nabla_{g_2}f\big)=df\big(\nabla_{g_2}f\big)=g_2\big(\nabla_{g_2}g,\nabla_{g_2}g\big).$$
The third term can be estimated using the assumption of the lemma
$$g_1\big(\nabla_{g_2}f,\nabla_{g_2}f\big) \leq c^2g_2\big(\nabla_{g_2}f,\nabla_{g_2}f\big).$$
Combining these three facts we obtain the estimate
\begin{eqnarray*}
0 &\leq& c^4 g_1\big(\nabla_{g_1}f,\nabla_{g_1}f\big)-2c^2 g_2\big(\nabla_{g_2}f,\nabla_{g_2}f\big)+
c^2 g_2\big(\nabla_{g_2}f,\nabla_{g_2}f\big)\\
&=&c^4 g_1\big(\nabla_{g_1}f,\nabla_{g_1}f\big)-c^2 g_2\big(\nabla_{g_2}f,\nabla_{g_2}f\big)
\end{eqnarray*}
implying
$$||\nabla_{g_2}f||_{g_2}^2 \leq c^2||\nabla_{g_1}f||^2_{g_1}.$$
This finishes the proof of the lemma. \hfill $\square$
\\ \\
The two lemmas imply the existence of an action-energy inequality, namely we have the following proposition.
\begin{prop}\label{acenex}
Maybe after shrinking the open neighbourhood $U$ of $0$ an action-energy inequality (\ref{actionenergy}) holds true. 
\end{prop}
\textbf{Proof: } By the Morse Lemma, i.e., Theorem~\ref{morselemma}, maybe after shrinking $U$ we can assume that 
$f$ is of the form (\ref{morseform}). It follows from Lemma~\ref{acen1} that an action-energy inequality holds on $U$
for the Euclidean metric. Given an arbitrary Riemannian metric we can assume maybe after shrinking $U$ further that
the metric is equivalent to the Euclidean metric. Now the action-energy inequality for the general Riemannian metric follows from the one for the Euclidean metric by Lemma~\ref{acen2}. This finishes the proof of the proposition. \hfill $\square$

\subsection{Proof of exponential decay}

In this paragraph we explain how exponential decay follows from an action-energy inequality. We first establish the following lemma.
\begin{lemma}\label{expolem1}
Suppose that $f$ and $g$ on $U$ satisfy the action-energy inequality (\ref{actionenergy}) on $U$ and $x \colon [s_-,s_+] \to U$ is a gradient flow line of $\nabla_g f$ satisfying $f(x(s))>0$ for every $s \in [s_-,s_+]$. Then the following inequality holds for the distance between $x(s_-)$ and $x(s_+)$
$$d_g\Big(x(s_-),x(s_+)\Big) \leq \frac{2}{\sqrt{\kappa}}\Big(\sqrt{f(x(s_-))}-\sqrt{f(x(s_+))}\Big)
\leq \frac{2}{\sqrt{\kappa}}\sqrt{f(x(s_-))}.$$
\end{lemma}
\textbf{Proof: } Note that by the assumption that $f(x)>0$ it follows that the gradient never vanishes since $f$ has the unique critical point $0$ on $U$ where $f$ vanishes. Therefore we estimate using the gradient flow equation and the action-energy inequality
\begin{eqnarray*}
d_g\Big(x(s_-),x(s_+)\Big) &\leq& \int_{s_-}^{s_+}||\partial_s x||_gds\\
&=&\int_{s_-}^{s_+}||\nabla_g f(x)||_g ds\\
&=&\int_{s_-}^{s_+}\frac{||\nabla_g f(x)||^2_g}{||\nabla_g f(x)||_g}ds\\
&\leq&\frac{1}{\sqrt{\kappa}}\int_{s_-}^{s_+}\frac{||\nabla_g f(x)||^2}{\sqrt{f(x)}}ds\\
&=&-\frac{1}{\sqrt{\kappa}}\int_{s_-}^{s_+}\frac{df(x)\partial_s x}{\sqrt{f(x)}}ds\\
&=&-\frac{1}{\sqrt{\kappa}}\int_{s_-}^{s_+}\frac{\partial_s f(x)}{\sqrt{f(x)}}ds\\
&=&-\frac{2}{\sqrt{\kappa}}\int_{s_-}^{s_+}\partial_s \sqrt{f(x)}ds\\
&=&\frac{2}{\sqrt{\kappa}}\Big(\sqrt{f(x(s_-))}-\sqrt{f(x(s_+))}\Big).
\end{eqnarray*}
This finishes the proof of the lemma. \hfill $\square$
\\ \\
Our next lemma establishes exponential decay of the action.
\begin{lemma}\label{expolem2}
Suppose that $f$ and $g$ on $U$ satisfy the action-energy inequality (\ref{actionenergy}) on $U$ and $x \colon [s_0,\infty) \to U$ is a gradient flow line of $\nabla_g f$ satisfying $f(x(s))> 0$ for every $s \geq s_0$. Then
\begin{equation}\label{expac}
f(x(s)) \leq f(x(s_0))e^{\frac{1}{\kappa}(s_0-s)}, \quad s \geq s_0.
\end{equation}
\end{lemma}
\textbf{Proof: } The action-energy inequality combined with the gradient flow equations gives rise to the following inequality
$$f(x(s)) \leq \ \kappa||\nabla f(x(s))||^2_g=-\kappa\frac{d}{ds}f(x(s)).$$
This implies
$$\frac{d}{ds} \ln f(x(s))=\frac{\frac{d}{ds}f(x(s))}{f(x(s))}\leq -\frac{1}{\kappa}.$$
Integrating this inequality we obtain
$$\ln f(x(s))=\int_{s_0}^s \frac{d}{ds}\ln f(x(s))ds+\ln f(x(s_0))\leq -\frac{s-s_0}{c}+\ln f(x(s_0)).$$
Taking exponents we get (\ref{expac}). \hfill $\square$
\\ \\
We are now ready to prove exponential decay.
\\ \\
\textbf{Proof of Theorem~\ref{expo}: } By Proposition~\ref{acenex} we can assume without loss of generality that there exists
a time $T$ such that $x(s) \in U$ for $s \geq T$, where on $U$ an action-energy inequality holds. Using Lemma~\ref{expolem1}
and Lemma~\ref{expolem2} we estimate for $T<s_1<s_2$
\begin{eqnarray*}
d_g\big(x(s_1),x(s_2)\big)&\leq&\frac{2}{\sqrt{\kappa}}\sqrt{f(x(s_1))}\\
&\leq&\frac{2}{\sqrt{\kappa}}\sqrt{f(x(s_0))}e^{\frac{1}{2\kappa}(s_0-s_1)}\\
&=&\frac{2}{\sqrt{\kappa}}\sqrt{f(x(s_0))}e^{\frac{s_0}{2\kappa}}e^{-\frac{1}{2\kappa}s_1}.
\end{eqnarray*}
Set
$$\mu_0:=\frac{2}{\sqrt{\kappa}}\sqrt{f(x(s_0))}e^{\frac{s_0}{2\kappa}}, \quad c:=\frac{1}{2\kappa}.$$
With this notation the above inequality becomes
$$d_g\big(x(s_1),x(s_2)\big) \leq \mu_0 e^{-cs_1}.$$
Note that the right-hand side does not depend on $s_2$. Because $s_2>s_1$ is arbitrary and $x(s)$ converges to $0$ as
$s$ goes to infinity, we obtain from this 
$$d_g\big(x(s_1),0\big) \leq \mu_0 e^{-cs_1}.$$
This shows that $x$ converges to the critical point exponentially. That the derivates converge to zero as well exponentially
follows the exponential convergence of $x$ by bootstrapping the gradient flow equation. This finishes the proof about
exponential decay. \hfill $\square$

\section{Transversality}\label{transversality}

\subsection{The Hilbert manifold of paths}

Assume that $M$ is a closed connected manifold and $x^-$ and $x^+$ are two points on $M$. In this paragraph we introduce the 
Hilbertmanifold
$$\mathcal{H}:=\mathcal{H}_{x^-,x^+}$$
of $W^{1,2}$-paths from $x^-$ to $x^+$. To define charts we choose $x \in C^\infty(\mathbb{R},M)$ for which there
exists $T>0$ with the property that 
$$\left\{\begin{array}{cc}
x(s)=x^+ & s \geq T\\
x(s)=x^- & s \leq -T.
\end{array}\right.$$
Pulling back the tangent bundle of $M$ with $x$ we obtain the bundle $x^*TM \to \mathbb{R}$ which we can trivialize since
$\mathbb{R}$ is contractible. Hence we choose a trivialization
$$\Phi \colon x^* TM \to \mathbb{R} \times \mathbb{R}^n$$
where $n$ is the dimension of $M$. Let 
$$V \subset x^* TM$$
be an open neighbourhood of the zero section which has the property that for a Riemannian metric $g$ on $M$ for all
$r \in \mathbb{R}$ the restriction of the exponential map
$$\mathrm{exp}_g \colon T_{x(r)}M \cap V \to M$$
is injective. Abbreviate
$$U_x:=\big\{\xi \in W^{1,2}(\mathbb{R},\mathbb{R}^n): (r,\xi(r)) \in \Phi(V),\,\,\forall\,\,r \in \mathbb{R}\big\}
\subset W^{1,2}(\mathbb{R},\mathbb{R}^n).$$
We then have charts 
$$\phi_x \colon U_x \to \mathcal{H}, \quad \xi \mapsto \mathrm{exp}_g(\Phi^{-1}(\xi)).$$
Sobolev theory tells us that chart transitions are smooth in $W^{1,2}(\mathbb{R},\mathbb{R}^n)$. 
\\ \\
If $x \in \mathcal{H}$ then the tangent space of $\mathcal{H}$ at $x$ is given by
$$T_x \mathcal{H}=W^{1,2}(\mathbb{R},x^* TM),$$
i.e., consists of $W^{1,2}$-vector fields along $x$. Apart from the tangent bundle $T\mathcal{H} \to \mathcal{H}$ we need 
further the $L^2$-bundle $\mathcal{E} \to \mathcal{H}$ whose fiber over $x \in \mathcal{H}$ is given by 
$$\mathcal{E}_x=L^2(\mathbb{R},x^* TM),$$
i.e., consists of $L^2$-vector fields along $x$. Note that $\mathcal{E}_x$ contains the tangent space $T_x \mathcal{H}$
as a dense subspace. 

\subsection{Gradient flow lines as the zero set of a section}

Recall that if $\pi \colon \mathcal{E} \to \mathcal{H}$ is a vector bundle then a section $s \colon \mathcal{H} \to \mathcal{E}$ is required to satisfy
$$\pi \circ s=\mathbb{1}_{\mathcal{H}},$$
i.e., it maps any point $x \in \mathcal{H}$ to a vector in the fiber $\mathcal{E}_x=\pi^{-1}(x)$. Assume that $f \colon M \to \mathbb{R}$ is a Morse function on a Riemannian manifolg $(M,g)$, $x^-$ and $x^+$ are critical points of $f$, $\mathcal{H}
=\mathcal{H}_{x^-,x^+}$ is the Hilbert manifold of $W^{1,2}$-paths from $x^-$ to $x^+$ and $\mathcal{E} \to \mathcal{H}$
is the $L^2$-bundle over $\mathcal{H}$. We define a section
\begin{equation}\label{flosec}
s \colon \mathcal{H} \to \mathcal{E}, \quad x \mapsto \partial_s x +\nabla f(x).
\end{equation}
\begin{figure}[h]
\begin{center}
 \includegraphics[scale=1]{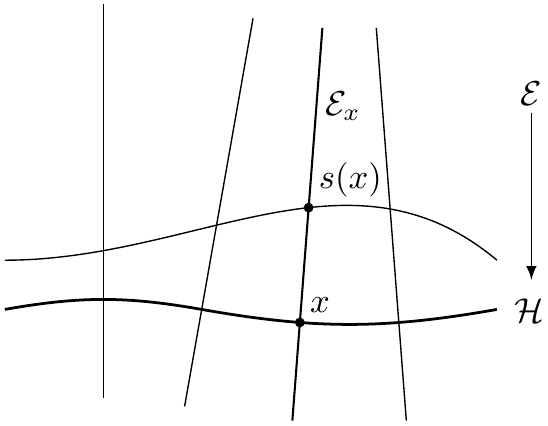}
 \caption{Section on a vector bundle.}
\end{center}
\end{figure}
Note that if $x \in s^{-1}(0)$, then $x$ is a $W^{1,2}$-gradient flow line from $x^-$ to $x^+$. In view of the gradient flow equation we can bootstrap $x$ to see that $x$ is actually a smooth gradient flow line. Moreover, in view of exponential decay 
established in Theorem~\ref{expo} every gradient flow line from $x^-$ to $x^+$ lies in $\mathcal{H}$ so that we can interpret the moduli space of parametrized gradient flow lines from $x^-$ to $x^+$ as the zero set of a section
$$\widetilde{\mathcal{M}}(f,g;x^-,x^+)=s^{-1}(0).$$

\subsection{The vertical differential}

Suppose that $\mathcal{E} \to \mathcal{H}$ is a vector bundle and $s \colon \mathcal{H} \to \mathcal{E}$ is a section.
If $x \in \mathcal{H}$, then the differential of the section is a linear map
$$ds(x) \colon T_x \mathcal{H} \to T_{s(x)}\mathcal{E}.$$
We identify $\mathcal{H}$ with the zero section of $\mathcal{E}$. Note that if $x \in \mathcal{H}$ we have a canonical splittling of the tangent space $T_x \mathcal{E}$ into horizontal and vertical subspaces
$$T_x \mathcal{E}=\mathcal{E}_x \oplus T_x \mathcal{H}.$$
It is worth pointing out that this splitting is only canonical if $x$ lies in the zero section. If $e$ is a general point on
$E$ we still have a canonical vertical subspace of the tangent space $T_e \mathcal{E}$ which we can canonically identify with
the fiber $\mathcal{E}_x$ but to define a horizontal subspace we need a connection.
\begin{figure}[h]
\begin{center}
 \includegraphics[scale=1]{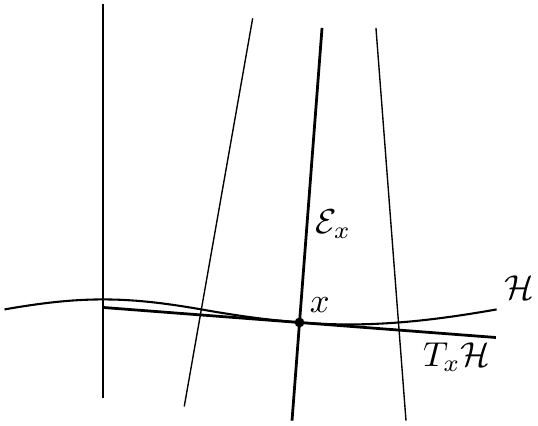} \quad
 \includegraphics[scale=1]{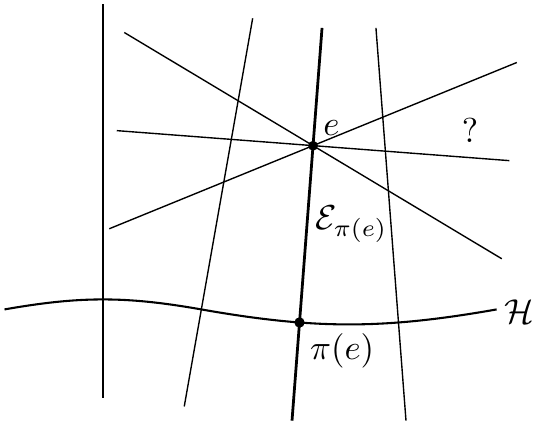}
 \caption{Canonical horizontal space at a point in the base (left) and ambiguity for general points in the bundle (right).}
\end{center}
\end{figure}
Suppose now that $x \in s^{-1}(0)$. We abbreviate
$$\pi \colon T_x\mathcal{E} \to \mathcal{E}_x$$
the projection along $T_x \mathcal{H}$ and define the \emph{vertical differential}
$$Ds(x):=\pi \circ ds(x) \colon T_x \mathcal{H} \to \mathcal{E}_x.$$
It is instructive to look at the case of a trivial vector bundle. Hence suppose that
$$\mathcal{E}=\mathcal{H} \times V$$
where $V$ is a vector space. If 
$$\sigma \colon \mathcal{H} \to V$$
is a smooth map we obtain a section
$$s \colon \mathcal{H} \to \mathcal{E}, \quad x \mapsto \big(x,\sigma(x)\big).$$
The differential of this section splits into
$$ds(x)=(\mathbb{1},d\sigma(x))$$
and the vertical differential becomes
$$Ds(x)=\sigma(x) \colon T_x \mathcal{H} \to V.$$
Using the vertical differential we can define the notion of a transverse section as well as the notion of a Fredholm section.  
\begin{fed}
A section $s \colon \mathcal{H} \to \mathcal{E}$ is called \emph{transverse to the zero section}, denoted by
$s \pitchfork 0$, if $Ds(x)$ is surjective for every $x \in s^{-1}(0)$. 
\end{fed}
\begin{fed}
A section $s \colon \mathcal{H} \to \mathcal{E}$ is called \emph{Fredholm section}, if
$Ds(x) \colon T_x \mathcal{H} \to \mathcal{E}_x$ is a Fredholm operator for every $x \in s^{-1}(0)$.
\end{fed}
The implicit function theorem, see for instance \cite[Appendix A]{mcduff-salamon} tells us the following.
\begin{thm}
Assume that $s \colon \mathcal{H} \to \mathcal{E}$ is a transverse Fredholm section, then
$s^{-1}(0)$ is a manifold and if $x \in s^{-1}(0)$, then the dimension of the component of $s^{-1}(0)$ containing $x$ is given by
\begin{equation}\label{indform}
\mathrm{dim}_x s^{-1}(0)=\mathrm{dim}(\mathrm{ker}Ds(x))=\mathrm{ind} Ds(x).
\end{equation}
\end{thm}
Note that the last equation in (\ref{indform}) follows from the definition of the Fredholm index since $Ds(x)$ is by the transversality assumption surjective and therefore its cokernel trivial. 
\\ \\
As a Corollary from Theorem~\ref{index} we have
\begin{cor}\label{dimunpar}
Suppose that $s \colon \mathcal{H}_{x^-,x^+} \to \mathcal{E}$ is given as in (\ref{flosec}) and $s \pitchfork 0$. 
Then $s^{-1}(0)=\widetilde{\mathcal{M}}(f,g;x^-,x^+)$ is a manifold and
$$\mathrm{dim}\big(\widetilde{\mathcal{M}}(f,g;x^-,x^+)\big)=\mu(x^-)-\mu(x^+).$$
\end{cor}
We point out here that the empty set is a manifold of arbitrary dimension. In particular, the Corollary does not give an
existence result for gradient flow lines.

\subsection{Generic metrics}\label{generic}

Suppose that $f \colon M \to \mathbb{R}$ is a Morse function on a closed manifold $M$ and $x^-$ and $x^+$ are critical points of $f$. The choice of a Riemannian metric $g$ on $M$ gives rise to the section as explained in (\ref{flosec})
$$s_g \colon \mathcal{H}_{x^-,x^+} \to \mathcal{E}, \quad x \mapsto \partial_s x+\nabla_g f(x).$$
If $s_g$ is not transverse to the zero section then we would like to say that we can perturb the metric a little bit to make the section transverse. To make this statement precise we introduce some language. Abbreviate by $\mathfrak{M}$ the space of
all smooth metric on $M$ endowed with the $C^\infty$-topology. We say that a subset $\mathfrak{U} \subset \mathfrak{M}$ is of
\emph{second category}, if it can be written as
$$\mathfrak{U}=\bigcap_{i=1}^\infty \mathfrak{M}_i$$
where $\mathfrak{M}_i \subset \mathfrak{M}$ is open and dense. It follows from Baire's Theorem that $\mathfrak{U}$ itself is dense in $\mathfrak{M}$. On the other hand not every dense subset is of second category. Indeed, a very useful feature of the property of being of second category is that this property is stable under countable intersection. Indeed, if
$\mathfrak{U}_n \subset \mathfrak{M}$ for $n \in \mathbb{N}$ is of second category than an immediate consequence of its definition is that
$$\bigcap_{n=1}^\infty \mathfrak{U}_n \subset \mathfrak{M}$$
is still of second category. This is definitely not true for any dense set - one just might think of the subsets of rational and irrational numbers inside the real numbers. 
\\ \\
Given a property of a metric one than says that this property is \emph{generic}, if the subset of metrics having this property is of second category. Because second category is stable under countable intersection, if one has countably many generic properties the property of having all this countably many properties together is still generic. 
\\ \\
The main result of this chapter is the following transversality result which tells us that for generic Riemannian metrics
the section $s_g$ is transverse to the zero section. 
\begin{thm}\label{transv}
There exists $\mathfrak{U} \subset \mathfrak{M}$ of second category such that
$$s_g \pitchfork 0, \quad \forall\,\,g \in \mathfrak{U}.$$
\end{thm}
The set $\mathfrak{U}$ might depend on the Morse function $f$ and the critical points $x^-$ and $x^+$, so that we write more precisely 
$$\mathfrak{U}(f;x^-,x^+) \subset \mathfrak{M}.$$
From Lemma~\ref{finite} we know that there are only finitely many critical points of $f$ and since the property of 
being of second category is stable under finite and even countable intersections we have that
$$\mathfrak{U}(f):=\bigcap_{(x^-,x^+) \in \mathrm{crit}f \times \mathrm{crit}f}\mathfrak{U}(f;x^-,x^+)$$
is still of second category in $\mathfrak{M}$. Metrics $g$ belonging to $\mathfrak{U}(f)$ have the property that
$s_g \colon \mathcal{H}_{x^-,x^+} \to \mathcal{E}$ are transverse to the zero section for all pairs of critical points
$(x^-,x^+)$ of $f$. Metrics with this property are referred to as \emph{Morse-Smale metrics} with respect to $f$. As a Corollary of Theorem~\ref{transv} we therefore have
\begin{cor}
There exists a set of second category $\mathfrak{U}(f) \subset \mathfrak{M}$ such that every $g \in \mathfrak{U}(f)$
is Morse-Smale with respect to the Morse function $f$. 
\end{cor}

\subsection{The universal moduli space}

In order to prove Theorem~\ref{transv} we consider the universal moduli space consisting of gradient flow lines for arbitrary metrics. The goal of this paragraph is to show that this universal moduli space is an (infinite dimensional) Banach manifold.
To prove Theorem~\ref{transv} we than consider the projection from the universal moduli space to the space of metrics. By Sard's theorem the set of regular values of this map is of second category and we show that the set of regular values gives us the desired set $\mathfrak{U}$. 
\\ \\
Unfortunately the space of all smooth metrics is not a Banach manifold but just a Fr\'echet manifold so that we cannot apply Sard's theorem directly. For this technical reason we consider therefore the space $\mathfrak{M}^k$ of all $C^k$-metrics on $M$ endowed with the $C^k$-topology. This space is a Banach manifold (even a cone in a Banach space). To define Morse homology we could actually use metrics of class $C^k$ for high enough regularity $k$ as well. On the other hand it is possible as well to obtain the smooth version of the transversality theorem from the $C^k$-version for all large enough integers $k$ due to an argument of Taubes discussed later. 
\\ \\
We consider the section
$$S \colon \mathfrak{M}^k \times \mathcal{H} \to \mathcal{E}, \quad (g,x) \mapsto s_g(x).$$
If $(g,x) \in S^{-1}(0)$ then $x$ is a gradient flow line from $x^-$ to $x^+$ with respect to $\nabla_g f$. 
The section $S$ has much more chance to be transverse to the zero section than the sections $s_g$ since we are now allowed to wiggle on the Riemannian metric as well to make the vertical differential surjective. The main result of this paragraph shows 
that transversality for this section is actually true. 
\\ \\
In the following we assume that $k$ is large so that we have enough regularity that the arguments in the following proofs go through. For example in the following proposition we should have $k$ to be at least one.
\begin{prop}\label{unitrans}
Assume that $(g,x) \in S^{-1}(0)$. Then $DS(g,x)$ is surjective. 
\end{prop} 
\textbf{Proof: } We consider two cases. 
\\ \\
\textbf{Case 1: } $x^- \neq x^+$.
\\ \\
In this case it is crucial to be able to wiggle on the Riemannian metric $g$. Suppose that $(g,x) \in S^{-1}(0)$, i.e.,
$$\partial_s x+\nabla_g f(x)=0.$$
A metric is symmetric and positive definite. The first condition is a closed condition while the second one is open. Therefore if $g \in \mathfrak{M}^k$ and $h \in T_g \mathfrak{M}^k$ this means that $h_y$ for every $y \in M$ is a symmetric bilinear form on $T_yM$ which depends $C^k$ on $y$. For $y \in M$ we introduce the linear map
$$L_{y,g} \colon \mathrm{Sym}(T_y M) \to T_y M, \quad h \mapsto \frac{d}{dr}\bigg|_{r=0} \nabla_{g+rh} f(y).$$
Using this notation we can write the vertical differential of $S$ at $(g,x)$ as
\begin{equation}\label{DS}
DS(g,x) \colon T_g \mathfrak{M}^k \oplus T_x \mathcal{H} \to \mathcal{E}_x, \quad
(h,\xi) \mapsto Ds_g(x)\xi+L_{x,g}h.
\end{equation}
Note in particular, that from this formula it follows that the image of $Ds_g(x)$ is contained in the image of
$DS(g,x)$ which implies that the image of $DS(g,x)$ is closed since $s_g(x)$ is Fredholm. Therefore to show surjectivity it suffices to prove that its orthogonal complement vanishes. We define the orthogonal complement with respect to the scalar product
$$\langle \eta_1,\eta_2 \rangle_g=\int_{-\infty}^\infty g_{x(s)}\big(\eta_1(s),\eta_2(s)\big)ds, \quad
\eta_1, \eta_2 \in \mathcal{E}_x.$$
Suppose that
$$\eta \in \mathrm{im} DS(g,x)^\perp=\mathrm{coker}(DS(g,x)).$$
Our goal is to show that $\eta$ vanishes. The assumption that $\eta$ lies in the orthogonal complement of the image of $DS(g,x)$ we can rephrase equivalently as
$$\big\langle DS(g,x)(h,\xi),\eta\big\rangle_g=0, \quad \forall\,\,(h,\xi) \in T_g \mathfrak{M}^k \oplus T_x \mathcal{H}.$$
In view of (\ref{DS}) this implies that
\begin{equation}\label{DS2}
\left\{\begin{array}{cc}
\langle Ds_g(x)\xi,\eta\rangle_g=0 & \forall\,\,\xi \in T_x \mathcal{H}\\
\langle L_{x,g}h,\eta\rangle_g=0 & \forall\,\,h \in T_g \mathfrak{M}^k. 
\end{array}\right.
\end{equation}
The first equation in (\ref{DS2}) implies that $\eta$ lies in the cokernel of $Ds_g(x)$. The cokernel of 
$Ds_g(x)$ coincides with the kernel of the adjoint of $Ds_g(x)$ and from the arguments in Section~\ref{adjoint} we infer that
$\eta$ is continuous. To take advantage of the second equation in (\ref{DS2}) we first note that by definition of the gradient we have
$$g(\nabla_g f(x),\eta)=df(x)\eta$$
and the righthand side is independent of the metric $g$. Hence differentiating with respect to the metric $g$ we obtain
$$0=\frac{d}{dr}\bigg|_{r=0}(g+rh)\big(\nabla_{g+rh}f(x),\eta\big)=h\big(\nabla_g f(x),\eta\big)+
g\big(L_{x,g} h,\eta\big).$$
Combining this with the second equation in (\ref{DS2}) we infer that for every $h \in T_g \mathfrak{M}^k$ 
\begin{equation}\label{DS3}
0=\big\langle L_{x,g} h,\eta\big \rangle_g=\int_{-\infty}^\infty g\big(L_{x,g} h,\eta\big)ds=
-\int_{-\infty}^\infty h\big(\nabla_g f(x),\eta\big)ds.
\end{equation}
We now assume by contradiction that there exists $s_0 \in \mathbb{R}$ with the property that
\begin{equation}\label{etas0}
\eta(s_0) \neq 0.
\end{equation}
To see how this contradicts (\ref{DS3}) we need the following lemma from linear algebra.
\begin{lemma}\label{linalg}
Suppose that $v,w \in \mathbb{R}^n \setminus \{0\}$. Then there exists a symmetric bilinear form such that
$$h(v,w)>0.$$
\end{lemma} 
\textbf{Proof: } Choose a complement $X$ of the one-dimensional subspace of $\mathbb{R}^n$ spanned by $v$ with the property that $w \notin X$. That implies that there exists $\lambda \neq 0$ such that
$$w=\lambda v \,\,\textrm{mod}\,\,X.$$
If $y_1,y_2 \in \mathbb{R}^n$ there exist unique $r_1,r_2 \in \mathbb{R}$ and $x_1,x_2 \in X$ such that
$$y_1=r_1v+x_1,  \qquad y_2=r_2 v+x_2.$$
We define
$$h(y_1,y_2):=\lambda r_1 r_2.$$
Then $h$ is symmetric and bilinear and moreover
$$h(v,w)=\lambda^2>0.$$
This finishes the proof of the Lemma. \hfill $\square$
\\ \\
\textbf{Proof of Proposition~\ref{unitrans} continued: } Since by assumption $x^- \neq x^+$ and gradient flow lines flow downhill by Lemma~\ref{downlem} we have
$$\nabla_g f(x(s)) \neq 0, \quad \forall\,\,s \in \mathbb{R}.$$
Using (\ref{etas0}) and Lemma~\ref{linalg} we conclude from this that there exists a symmetric bilinear form $h_0$ on
$T_{x(s_0)} M$ such that
$$h_0\big(\nabla_g f(x(s_0)),\eta(s_0)\big)>0.$$
Consider a chart $U$ around $x(s_0)$. In a chart all tangent spaces are naturally identified and $h_0$ canonically extends to all tangent spaces of points lying in the chart. Since we have already seen that $\eta$ is continuous, we infer that there 
exists $\epsilon>0$ such that 
$$h_0\big(\nabla_g f(x(s)),\eta(s)\big)>0, \quad s \in (s_0-\epsilon, s_0+\epsilon).$$
We set
$$V:=U \cap f^{-1}\big((f\circ x(s_0+\epsilon),f\circ x(s_0-\epsilon)\big).$$
Since $x$ is a nonconstant gradient flow line it holds that $f \circ x$ is strictly monotone decreasing and therefore
$$x^{-1}(V)=(s_0-\epsilon,s_0+\epsilon).$$
Choose a cutoff function $\beta \in C^\infty(M,[0,1])$ satisfying
$$\mathrm{supp}\beta \subset V, \quad \beta(x(s_0))=1$$
and define
$$h:=\beta \cdot h_0.$$
Then
$$h\big(\nabla_g f(x),\eta\big)(s) \geq 0,\quad \forall\,\,s \in \mathbb{R},\qquad
h\big(\nabla_g f(x),\eta\big)(s_0)>0.$$
In particular,
$$\int_{-\infty}^\infty h\big(\nabla_g f(s),\eta\big) ds>0.$$
This contradicts (\ref{DS3}) and therefore $\eta$ has to vanish identically. This finishes the proof of Case\,1. 
\\ \\
\textbf{Case\,2: } $x^-=x^+$.
\\ \\
In this case we do not need to wiggle at the Riemannian metric. Because gradient flow lines flow downhill by Lemma~\ref{downhill} the only gradient flow line from $x^-$ to $x^+$ is for all Riemannian metrics the constant one
$$x(s)=x^-=x^+, \quad \forall\,\,s \in \mathbb{R}.$$
In other words
$$s_g^{-1}(0)=\{x\}, \quad \forall\,\,g \in \mathfrak{M}^k.$$
We show that $Ds_g(x)$ is surjective for all Riemannian metrics. Since $x$ is constant the linearisation of the gradient flow equation is given by
$$D_A \xi=\partial_s \xi+A\xi$$
where $A$ is constant, namely the Hessian of $f$ and $x^+=x^-$. The Hessian is symmetric and moreover, since $f$ is Morse it
is nondegenerate as well and after changing coordinates we can assume that it is diagonal
$$A=\left(\begin{array}{ccc}
a_1 & \cdots & 0\\
\vdots & \ddots &\vdots\\
0 & \cdots & a_n
\end{array}\right)$$
such that
$$a_i \neq 0, \quad 1 \leq i \leq n.$$
Suppose that that $\xi=(\xi_1,\ldots,\xi_n)$ is in the kernel of $D_A$ i.e.,
$$\partial_s \xi+A\xi=0.$$
This implies that
$$\xi_i(s)=\xi_i(0)e^{-a_i s}, \quad 1 \leq i \leq n$$
On the other hand elements in the kernel of $D_A$ have to lie in $W^{1,2}$ and therefore the only element in the kernel is the trivial one so that we have
$$\ker D_A=\{0\}.$$
The same reasoning of course applies to $-A$ instead of $A$ and therefore
$$\mathrm{coker}D_A=\mathrm{ker}D_{-A}=\{0\}.$$
This implies that $D_A$ is surjective and the proposition is proven. \hfill $\square$

\subsection{Regular values}

As a consequence of Proposition~\ref{unitrans} and the implicit function theorem we have that the universal moduli space
$S^{-1}(0)$ is a Banach manifold. We now consider the map
$$\Pi \colon S^{-1}(0) \to \mathfrak{M}^k, \quad (g,x) \mapsto g.$$
Recall that $g$ is called a \emph{regular value} of $\Pi$ is for all $y \in \Pi^{-1}(g)$ the differential
$$d \Pi(y) \colon T_y S^{-1}(0) \to T_g \mathfrak{M}^k$$
is surjective. Since $S^{-1}(0)$ is a Banach manifold, Sard's theorem (see for instance \cite[Appendix A]{mcduff-salamon})
tells us that 
$$\mathfrak{M}^k_{\mathrm{reg}}=\big\{g \in \mathfrak{M}^k: g\,\,\textrm{regular value of}\,\,\Pi\big\}$$
is of second category in $\mathfrak{M}^k$. The next lemma tells us that the regular values of $\Pi$ are precisely the metrics
for which the section $s_g$ is transverse to the zero section. 
\begin{lemma}\label{salem}
If $g \in \mathfrak{M}^k_{\mathrm{reg}}$, then $Ds_g(x)$ is surjective for every $(g,x) \in \Pi^{-1}(g)$.
\end{lemma}
\textbf{Proof: } If $(g,x) \in S^{-1}(0)$ we can describe its tangent space using (\ref{DS}) by
\begin{eqnarray*}
T_{(g,x)}S^{-1}(0)&=&\big\{(h,\xi) \in T_g \mathfrak{M}^k \oplus T_x \mathcal{H}: DS(g,x)(h,\xi)=0\big\}\\
&=&\big\{(h,\xi) \in T_g \mathfrak{M}^k \oplus T_x \mathcal{H}: Ds_g(x)\xi=-L_{x,g}h\big\}.
\end{eqnarray*}
Choose $\eta \in \mathcal{E}_x$. By Proposition~\ref{unitrans} there exists 
$$(h_0,\xi_0) \in T_g \mathfrak{M}^k \oplus T_x \mathcal{H}$$
such that
\begin{equation}\label{sa1}
DS(g,x)(h_0,\xi_0)=\eta.
\end{equation}
Since $g \in \mathfrak{M}^k_{\mathrm{reg}}$ there exists
$$(h_1,\xi_1) \in T_{(g,x)}S^{-1}(0)$$
such that 
\begin{equation}\label{sa2}
d\Pi(g,x)(h_1,\xi_1)=h_0.
\end{equation}
Because $(h_1,\xi_1) \in T_{(g,x)} S^{-1}(0)$ it holds that
\begin{equation}\label{sa3}
Ds_g(x)\xi_1=-L_{x,g} h_1.
\end{equation}
Moreover, by definition of $\Pi$ we have
\begin{equation}\label{sa4}
d\Pi(g,x)(h_1,\xi_1)=h_1.
\end{equation}
Combining (\ref{sa2}) and (\ref{sa4}) we obtain
\begin{equation}\label{sa5}
h_0=h_1.
\end{equation}
Using (\ref{sa1}), (\ref{sa3}), and (\ref{sa5}) we compute
\begin{eqnarray*}
\eta&=&DS(g,x)(h_0,\xi_0)\\
&=& Ds_g(x)\xi_0+L_{x,g}h_0\\
&=&Ds_g(x)\xi_0+L_{x,g} h_1\\
&=&Ds_g(x)(\xi_0-\xi_1).
\end{eqnarray*}
Setting
$$\xi:=\xi_0-\xi_1 \in T_x \mathcal{H}$$
this becomes
$$Ds_g(x)\xi=\eta$$
showing that $Ds_g(x)$ is surjective. This finishes the proof of the lemma. \hfill $\square$
\\ \\
Combining Lemma~\ref{salem} with Sard's theorem we have the following Corollary
\begin{cor}\label{pretaubes}
For all integers $k$ large enough there is a subset $\mathfrak{U}^k \subset \mathfrak{M}^k$ of the second category such that
$$s_g \pitchfork 0, \quad \forall\,\,g \in \mathfrak{U}^k.$$
\end{cor}
This Corollary gives us almost the statement of Theorem~\ref{transv}, except that the metrics are not necessarily smooth but just of arbitrary high regularity. For practical purposes to define Morse homology this is enough, since one does not necessarily need smooth metrics but metrics of high regularity are enough. To get smooth metrics one needs an additional argument by Taubes explained in the following paragraph. Readers who are satisfied with only finitely many times differentiable metrics can skip this paragraph.

\subsection{An argument by Taubes}

In this paragraph we finish the proof of Theorem~\ref{transv}.
\\ \\
\textbf{Proof of Theorem~\ref{transv}: }For $g \in \mathfrak{M}^k$ we abbreviate by
$$\widetilde{\mathcal{M}}(g):=\widetilde{\mathcal{M}}(f,g;x^-,x^+)$$
the moduli space of all parametrised gradient flow lines of $\nabla_g f$ from $x^-$ to $x^+$. 
For $c>0$ we introduce
$$\widetilde{\mathcal{M}}^c(g):=\Big\{x \in \widetilde{\mathcal{M}}(g): ||\partial_s x(s)||_g \leq c e^{-\frac{|s|}{c}}, 
\,\,s \in \mathbb{R}\Big\}$$
the subspace of parametrised gradient flow lines satisfying a uniform exponential decay. In view of Theorem~\ref{expo} about exponential decay it holds that
\begin{equation}\label{expcor}
\widetilde{\mathcal{M}}(g)=\bigcup_{c \in \mathbb{N}} \widetilde{\mathcal{M}}^c(g).
\end{equation}
Gradient flow lines satisfying a uniform exponential decay cannot break and therefore for each $c>0$ the space
$\widetilde{\mathcal{M}}^c(g)$ is compact. For a large integer $k$ and $c>0$ abbreviate
$$\mathfrak{U}^{k,c}:=
\Big\{g \in \mathfrak{M}^k: Ds_g(x)\,\,\textrm{surjective},\,\,\forall\,\,x \in \widetilde{\mathcal{M}}^c(g)\Big\}.$$
Note that
$$\mathfrak{U}^k \subset \mathfrak{U}^{k,c}$$
so that by Corollary~\ref{pretaubes} the space $\mathfrak{U}^{k,c}$ is dense in $\mathfrak{M}$. On the other hand, since
$\widetilde{\mathcal{M}}^c(g)$ is compact it is open as well. Therefore
$$\mathfrak{U}^{\infty,c}:=\mathfrak{U}^{k,c} \cap \mathfrak{M}$$
is dense in $\mathfrak{M}^k$ with respect to the $C^k$-topology for every large integer $k$. In particular, it is dense with respect to the $C^\infty$-topology. Taking once more advantage of the fact that the spaces $\widetilde{\mathcal{M}}^c(g)$
are compact it follows that $\mathfrak{U}^{\infty,c}$ is as well $C^\infty$-open. We have shown that
$$\mathfrak{U}:=\bigcap_{c \in \mathbb{N}}\mathfrak{U}^{\infty,c}$$
is a countable intersection of open and dense sets and therefore of second category. From (\ref{expcor}) we infer that if
$g \in \mathfrak{U}$ it follows that $Ds_g(x)$ is surjective for every $x \in \widetilde{\mathcal{M}}(c)$ so that
$s_g$ is transverse to the zero section. This finishes the proof of the theorem. \hfill $\square$

\subsection{Morse-Smale metrics}

Recall that if $f \colon M \to \mathbb{R}$ is a Morse function on a closed manifold $M$, then a Morse-Smale metric
with respect $f$ is a Riemannian metric $g$ such that for every pair of critical points $(x^-,x^+)$ of $f$ the section
$s_g \colon \mathcal{H}_{x^-,x^+}$ is transverse to the zero section. In particular, by Corollary~\ref{dimunpar}
the moduli spaces of parametrised gradient $\widetilde{\mathcal{M}}(f,g;x^-,x^+)$ are manifolds of dimension
\begin{equation}\label{dimun}
\mathrm{dim}\widetilde{\mathcal{M}}(f,g;x^-,x^+)=\mu(x^-)-\mu(x^+).
\end{equation}
In this paragraph we establish some structural results about the moduli spaces of unparametrised gradient flow lines and its completions for Morse-Smale metrics.
\begin{thm}\label{dimpar}
Assume that $g$ is Morse-Smale and $x^-$ and $x^+$ are two different critical points of $f$. Then the moduli space of
unparametrised gradient flow lines $\mathcal{M}(f,g;x^-,x^+)$ is a manifold of dimension
$$\mathrm{dim}\mathcal{M}(f,g;x^-,x^+)=\mu(x^-)-\mu(x^+)-1.$$
\end{thm} 
\textbf{Proof: } Since $x^- \neq x^+$ the action of $\mathbb{R}$ on $\widetilde{\mathcal{M}}:=\widetilde{\mathcal{M}}(f,g;x^-,x^+)$ by timeshift is free. We check that it is proper as well. This means that 
$$\big\{(x,r_*x): x \in \widetilde{\mathcal{M}},\,\,r \in \mathbb{R}\big\} \subset \widetilde{\mathcal{M}} \times
\widetilde{\mathcal{M}}$$
is closed. Suppose that $x_\nu \in \widetilde{\mathcal{M}}$ and $r_\nu \in \mathbb{R}$ are sequences for which there exist
$x,y \in \widetilde{\mathcal{M}}$ such that
$$\lim_{\nu \to \infty}\big(x_\nu, (r_\nu)_*x_\nu\big)=(x,y).$$
Since $x$ and $y$ have the same asymptotics it follows from Theorem~\ref{unilim} that they are equivalent, i.e., that there exists $r \in \mathbb{R}$ such that
$$y=r_* x.$$
This proves that the $\mathbb{R}$-action is proper. 
\\ \\
If a Lie group $G$ acts freely and properly on a manifold $N$ its quotient space $N/G$
is a manifold, see for instance \cite[Theorem 3.5.25]{abraham-marsden-ratiu}, and its dimension is given by
$$\mathrm{dim} (N/G)=\mathrm{dim}(N)-\mathrm{dim}(G).$$
Applying this to $\mathcal{M}=\widetilde{\mathcal{M}}/\mathbb{R}$ the assertion of the theorem follows. \hfill $\square$
\\ \\
Recall from Section~\ref{moduli} that if $x^-$ and $x^+$ are two different critical points the space of unparametrised gradient flow lines can be compactified to
$$\overline{\mathcal{M}}(f,g;x^-,x^+)=\bigsqcup_{k=1}^\infty \mathcal{M}^b_k(f,g;x^-,x^+)$$
where $\mathcal{M}^b_k(f,g;x^-,x^+)$ is the moduli space of unparametrised $(k-1)$-fold broken gradient flow lines from
$x^-$ to $x^+$. In particular
$$\mathcal{M}(f,g;x^-,x^+)=\mathcal{M}^b_1(f,g;x^-,x^+).$$
When the metric is Morse-Smale we obtain additional information about this stratification.
\begin{cor}
Assume that $g$ is Morse-Smale and $x^-$ and $x^+$ are two critical points of $f$ such that
$$\ell:=\mu(x^-)-\mu(x^+) \in \mathbb{N}.$$
Then
$$\mathcal{M}^b_k(f,g;x^-,x^+)=\emptyset, \quad k >\ell.$$
\end{cor}
\textbf{Proof: } Suppose by contradiction that $x=(x_1,\ldots,x_k) \in \overline{\mathcal{M}}(f,g;x^-,x^+)$ for
$k>\ell$. Then there exists $j \in \{1,\ldots, k\}$ such that the asymptotics
$$x^\pm_j:=\lim_{s \to \pm \infty}x_j(s)$$
satisfy
$$\mu(x^-_j) \leq \mu(x^+_j).$$
But this is impossible for a Morse-Smale metric by Theorem~\ref{dimpar} and the Corollary follows. \hfill $\square$
\\ \\
Under the assumption of the Corollary we can write
$$\overline{\mathcal{M}}(f,g;x^-,x^+)=\bigsqcup_{k=1}^\ell \mathcal{M}^b_k(f,g;x^-,x^+).$$
Especially, if $\ell=1$ we have
$$\overline{\mathcal{M}}(f,g;x^-,x^+)=\mathcal{M}(f,g;x^-,x^+).$$
Therefore, since $\overline{\mathcal{M}}(f,g;x^-,x^+)$ is compact we obtain as a further Corollary
\begin{cor}\label{dd2f}
Assume that $g$ is Morse-Smale and $x^-$ and $x^+$ are critical points of $f$ such that
$$\mu(x^-)-\mu(x^+)=1.$$
Then $\mathcal{M}(f,g;x^-,x^+)$ is a compact zero-dimensional manifold, i.e., a finite set. 
\end{cor}

\section{Stable and unstable manifolds}\label{stableunstable}
If $M$ is a closed manifold and $q \in M$ we denote by $\mathcal{H}_q^-$ the Hilbert manifold of all paths $x \in W^{1,2}\big((-\infty,0],M\big)$ satisfying $\lim_{s\to -\infty}x(s)=q$ and similarly by $\mathcal{H}_q^+$ the Hilbert manifold of paths
$x \in W^{1,2}\big([0,\infty),M\big)$ satisfying $\lim_{s\to \infty}x(s)=q$. We have a smooth evaluation map
$$\mathrm{ev}\colon \mathcal{H}^\pm_q  \to M, \quad x \mapsto x(0).$$
We further denote by
$$\mathcal{E} \to \mathcal{H}^\pm_q$$
the $L^2$-bundle whose fiber over a path $x \in \mathcal{H}^-_q$ is given by $L^2$-vector fields along $x$
$$\mathcal{E}_x=L^2\big((-\infty,0], x^* TM\big)$$
and similarly for $x \in \mathcal{H}^+_q$
$$\mathcal{E}_x=L^2\big([0,\infty), x^* TM\big).$$ 
If $f \colon M \to \mathbb{R}$ is a Morse function on $M$ and $g$ a Riemannian metric on $M$ giving rise to a gradient
$\nabla f$ and $c \in M$ is a critical point of  $f$ we define sections
$$\sigma^\pm_c \colon \mathcal{H}^\pm_c \to \mathcal{E}, \quad x \mapsto \partial_s x+\nabla f(x).$$
We abbreviate its zero set by
$$W^\pm_c:=(\sigma^\pm_c)^{-1}(0)$$
and refer to $W^-_c$ as the \emph{unstable manifold} of the critical point $c$ and $W^+_c$ as its \emph{stable manifold}.
In view of exponential decay established in Theorem~\ref{expo} elements of the unstable manifold are all partial gradient flow lines, i.e., 
$$x \in C^\infty\big((-\infty,0],M)$$
solving the ODE
$$\partial_s x(s)+\nabla f(x(s)), \quad s \in (-\infty,0],$$
which meet the asymptotic condition
$$\lim_{s \to -\infty}x(s)=c.$$
Similarly the stable manifold consists of partial gradient flow lines 
$$x \in C^\infty\big([0,\infty),M\big)$$
whose positive asymptotics is $c$. Because a solution of an ODE is uniquely determined by its initial value the restriction of the evaluation map
$$\mathrm{ev}|_{W^{\pm}_c} \colon W^{\pm}_c \to M$$
is injective. We denote by $\phi^s_{\nabla f}$ the flow of the gradient vector field of $f$, i.e.,\,for $s \in \mathbb{R}$
we have a smooth family of diffeomorphisms of $M$ determined by the conditions
$$\phi^0_{\nabla f}=\mathrm{id}_M, \quad \frac{d}{ds}\phi^s_{\nabla f}=\nabla f \circ \phi^s_{\nabla f}.$$
With this notation the image of the unstable manifold under the evaluation map can be described as
$$\mathrm{ev}(W^-_c)=\big\{q \in M: \lim_{s \to -\infty}\phi^s_{\nabla f}(q)=c\big\}$$ 
and similarly for the image of the stable manifold
$$\mathrm{ev}(W^+_c)=\big\{q \in M: \lim_{s \to \infty}\phi^s_{\nabla f}(q)=c\big\}.$$
In fact, it is actually more common to define stable and unstable manifolds directly by the expressions on the righthand side and interpret them as subsets of the manifold $M$. Our approach of thinking of them as subsets of path spaces is 
motivated by the work of Simcevic \cite{simcevic}. In this work Simcevic showed that stable and unstable manifolds can 
as well be defined for action functionals living on infinite dimensional spaces which appear in Floer homology. In this case however, the evaluation maps are much more complicated and have targets in suitable interpolation spaces. Thinking
of stable and unstable manifolds as subsets of path spaces also gives rise to a quick proof that they are actually manifolds by applying the powerful infinite dimensional implicit function theorem. Namely the following theorem holds.
\begin{figure}[h]
\begin{center}
 \includegraphics[scale=1]{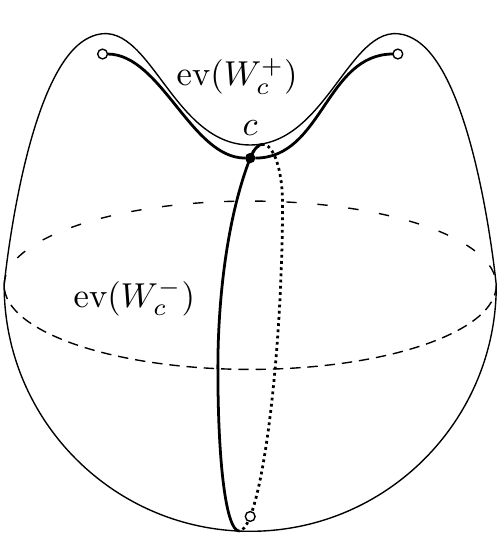}
 \caption{Stable and unstable manifolds on the distorted sphere.}
\end{center}
\end{figure}
\begin{thm}\label{stabthm}
The stable and unstable manifolds are smooth manifolds and their dimension is given by
$$\mathrm{dim}W^-_c=\mu(c), \quad \mathrm{dim}W^+_c=n-\mu(c)$$
where $\mu(c)$ is the Morse index of the critical point $c$ and $n$ is the dimension of the manifold $M$. 
\end{thm} 
\textbf{Proof: }By Theorem~\ref{index2} we know that for $x \in W^\pm_c$ the vertical differential 
$$D \sigma_c^\pm(x) \colon T_x \mathcal{H}_c^\pm \to \mathcal{E}_x$$
is a surjective Fredholm operator. If $x$ lies in the unstable manifold $W^-_c$ its index is given by  
$$\mathrm{ind}\big(D \sigma^-_c(x)\big)=\mu(c),$$
respectively, if it lies in the stable manifold $W^+_c$ the index satisfies
$$\mathrm{ind}\big(D \sigma^+_c(x)\big)=n-\mu(c).$$
Hence the theorem is a consequence of the Implicit Function Theorem. \hfill $\square$
\\ \\
Note that a special element in $W^\pm_c$ is the critical point $c$ itself interpreted as a 
constant partial gradient flow line. We denote by 
$$c_\pm \in W^\pm_c$$ 
the constant map to $c$, which serves as a natural base point in $W^\pm_c$.
The two manifolds $W^\pm_c$ naturally retract to $c_\pm$.
To see that we observe that one can endow them with an action of the monoid $\mathbb{R}^+$. Namely $\mathbb{R}^+$ as a set is the nonnegative real halfline $[0,\infty)$ with operation given by addition. The monoid $\mathbb{R}^+$ acts on $\mathcal{H}^\pm_c$ by timeshift, namely if $x \in \mathcal{H}^-_c$ and $r \in \mathbb{R}^+$ the action is by negative timeshift
$$r_*x(s)=x(s-r), \quad s \in (-\infty,0]$$
and if $x \in \mathcal{H}^+_c$ the action is by positive timeshift
$$r_*x(s)=x(s+r), \quad s \in [0,\infty).$$
Note that $c$ is the only fixed point of the action of $\mathbb{R}^+$ on $W^\pm_c$ and on the complement of $c$ the action
is free. 
Observe that for every $x \in W^\pm_c$ we have
$$\lim_{r \to \infty}r_* x=c_\pm.$$
In particular we have
\begin{lemma}
Stable and unstable manifolds are contractible.
\end{lemma}
\textbf{Remark: } One can show that the stable and unstable manifolds are diffeomorphic to Euclidean vector space of the corresponding dimension, see for instance \cite[Chapter II.3]{zehnder}.
\begin{thm}\label{embedding}
The evaluation maps $\mathrm{ev}|_{W^\pm_c} \colon W^\pm_c \to M$ are smooth embeddings. 
\end{thm}
\textbf{Proof: } That the evaluation maps are injective follows from the fact that an ODE is uniquely determined by its initial value. To see that they are immersions we note, that if $x \in W^\pm_c$ then its tangent space is given by
$$T_x W^\pm_c=\mathrm{ker}D \sigma^\pm_c$$
which means that their elements are itself solutions of a first order linear ODE. The differential of the evaluation map is given by 
$$d \mathrm{ev}(x) \colon T_x W^\pm_c \to T_{\mathrm{ev}(x)}M=T_{x(0)}M, \quad \xi \mapsto \xi(0).$$
Because $\xi$ is a solution of a first order linear ODE if the initial condition $\xi(0)$ vanishes, it follows that $\xi$
vanishes identically. This implies that the evaluation maps are immersions. 
\\ \\
Having established that the evaluation maps are injective immersions it remains to show that if $x_\nu$ is a sequence in
$W^\pm_c$ such that there exists $x \in W^\pm_c$ with the property that
\begin{equation}\label{conhaben}
\lim_{\nu \to \infty}\mathrm{ev}(x_\nu)=\mathrm{ev}(x)
\end{equation}
it follows that 
\begin{equation}\label{conmochten}
\lim_{\nu \to \infty}x_\nu=x.
\end{equation}
We discuss this for the stable case $W^+_c$. The unstable case can be discussed similarly and follows as well from the stable case by replacing $f$ with $-f$. Because we already know that the evaluation map is an injective immersion we can find
an open neighbourhood $U$ of the constant trajectory $c$ such that $\mathrm{ev}|_U \colon U \to M$ is an embedding. We choose an even smaller neighbourhood $U_0$ of $c$ with the property that its closure is contained in $U$, i.e.,
$$\overline{U}_0 \subset U.$$
Since on $U$ the evaluation map is already an embedding we can assume without loss of generality that
$$x_\nu \notin U_0,\quad \forall\,\,\nu \in \mathbb{N}.$$
We use the
$\mathbb{R}^+$ action on $W^+_c$. Namely there exist  nonnegative real numbers $r_\nu$ for every $\nu \in \mathbb{N}$ with the property that
$$(r_\nu)_*x_\nu \in \partial U_0,\quad \forall\,\,\nu \in \mathbb{N}.$$
Because $\overline{U}_0$ is compact and hence as well $\partial U_0$ we can assume maybe after transition to a subsequence that there exists 
$$y \in \partial U_0$$
with the property that
\begin{equation}\label{con1}
\lim_{\nu \to \infty}(r_\nu)_*x_\nu=y.
\end{equation}
We first consider the case where the sequence $r_\nu$ is bounded. In this case we can maybe after transition to a further subsequence assume that there exists $r_\infty \in \mathbb{R}^+$ such that
$$\lim_{\nu \to \infty}r_\nu=r_\infty.$$
From (\ref{con1}) we infer that
$$\lim_{\nu \to \infty}(r_\infty)_*x_\nu=y$$
and therefore
$$\mathrm{ev}(y)=\mathrm{ev}\big(\lim_{\nu \to \infty}(r_\infty)_*x_\nu\big)=\mathrm{ev}\big((r_\infty)_* x\big).$$
Because the evaluation map on $U$ is an embedding we deduce from this that
$$(r_\infty)_* x=y=\lim_{\nu \to \infty}(r_\infty)_*x_\nu$$
which implies (\ref{conmochten}).
\\ \\
To finish the proof we exclude the case that the sequence $r_\nu$ is not bounded. We argue by contradiction and assume that it is not bounded. In this case maybe after transition to a subsequence we can assume that
$$\lim_{\nu \to \infty}r_\nu=\infty.$$
Since $y$ lies on the boundary of $U_0$ it differs in particular from $c$ and using that $f \circ \mathrm{ev}$ on 
the stable manifold $W^+_c$ attains its unique global minimum at $c$ we conclude that 
$$f(y(0))=f \circ \mathrm{ev}(y)>f(c).$$
Abbreviate
$$\Delta:=f(y(0))-f(c)>0.$$
We choose open neighbourhoods $V_c \subset M$ of c satisfying
$$V_c \subset f^{-1}\Big(\big(-\infty,f(c)+\tfrac{\Delta}{2}\big)\Big).$$
This requirement guarantees that $y(0) \notin \overline{V}_c$, and moreover by (\ref{con1}) and the fact that gradient flow lines flow downhill there exists $\nu_0$ such that for every $\nu \geq \nu_0$ 
\begin{equation}\label{nichtdrin}
x_\nu(r) \notin V_c, \quad r \in [0,r_\nu].
\end{equation}
On the other hand since $x \in W^+_c$ there exists $R>0$ with the property that
$$x(r) \in V_c \quad r \geq R.$$
By (\ref{conhaben}) it holds that
$$\lim_{\nu \to \infty}x_\nu(0)=x(0)$$
and since gradient flow lines as solutions of an ODE depend continuously on their initial condition we conclude from that
$$\lim_{\nu \to \infty}x_\nu(R)=x(R).$$
In particular, there exists $\nu_1 \in \mathbb{N}$ such that
$$x_\nu(R) \in V_c, \quad \nu \geq \nu_1.$$
But since the sequence $r_\nu$ goes to infinity this contradicts (\ref{nichtdrin}) and the Theorem is proven. \hfill $\square$

\section{Gluing}

\subsection{The pregluing map}

Suppose that $U \subset \mathbb{R}^n$ is an open subset which is starshaped with respect to the origin, i.e., if $u \in U$, then the whole line segment from the origin to $u$ is contained in $U$, in formulas
$$u \in U \quad \Longrightarrow \quad \lambda u \in U,\,\,\forall\,\,\lambda \in [0,1].$$
Following Section~\ref{stableunstable} we abbreviate by 
$$\mathcal{H}_0^-:=\mathcal{H}_0^-(U)$$ 
the Hilbert manifold of all paths $x \in W^{1,2}\big((-\infty,0],U\big)$ satisfying $\lim_{s\to -\infty}x(s)=0$ and by 
$$\mathcal{H}_0^+:=\mathcal{H}^+_0(U)$$ 
the Hilbert manifold of paths
$x \in W^{1,2}\big([0,\infty),U\big)$ satisfying $\lim_{s\to \infty}x(s)=0$. Moreover, for $T>2$ we set
$$\mathcal{H}^T:=\mathcal{H}^T(U):=W^{1,2}\big([-T,T],U\big).$$
We further choose a smooth monotone cutoff function $\beta \in C^\infty\big(\mathbb{R},[0,1]\big)$ satisfying
$$\beta(s)=\left\{\begin{array}{cc}
0 & s \leq 0,\\
1 & s \geq 1.
\end{array}\right.$$
We define the \emph{pregluing map}
$$\#_0^T:=\#_{0,\beta}^T \colon \mathcal{H}^+_0 \times \mathcal{H}^-_0 \to \mathcal{H}^T$$
for $x \in \mathcal{H}^+_0$ and $y \in \mathcal{H}^-_0$ by the formula
\begin{equation}\label{pregl}
\#_0^T(x,y)(s)=
\Big(1-\beta\big(s+\tfrac{T}{2}+1\big)\Big)x(s+T)+
\beta\big(s-\tfrac{T}{2}\big)y(s-T), \quad s \in [-T,T].
\end{equation}
Observe that
$$\#_0^T(x,y)(s)=0, \quad s \in \big[-\tfrac{T}{2},\tfrac{T}{2}\big].$$
We further have
$$\#_0^T(x,y)(s)=x(s), \quad s \in \big[-T,-\tfrac{T}{2}-1\big]$$
and
$$\#_0^T(x,y)(s)=y(s), \quad s \in \big[\tfrac{T}{2}+1,T\big].$$
Moreover, for all negative times $\#_0^T(x,y)(s)$ lies in the line segment from the origin to $x(s)$ and for all
positive times $\#_0^T(x,y)(s)$ lies in the line segment from the origin to $y(s)$.
In particular, due to the assumption that $U$ is starshaped with respect to the origin the pregluing map is well defined.
\\ \\
We have evaluation maps
\begin{equation}\label{eva1}
\mathrm{ev} \colon \mathcal{H}^+_0 \times \mathcal{H}^-_0 \to U \times U, \quad (x,y) \mapsto \big(x(0),y(0)\big)
\end{equation}
and 
\begin{equation}\label{eva2}
\mathrm{ev} \colon \mathcal{H}^T \to U \times U, \quad z \mapsto \big(z(-T),z(T)\big).
\end{equation}
Note that
\begin{equation}\label{evainv}
\mathrm{ev}=\mathrm{ev} \circ \#_0^T.
\end{equation}
We further abbreviate by
$$0^+ \in \mathcal{H}^+_0, \quad 0^- \in \mathcal{H}^-_0, \quad 0^T \in \mathcal{H}^T$$
the constant maps to $0$ in the different path spaces. Then the following formula holds
\begin{equation}\label{nuller}
\#_0^T(0^+,0^-)=0^T.
\end{equation}

\subsection{Infinitesimal gluing}

We continue the notation of the previous subsection but assume in addition that we have given a Morse function
$f \colon U \to \mathbb{R}$ which has a critical point at the origin and a Riemannian metric $g$ on $U$. We further assume that the Riemannian metric $g$ at the origin is given by the standard inner product on $\mathbb{R}^n$ and the metric Hessian of $f$ at the origin is diagonal, i.e.,
\begin{equation}\label{Hes}
H:=H^g_f(0)=\left(\begin{array}{ccc}
a_1 & \cdots & 0\\
\vdots & \ddots &\vdots\\
0 & \cdots & a_n
\end{array}\right),
\end{equation}
where in addition the eigenvalues are increasing with increasing index, i.e.,
$$a_1 \leq a_2 \leq \cdots \leq a_n.$$
Because for a symmetric matrix there always exists an orthogonal basis for which it is diagonal we can always arrange $H$ to be of this form maybe after rotating $U$ which does not violate its property to be starshaped. Note that if
$$\mu=\mu_f(0) \in \{0,\ldots,n\}$$
is the Morse index of $f$ at $0$ we have
$$a_i<0,\quad i \in \{1,\ldots,\mu\},\qquad a_i>0,\quad i \in \{\mu+1,\ldots, n\}.$$
Let 
$$W^\pm_0=(\sigma_0^\pm)^{-1}(0)$$
be the stable and unstable manifolds of $0$ as explained in Section~\ref{stableunstable}. If $\mathcal{E}$ is the $L^2$-bundle we define analogousy sections
$$\sigma^T \colon \mathcal{H}^T \to \mathcal{E}, \quad x \mapsto \partial_s x+\nabla f(x)$$
and we denote its zero set by 
$$W^T:=(\sigma^T)^{-1}(0).$$
Elements of $W^T$ are segments of gradient flow lines defined on the interval $[-T,T]$. By Theorem~\ref{index2}
the vertical differential of $\sigma^T$ at $x \in W^T$ has index $n$ and is always surjective. In particular, 
$W^T$ is a smooth manifold satisfying
$$\mathrm{dim}W^T=n.$$
This can be seen alternatively as follows. Abbreviate by
$$\mathrm{ev}_\pm \colon \mathcal{H}^T \to U, \quad x \mapsto x(\pm T)$$
the evaluation maps at the starting point and end point. Because the gradient flow equation is a first order ODE a gradient flow line is uniquely determined by an initial condition and therefore both evaluation maps $\mathrm{ev}_\pm$ can be used
to identify an open neighbourhood of a point in $W^T$ with an open subset of $U$. Because $U$ is open, the gradient flow does not need to exist for all time. However, if we considered $W^T$ on a closed manifold $M$ the manifold $W^T$ would be
diffeomorphic to $M$ and both evaluation maps $\mathrm{ev}_\pm$ provide a diffeomorphism. 
\\ \\
We abbreviate
$$W^\infty:=W^+_0 \times W^-_0.$$
By Theorem~\ref{stabthm} we have
$$\mathrm{dim}W^-_0=\mu, \quad \mathrm{dim}W^+_0=n-\mu$$
so that
$$\mathrm{dim}(W^\infty)=n=\mathrm{dim}W^T.$$
Note that 
$$0_\infty:=(0_+,0_-) \in W^\infty, \quad 0_T \in W^T$$
so that by $(\ref{nuller})$ the pregluing map maps one point of $W^\infty$ to a point of $W^T$. On the other hand
the differential
$$d \#_0^T(0_\infty) \colon T_{0_+}\mathcal{H}^+_0 \times T_{0_-} \mathcal{H}^-_0 \to T_{0_T}\mathcal{H}^T$$
does not restrict to an isomorphism from $T_{0_\infty}W^\infty$ to $T_{0_T}W^T$. Therefore to define the gluing
map at an infinitesimal level at the origin we need to project $T_{0_T}\mathcal{H}^T$ to $T_{0_T}W^T$. In order to define
such a projection we have to choose a complement of $T_{0_T}W^T$ in $T_{0_T}\mathcal{H}$. Note that 
$$T_{0_T}\mathcal{H}^T=W^{1,2}\big([-T,T],\mathbb{R}^n\big).$$
To define a complement there are different options. For example one could choose the $L^2$-complement, see
\cite{schwarz} or the $W^{1,2}$-complement, see \cite{salamon}. We choose as complement the following. Abbreviate
$$\pi \colon \mathbb{R}^n \to \mathbb{R}^\mu$$
the orthogonal projection, where we recall that $\mu$ is the Morse index. We set
\begin{equation}\label{compl}
K^T:=\Big\{\xi \in W^{1,2}\big([-T,T],\mathbb{R}^n\big): \pi(\xi(T))=0,\,\,(\mathbb{1}-\pi)(\xi(-T))=0\Big\}.
\end{equation}
Observe that
$$T_{0_T}W^T \oplus K^T=T_{0_T}\mathcal{H}^T.$$
Denote by
$$\Pi^T \colon T_{0_T}\mathcal{H}^T \to T_{0_T} W^T$$
the projection along $K^T$. We define now the infinitesimal gluing map at the origin as 
$$\Gamma^T:=\Pi^T \circ d \#_0^T(0_\infty)|_{T_{0_\infty}W^\infty} \colon T_{0_\infty}W^\infty \to T_{0_T}W^T.$$ 
Our next goal is to describe the infinitesimal gluing map explicitly. We abbreviate
$$D^T:=D\sigma^T(0_T) \colon T_{0_T}\mathcal{H}^T \to \mathcal{E}_{0_T}$$
the vertical differential of the section $\sigma^T$ at $0_T$. Observe that $D^T$ is surjective and its kernel is given by
$$\mathrm{ker}D^T=T_{0_T}W^T.$$
In concrete terms the linear operator $D^T$ is given by
$$D^T \colon W^{1,2}\big([-T,T],\mathbb{R}^n\big) \to L^2\big([-T,T],\mathbb{R}^n\big), \quad \xi \mapsto \partial_s \xi
+H\xi$$
where $H$ is the Hessian at the origin as explained in (\ref{Hes}). Therefore $\xi \in \mathrm{ker}D^T$ if and only if 
$\xi$ is of the form
\begin{equation}\label{expon}
\xi(s)=e^{-Hs}\xi_0, \quad s \in [-T,T]
\end{equation}
for some $\xi_0 \in \mathbb{R}^n$. The differential of the evaluation map (\ref{eva2}) is given by
$$d\mathrm{ev}(0_T) \colon T_{0_T}W^T \to \mathbb{R}^n \times \mathbb{R}^n, \quad \xi \mapsto \big(\xi(-T),\xi(T)\big).$$
If we introduce the $n$-dimensional subspace
$$V^T:=\big\{(\xi_1,\xi_2) \in \mathbb{R}^n \times \mathbb{R}^n: \xi_1=e^{2TH}\xi_2\big\}$$
We learn from (\ref{expon}) that 
$$d\mathrm{ev}(0_T) \colon T_{0_T}W^T \to V^T$$
is a vector space isomorphism. 
We further abbreviate
$$V^\infty:=\{0\}\times \mathbb{R}^{n-\mu} \times \mathbb{R}^\mu \times \{0\} \subset \mathbb{R}^n \times \mathbb{R}^n.$$
Similarly, we have a vector space isomorphism
$$d\mathrm{ev}(0_\infty) \colon T_{0_\infty}W^\infty \to V^\infty.$$
The vector space $\mathbb{R}^{n-\mu}$ parametrizes the stable directions and the vector space $\mathbb{R}^\mu$ the unstable ones. Set
$$\widetilde{\Gamma}^T:=d\mathrm{ev}(0_T) \circ \Gamma^T \circ d\mathrm{ev}(0_\infty)^{-1} \colon V^\infty \to V^T.$$
Since $\widetilde{\Gamma}^T$ is conjugated to $\Gamma^T$ the study of both maps is equivalent. However, it is notationally easier to write down an explicit formula for $\widetilde{\Gamma}^T$ then for $\Gamma^T$ because both
vector spaces $V^T$ and $V^\infty$ are subspaces of $\mathbb{R}^n \times \mathbb{R}^n$. In order to do that
we further decompose the Hessian
$$H=H_- \oplus H_+$$
according to the splitting $\mathbb{R}^n=\mathbb{R}^\mu \oplus \mathbb{R}^{n-\mu}$. 
Then $H_-$ is the diagonal negative definite matrix
$$H_-=\left(\begin{array}{ccc}
a_1 & \cdots & 0\\
\vdots & \ddots &\vdots\\
0 & \cdots & a_\mu
\end{array}\right),$$
and $H_+$ is the diagonal positive definite matrix
$$H_+=\left(\begin{array}{ccc}
a_{\mu+1} & \cdots & 0\\
\vdots & \ddots &\vdots\\
0 & \cdots & a_n
\end{array}\right).$$
Using these notions we can now give an explicit formula for the infinitesimal gluing map. 
\begin{lemma}
The map $\widetilde{\Gamma}^T$ is given for $\xi_+ \in \mathbb{R}^{n-\mu}$ and $\xi_- \in \mathbb{R}^\mu$ by the formula
\begin{equation}\label{infglu}
\widetilde{\Gamma}^T\big(0,\xi_+,\xi_-,0\big)=\big(e^{2T H_-}\xi_-,\xi_+,\xi_-,e^{-2T H_+}\xi_+\big).
\end{equation}
\end{lemma}
\begin{figure}[h]
\begin{center}
 \includegraphics[scale=1]{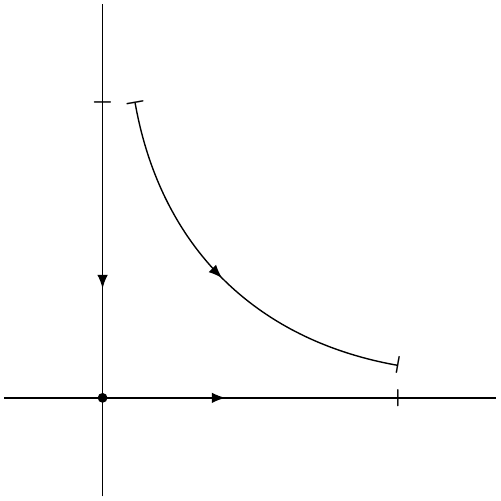}
 \caption{Infinitessimal gluing map.}
\end{center}
\end{figure}
\textbf{Proof: } Since the evaluation maps are invariant under pregluing by (\ref{evainv}) we have
$$d\mathrm{ev}(0_T) \circ d\#^T_0(0_\infty) \circ d\mathrm{ev}(0_\infty)^{-1}\big(0,\xi_+,\xi_-,0\big)=\big(0,\xi_+,\xi_-,0\big).$$
Now formula (\ref{infglu}) follows from the definition of the complement $K^T$ in (\ref{compl}). \hfill $\square$
\\ \\
\textbf{Remark: } It is interesting to note, that by (\ref{infglu}) the infinitesimal gluing map does not depend on the
choice of the cutoff function $\beta$ used to define the pregluing map in (\ref{pregl}). 

\subsection{The local gluing map}

\begin{thm}\label{locglu}
There exists an open neighbourhood $\mathcal{U} \subset W^\infty$ of $0_\infty$ and $T_0>2$ for which there exists a smooth family
of smooth maps
$$\#^T \colon \mathcal{U} \to W^T, \quad T \geq T_0$$
diffeomorphic to their image with the property that 
$$\lim_{T \to \infty}\mathrm{ev} \circ \#^T=\mathrm{ev}.$$
\end{thm}
\textbf{Proof: } Because $D^T$ is surjective it follows that 
$$D^T|_{K^T} \colon K^T \to \mathcal{E}_{0_T}$$
is bijective and therefore by the open mapping theorem it has a continuous inverse
$$R^T \colon \mathcal{E}_{0_T} \to K^T.$$ 
Note that $R^T$ is a right-inverse to $D^T$, i.e.,
\begin{equation}\label{rightinv}
D^T R^T=\mathbb{1} \colon \mathcal{E}_{0_T} \to \mathcal{E}_{0_T}.
\end{equation}
From (\ref{rightinv}) follows
$$(R^T D^T)^2=R^T(D^T R^T)D^T=R^T D^T$$
so that $R^T D^T$ is a projection, namely the projection from $T_{0_T} \mathcal{H}^T$ to $K^T$ along
$\mathrm{ker}D^T=T_{0_T}W^T$. Therefore the projection from $T_{0_T} \mathcal{H}^T$ to $T_{0_T}W^T$ along $K^T$ is given by
$$\Pi^T=\mathbb{1}-R^T D^T.$$
By the arguments in the proof of Step\,2 in Lemma~\ref{semifredi}, there exists a $T$-independent constant $c$ such that 
$$||\xi||_{1,2} \leq c||D^T \xi||_2, \quad \xi \in K^T.$$
Therefore the operator norms of the right inverses are uniformly bounded
$$||R^T|| \leq c.$$
The Theorem now follows by Newton-Picard iteration. \hfill $\square$

\subsection{Gluing gradient flow lines}

Assume that $M$ is a closed $n$-dimensional manifold, $f \colon M \to \mathbb{R}$ is a Morse function, $g$ is a Riemannian metric on $M$ which is Morse-Smale with respect to $f$, and $c_1,c_2,c_3$ are three critical points of $f$ whose Morse indices satisfy
$$\mu(c_1)=\mu(c_2)+1=\mu(c_3)+2.$$
Abbreviate by $W^\infty$ the product of the stable and unstable manifold of $c_2$, i.e.,
$$W^\infty=W^+_{c_2}\times W^-_{c_2}.$$
Observe that
$$\mathrm{dim}(W^\infty)=\mathrm{dim}(W^+_{c_2})+\mathrm{dim}(W^-_{c_2})=n-\mu(c_2)+\mu(c_2)=n.$$
Choose a family of gluing maps as in Theorem~\ref{locglu} and let $\mathcal{U} \subset W^\infty$ and $T_0>2$ be as in the Theorem. Suppose that $x^1$ is a gradient flow line from $c_1$ to $c_2$ and $x^2$ is a gradient flow line from $c_2$ to 
$c_3$. After reparametrization we can assume that
$$\big(x^1|_{[0,\infty)},x^2|_{(-\infty,0]}\big) \in \mathcal{U}.$$
Abbreviate
$$a:=f(x^1(0)), \quad b:=f(x^2(0)).$$
Note that due to the fact that nonconstant gradient flow lines flow downhill by Lemma~\ref{downlem}, we have the string of inequalities
\begin{equation}\label{inestring}
f(c_1)>a>f(c_2)>b>f(c_3).
\end{equation}
We introduce further
$$W^{-,a}_{c_1}:=\big\{x \in W^-_{c_1}: f(x(0))=a\big\}$$
and
$$W^{+,b}_{c_3}:=\big\{x \in W^+_{c_3}: f(x(0))=b\big\}.$$
Note that 
$$x^1|_{(-\infty,0]} \in W^{-,a}_{c_1}, \quad x^2|_{[0,\infty)} \in W^{+,b}_{c_3}.$$
If we consider the map
$$f \circ \mathrm{ev} \colon W^-_{c_1} \to \mathbb{R}$$
we can alternatively write
$$W^{-,a}_{c_1}=(f \circ \mathrm{ev})^{-1}(a).$$
Because $x^1$ is nonconstant it follows again from Lemma~\ref{downlem}, namely that gradient flow lines flow downhill, that
$$d(f \circ \mathrm{ev})\big(x^1|_{(-\infty,0]}\big) \neq 0$$
so that $W^{-,a}_{c_1}$ is locally around $x^1|_{(-\infty,0]}$ a codimension one submanifold of 
$W^-_{c_1}$ and in particular a manifold of dimension 
$$\mathrm{dim}(W^{-,a}_{c_1})=\mu(c_1)-1=\mu(c_2).$$ 
In the same way it holds that $W^{+,b}_{c_3}$
is locally around $x^2|_{[0,\infty)}$ a codimension one submanifold of $W^+_{c_3}$, hence a manifold of dimension
$$\mathrm{dim}W^{+,b}_{c_3}=n-\mu(c_3)-1=n-\mu(c_2).$$ 
Since the Riemannian metric $g$ is Morse-Smale it follows that we have transverse intersections
\begin{equation}\label{tran}
\mathrm{ev}(W^\infty) \pitchfork \mathrm{ev}\big(W^{-,a}_{c_1} \times W^{+,b}_{c_3}\big).
\end{equation}
Observe that
\begin{eqnarray*}
\mathrm{dim}(W^\infty)+\mathrm{dim}(W^{-,a}_{c_1}\times W^{+,b}_{c_3})&=&\mathrm{dim}(W^\infty)+\mathrm{dim}(W^{-,a}_{c_1})+
\mathrm{dim}(W^{+,b}_{c_3})\\
&=&n+\mu(c_2)+n-\mu(c_2)\\
&=&2n\\
&=&\mathrm{dim}(M \times M).
\end{eqnarray*}
Hence (\ref{tran}) and the fact that the evaluation maps are embeddings by Theorem~\ref{embedding} imply that
$$\big(x^1(0), x^2(0)\big) \in \mathrm{ev}(W^\infty) \cap \mathrm{ev}\big(W^{-,a}_{c_1} \times W^{+,b}_{c_3}\big)$$
is isolated. 
In view of Theorem~\ref{locglu} and (\ref{tran}) there exists $T_1 \geq T_0$ such that 
$$\mathrm{ev} \circ \#^T(\mathcal{U}) \pitchfork \mathrm{ev}\big(W^{-,a}_{c_1} \times W^{+,b}_{c_3}\big), \quad T \geq T_1.$$
Therefore it is possible to find an open neighbourhood $\mathcal{V} \subset \mathcal{U}$ of $\big(x^1|_{[0,\infty)},x^2|_{(-\infty,0]}\big)$ for which there exists $T_2 \geq T_1$ with the property that 
$$\#\Big(\mathrm{ev} \circ \#^T(\mathcal{V}) \cap \mathrm{ev}\big(W^{-,a}_{c_1} \times W^{+,b}_{c_3}\big)\Big)=1, \quad T \geq T_2,$$
where $\#$ denotes the cardinality of the set. This means that there exists a unique gradient flow line $x_T$ from
$c_1$ to $c_3$ meeting the requirements 
\begin{itemize}
 \item[(i)] $x_T|_{[-T,T]} \in \#^T(\mathcal{V})$,
 \item[(ii)] $f(x_T(-T))=a$,
 \item[(iii)] $f(x_T(T))=b$.
\end{itemize}
We define the gluing map
$$\rho:=\rho_{(x^1,x^2)} \colon (T_2,\infty) \to \mathcal{M}(c_1,c_3), \quad T \mapsto \big[x_T\big].$$
\begin{figure}[h]
\begin{center}
 \includegraphics[scale=1]{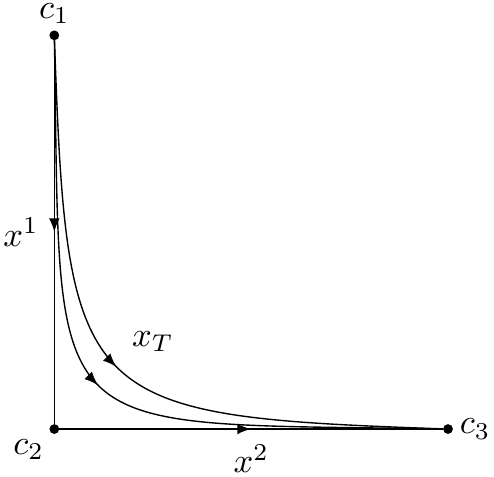}
 \caption{Global gluing.}
\end{center}
\end{figure}
\begin{lemma}\label{injective}
The gluing map $\rho$ is injective. 
\end{lemma}
\textbf{Proof: } We follow \cite{simcevic} and construct a smooth left inverse for the gluing map. By (\ref{inestring}) and the fact that gradient flow lines flow downhill established in Lemma~\ref{downlem}, there exist smooth maps
$$\tau_a \colon \widetilde{\mathcal{M}}(c_1,c_3) \to \mathbb{R}, \quad \tau_b \colon \widetilde{\mathcal{M}}(c_1,c_3) \to \mathbb{R}$$
which for $x \in \widetilde{\mathcal{M}}(c_1,c_3)$ are implicitly determined by the conditions
$$f(x(\tau_a(x)))=a, \qquad f(x(\tau_b(x)))=b.$$
Note that if $r \in \mathbb{R}$ then
$$\tau_a(r_*x)=\tau_a(x)-r,\qquad \tau_b(r_*x)=\tau_b(x)-r$$
so that the map
$$\tau=\frac{1}{2}(\tau_b-\tau_a)$$
satisfies
$$\tau(r_* x)=\tau(x)$$
and therefore it can be interpreted as a smooth map on the quotient of unparametrized gradient flow lines 
$$\tau \colon \mathcal{M}(c_1,c_3) \to \mathbb{R}.$$
By properties (ii) and (iii) of $x_T$ it holds that
$$\tau_b\big(x_T\big)=T,\qquad \tau_a\big(x_T\big)=-T.$$
Therefore for $T \in (T_2,\infty)$ we compute
$$\tau \circ \rho(T)=\tau\big(\big[x_T\big]\big)=\frac{1}{2}\big(\tau_b\big(x_T\big)-\tau_a\big(x_T\big)\big)=T$$
so that
$$\tau\circ \rho=\mathbb{1}|_{(-T_2,\infty)}.$$
This shows that the gluing map $\rho$ has a smooth left inverse and in particular is injective. \hfill $\square$
\begin{lemma}\label{kleben}
Assume that $y_\nu \in \widetilde{\mathcal{M}}(c_1,c_3)$ is a sequence of gradient flow lines from $c_1$ to $c_3$ such that
$$y_\nu \xrightarrow
 {
  \substack{\mathrm{Floer-Gromov}}
 }(x^1,x^2).$$
Then there exists $\nu_0$ such that for every $\nu \geq \nu_0$ there exists $\tau_\nu \in (T_2,\infty)$ such that
$$[y_\nu]=[x_{\tau_\nu}] \in \mathcal{M}(c_1,c_3).$$
\end{lemma}
\textbf{Proof: } We use the map $\tau \colon \mathcal{M}(c_1,c_3) \to \mathbb{R}$ which we introduced in the proof of Lemma~\ref{injective} and set
$$\tau_\nu:=\tau([y_\nu]).$$
Maybe after reparametrization we can assume that 
$$f(y_\nu(-\tau_\nu))=a$$
so that it follows from the definition of $\tau$ that
$$f(y_\nu(\tau_\nu))=b.$$
Since $y_\nu$ Floer-Gromov converges to $(x^1,x^2)$ there exists $\nu_0$ with the property that
$$y_\nu|_{[-\tau_\nu,\tau_\nu]} \in \#^{\tau_\nu}(\mathcal{V}), \quad \nu \geq \nu_0.$$
Therefore
$$y_\nu=x_{\tau_\nu}, \quad \nu \geq \nu_0$$
and the Lemma follows. \hfill $\square$
\\ \\
We are now ready to state the precise version of Theorem~\ref{dd3} and prove it. 
\begin{thm}\label{dd3f}
The compactified moduli space $\overline{\mathcal{M}}(c_1,c_3)$ endowed with the Floer-Gromov topology has the structure
of a one-dimensional manifold with boundary, where the boundary consists of once broken gradient flow lines. 
\end{thm}
\textbf{Proof: } The noncompactified moduli space $\mathcal{M}(c_1,c_3)$ is a (not necessarily compact) one-dimensional manifold by Theorem~\ref{dimpar}. By Lemma~\ref{injective} and Lemma~\ref{kleben} the gluing map provides local charts at the broken gradient flow lines so that broken gradient flow lines become boundary points. \hfill $\square$
\\ \\
As explained in Paragraph~\ref{strategy} Theorem~\ref{dd3f} implies
\begin{cor}
A boundary has no boundary, i.e. $\partial^2=0$.
\end{cor} 
In view of the Corollary we can now for a Morse function $f$ on a closed manifold $M$ and a Morse-Smale metric $g$ with respect to $f$ associate the $\mathbb{Z}_2$-vector space
$$HM_*(f,g)=\frac{\mathrm{ker}\partial}{\mathrm{im}\partial}.$$

\section{Invariance of Morse homology}\label{invariance}

\subsection{Morse-Smale pairs}

What we learned so far is that if $M$ is a closed manifold, $f \colon M \to \mathbb{R}$ is a Morse function on $M$ and
$g$ is a Morse-Smale metric with respect to $f$ we have a graded vector space $CM_*(f)$ on which we can define a boundary operator $\partial=\partial_{f,g}$ which leads to Morse homology $HM_*(f,g)$. We refer to the tuple $(f,g)$ as a 
\emph{Morse-Smale pair}. The goal of this section is to show that up to canonical isomorphism the graded $\mathbb{Z}_2$-vector space $HM_*(f,g)$ is independent of the choice of the Morse-Smale pair. Namely we prove 
\begin{thm}\label{caniso}
Suppose that $(f_\alpha,g_\alpha)$ and $(f_\beta,g_\beta)$ are two Morse-Smale pairs on $M$. Then there exists a canonical
isomorphism
$$\Phi^{\beta \alpha} \colon HM_*(f_\alpha,g_\alpha) \to HM_*(f_\beta,g_\beta).$$
\end{thm}

\subsection{Construction of the canonical isomorphism}

We consider gradient flow lines for time dependent functions and time dependent Riemannian metrics interpolation between
$f_\alpha$ and $f_\beta$, respectively $g_\alpha$ and $g_\beta$. Namely for $s \in \mathbb{R}$ we choose a smooth family of functions $f_s$ for which there exists $T>0$ such that
$$f_s=\left\{\begin{array}{cc}
f_\alpha & s \leq -T\\
f_\beta & s \geq T,
\end{array}\right.$$
and similarly a smooth family of Riemannian metrics $g_s$ such that
$$g_s=\left\{\begin{array}{cc}
g_\alpha & s \leq -T\\
g_\beta & s \geq T.
\end{array}\right.$$
Similarly as in the case of the boundary operator we define a linear map
$$\phi^{\beta \alpha} \colon CM_*(f_\alpha) \to CM_*(f_\beta)$$
which for a basis vector $c_1 \in \mathrm{crit}(f_\alpha)$ is given by 
$$\phi^{\beta \alpha}(c_1)=\sum_{\substack{c_2 \in \mathrm{crit}(f_\beta)\\ \mu(c_2)=\mu(c_1)}}
\#_2\big\{\textrm{gradient flow lines from}\,\,c_1\,\,\textrm{to}\,\,c_2\big\}c_2$$
Here with gradient flow lines from $c_1$ to $c_2$ we mean solutions $x \in C^\infty(\mathbb{R},M)$ of the ODE
\begin{equation}\label{timegrad}
\partial_s x(s)+\nabla_{g_s} f_s(x(s))=0, \quad s \in \mathbb{R}
\end{equation}
satisfying the asymptotic conditions
$$\lim_{s \to- -\infty}x(s)=c_1, \qquad \lim_{s \to \infty}x(s)=c_2.$$
Different from the time independent case the ODE (\ref{timegrad}) is not invariant anymore under timeshift. As in the case of the boundary operator we have to choose the homotopy $(f_s,g_s)$ generic so that we get a well defined count of solutions. In order to guarantee that this count generically is finite we need an analogon of Floer-Gromov compactness for time dependent gradient flow lines. 

\subsection{Floer-Gromov convergence for time dependent gradient flow lines}

We continue with the notations of the previous paragraph. The analogous notions of broken gradient flow line and Floer-Gromov convergence to a broken gradient flow line are the following. 
\begin{fed}
A \emph{broken gradient flow line} from $c_1$ to $c_2$ of $\nabla_{g_s}f_s$ is a tuple
$$y=\{x^k\}_{1 \leq k \leq n}, \quad n \in \mathbb{N}$$
such that the following holds true.
\begin{description}
 \item[(i)] There exists $k_0=k_0(y) \in \{1,\ldots,n\}$ such that $x^{k_0}$ is a gradient flow line of $\nabla_{g_s}f_s$.
 \item[(ii)] $x^k$ for $1 \leq k \leq k_0$ is a nonconstant gradient flow line of $\nabla_{g_\alpha}f_\alpha$ and
  $x^k$ for $k_0<k\leq n$ is a nonconstant gradient flow line of $\nabla_{g_\beta}f_\beta$.
 \item[(iii)] $\lim_{s \to-\infty}x^1(s)=c_1$,\\
              $\lim_{s \to \infty}x^k(s)=\lim_{s \to -\infty}x^{k+1}(s)$ for $k \in \{1,\ldots,n-1\}$,\\
              $\lim_{s \to \infty} x^n(s)=c_2$. 
\end{description}
\end{fed}
\begin{fed}
Assume that $x_\nu \in C^\infty(\mathbb{R},M)$ is a sequence of gradient flow lines of $\nabla_{g_s}f_s$ from $c_1$ to $c_2$ and $y=\{x^k\}_{1 \leq k \leq n}$ is a broken gradient flow line of $\nabla_{g_s}f_s$ from $c_1$ to $c_2$. We say
$x_\nu$ \emph{Floer-Gromov converges} to $y$ if for $1 \leq k \leq n$ there exists a sequence $r_\nu^k \in \mathbb{R}$ with $r_\nu^{k_0}=0$ such that
$$(r_\nu^k)_*x_\nu \xrightarrow
 {
  \substack{C^\infty_{\mathrm{loc}}}
 }
 x^k.
$$
\end{fed}
The analogon of Theorem~\ref{floergromov} about Floer-Gromov compactness in the time dependent case follows similarly. There is one interesting difference to notice. For $x \in C^\infty(\mathbb{R},M)$ we define the energy in the time-dependent case
by
$$E(x):=E_{g_s}(x):=\int_{-\infty}^\infty ||\partial_s x||_{g_s}^s ds \in [0,\infty].$$
While by Corollary~\ref{enbound} in the time independent case the energy of a gradient flow line can be estimated in terms of the asymptotic action values this is not true anymore in the time dependent case. Nevertheless there is still a uniform bound on the energy as the following lemma tells us. 
\begin{lemma}\label{timeenest}
There is a constant $c=c(f_s)>0$ such that for all gradient flow lines $x$ of $\nabla_{g_s}f_s$ from $c_1$ to $c_2$ the energy $E(x)$ can be uniformly estimated from above by $c$.
\end{lemma}
\textbf{Proof: } We estimate using (\ref{timegrad})
\begin{eqnarray*}
E(x)&=&\int_{-\infty}^\infty ||\partial_s x||_{g_s}^2 ds\\
&=&-\int_{-\infty}^\infty g_s\big(\nabla_{g_s} f_s(x),\partial_s x\big) ds\\
&=&-\int_{-\infty}^\infty df_s(x)\partial_s x ds\\
&=&-\int_{-\infty}^\infty \frac{d}{ds}(f_s(x(s)))ds+\int_{-\infty}^\infty (\partial_s f_s)(x(s))ds\\
&=&f_\alpha(c_1)-f_\beta(c_2)+\int_{-\infty}^\infty (\partial_s f_s)(x(s))ds\\
&\leq&||f_\alpha||_\infty+||f_\beta||_\infty+2T||\partial_s f_s||_\infty.
\end{eqnarray*}
This proves the lemma. \hfill $\square$

\subsection{Generic breaking of time dependent gradient flow lines}

Abbreviate by $\mathcal{N}(f_s,g_s;c_1,c_2)$ the moduli space of gradient flow lines of $\nabla_{g_s} f_s$ from $c_1$ to $c_2$. In view of Floer-Gromov compactness the analogons of the generic breaking behaviour discussed in Section~\ref{strategy}
in the time dependent case become
\begin{prop}\label{chain1}
If $\mu(c_1)=\mu(c_2)$, then for a generic homotopy $(f_s,g_s)$ the moduli space $\mathcal{N}(f_s,g_s;c_1,c_2)$ is a finite
set.
\end{prop}
\begin{prop}\label{chain2}
If $\mu(c_1)=\mu(c_2)+1$, then for a generic homotopy $(f_s,g_s)$ the moduli space $\mathcal{N}(f_s,g_s;c_1,c_2)$ can be
compactified to a one-dimensional manifold with boundary $\overline{\mathcal{N}}(f_s,g_s,c_1,c_2)$ such that the boundary is given by
\begin{eqnarray}\label{chainform}
\partial \overline{\mathcal{N}}(f_s,g_s;c_1,c_2)&=&\bigsqcup_{\substack{c \in \mathrm{crit}(f_\alpha)\\ \mu(c)=\mu(c_1)-1}}
\mathcal{M}(f_\alpha,g_\alpha;c_1,c) \times \mathcal{N}(f_s,g_s;c,c_2)\\ \nonumber
& & \sqcup
\bigsqcup_{\substack{c \in \mathrm{crit}(f_\beta)\\ \mu(c)=\mu(c_1)}}
\mathcal{N}(f_s,g_s;c_1,c) \times \mathcal{M}(f_\beta,g_\beta;c,c_2).
\end{eqnarray}
\end{prop}
The first term on the righthand side of (\ref{chainform}) occurs since a gradient flow line of the time independent Morse function $f_\alpha$ can break up at the left, while the second one takes account of gradient flow lines of $f_\beta$ breaking up at the right.
\begin{figure}[h]
\begin{center}
\includegraphics[scale=1]{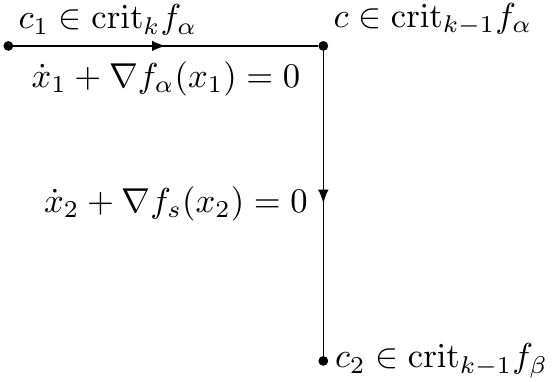} \quad
\includegraphics[scale=1]{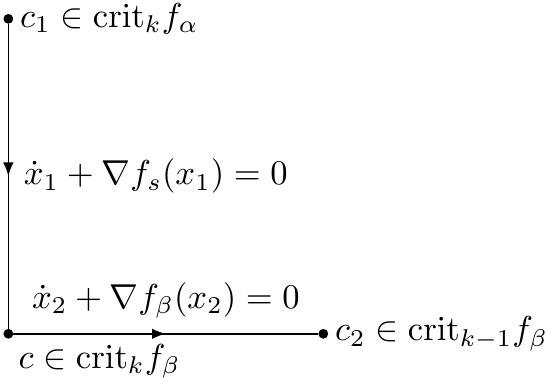}
\caption{Two ways of breaking for gradient flow lines of homotopies.}
\end{center}
\end{figure}

\subsection{A chain map}

From Proposition~\ref{chain1} we obtain that the linear map
$$\phi^{\beta \alpha} \colon CM_*(f_\alpha) \to CM_*(f_\beta)$$
is well defined. Because modulo two the number of boundary points of a compact one-dimensional manifold vanishes we derive from (\ref{chainform}) the algebraic equation
$$\partial^\beta \phi^{\beta \alpha}+\phi^{\beta \alpha} \partial^\alpha=0$$
where $\partial^\alpha$ is the Morse boundary operator on $CM_*(f_\alpha)$ and $\partial^\beta$ the one on $CM_*(f_\beta)$.
Modulo two this equation is equivalent to
\begin{equation}\label{chainmap}
\partial^\beta \phi^{\beta \alpha}=\phi^{\beta \alpha}\partial^\alpha,
\end{equation}
i.e. the linear map $\phi^{\beta \alpha}$ interchanges the two boundary operators. In terms of homological algebra such a map is called a \emph{chain map}. The interesting point about a chain map is that it induces a linear map on homology
$$\Phi^{\beta \alpha} \colon HM_*(f_\alpha,g_\alpha) \to HM_*(f_\beta,g_\beta), \quad
[\xi] \mapsto [\phi^{\beta \alpha} \xi].$$
To see that this is well defined we observe the following. Since $\xi \in CM_*(f_\alpha)$ represents a homology class in
$HM_*(f_\alpha,g_\alpha)$ it has to lie in the kernel of $\partial^\alpha$, i.e.,
$$\partial^\alpha \xi=0.$$
Using (\ref{chainmap}) we compute
$$\partial^\beta \phi^{\beta \alpha}\xi=\phi^{\beta \alpha}\partial^\alpha \xi=0$$
so that $\phi^{\beta \alpha}\xi$ lies in the kernel of $\partial^\beta$ and actually represents a homology class in
$HM_*(f_\beta,g_\beta)$. We next check that the class $[\phi^{\beta \alpha} \xi] \in HM_*(f_\beta,g_\beta)$ does not depend
on the choice of $\xi$ as a representative of its homology class. To see this observe that if
$$[\xi]=[\xi'] \in HM_*(f_\alpha,g_\alpha)$$
there exists $\eta \in CM_{*+1}(f_\alpha)$ such that 
$$\xi'=\xi+\partial^\alpha \eta.$$
Using (\ref{chainmap}) once more we compute
$$\phi^{\beta \alpha} \xi'=\phi^{\beta \alpha} \xi+\phi^{\beta \alpha} \partial^\alpha \eta=
\phi^{\beta \alpha} \xi+\partial ^\beta \phi^{\beta \alpha} \eta$$
so that we have
$$[\phi^{\beta \alpha}\xi]=[\phi^{\beta \alpha}\xi'] \in HM_*(f_\beta,g_\beta).$$
This proves that $\Phi^{\beta \alpha}$ is well defined.

\subsection{Homotopy of Homotopies}

What we achieved so far is that after choosing a generic homotopy $(f_s,g_s)$ we defined a linear map
$$\Phi^{\beta\alpha}=\Phi^{\beta\alpha}_{f_s,g_s} \colon HM_*(f_\alpha,g_\alpha) \to HM_*(f_\beta,g_\beta).$$
Our next goal is to show that this map is canonical.
\begin{prop}
The map $\Phi^{\beta \alpha}$ does not depend on the homotopy $(f_s,g_s)$.
\end{prop}
\textbf{Proof: } Given two generic homotopies $(f_s^0,g_s^0)$ and $(f_s^1,g_s^1)$ we interpolate between them by considering a smooth homotopy of homotopies $(f_s^r,g_s^r)$ for $r \in [0,1]$. If $c_1 \in \mathrm{crit}(f_\alpha)$ and
$c_2 \in \mathrm{crit}(f_\beta)$ we introduce the moduli space
\begin{eqnarray*}
& &\mathcal{R}(f_s^r,g_s^r;c_1,c_2):=\\
& &\Big\{(x,r): r \in [0,1],\,\, \partial_s x+\nabla_{g_s^r}f^r_s(x)=0,\,\,
\lim_{s \to -\infty}x(s)=c_1,\,\,\lim_{s \to \infty}x(s)=c_2\Big\}.
\end{eqnarray*}
Assume that $\mu(c_1)=\mu(c_2)$. We choose the homotopy of homotopies generic, so that we can compactify 
the moduli space $\mathcal{R}(f_s^r,g_s^r;c_1,c_2)$ to a one-dimensional manifold with boundary
$\overline{\mathcal{R}}(f_s^r,g_s^r;c_1,c_2)$. There are two obvious boundary components, namely the starting point of
the homotopy $r=0$ and the endpoint of the homotopy $r=1$. However, there might be additional contributions.
While for a generic homotopy there are no gradient flow lines when the difference of
Morse indices between the left and right asymptotic is minus one, in a one parameter family such gradient flow lines can occur at isolated points in the homotopy, even when the one parameter family of homotopies is chosen generic. Therefore therefore there is a finite subset $\mathcal{S} \subset (0,1)$ such that
\begin{eqnarray}\label{homhom}
\partial \mathcal{R}(f_s^r,g_s^r;c_1,c_2)&=&\mathcal{N}(f_s^0,g_s^0;c_1,c_2)\sqcup \mathcal{N}(f_s^1,g_s^1;c_1,c_2)\\ 
\nonumber
& &\sqcup \bigsqcup_{\substack{r_0 \in \mathcal{S}\\ c \in \mathrm{crit}(f_\alpha)\\ \mu(c)=\mu(c_1)-1}}
\mathcal{M}(f_\alpha,g_\alpha;c_1,c)\times \mathcal{N}(f_s^{r_0},g_s^{r_0};c,c_2)\\ \nonumber
& &\sqcup \bigsqcup_{\substack{r_0 \in \mathcal{S}\\ c \in \mathrm{crit}(f_\beta)\\ \mu(c)=\mu(c_2)+1}}
\mathcal{N}(f_s^{r_0},g_s^{r_0};c_1,c) \times \mathcal{M}(f_\beta,g_\beta;c,c_2).
\end{eqnarray}
\begin{figure}
\begin{center}
\includegraphics[scale=1]{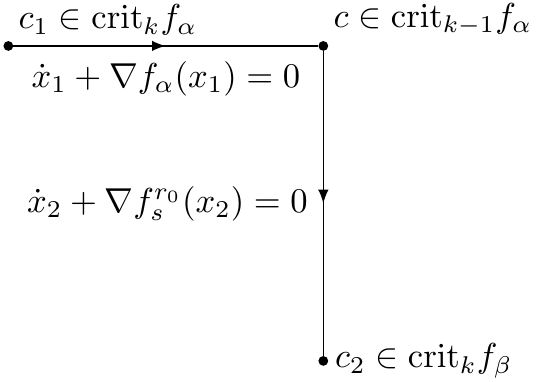} \quad
\includegraphics[scale=1]{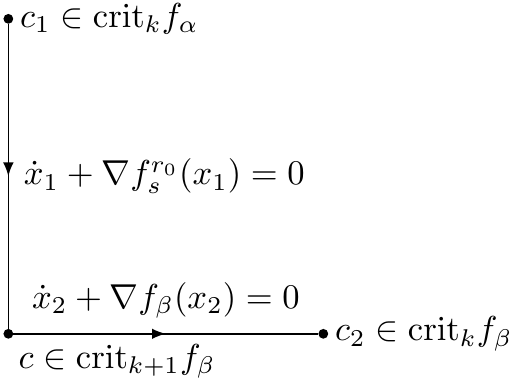}
\caption{Breaking gradient flow lines for homotopies of homotopies.}
\end{center}
\end{figure}
Let
$$T \colon CM_*(f_\alpha) \to CM_{*+1}(f_\beta)$$
be the linear map which on a basis vector $c_1 \in \mathrm{crit}(f_\alpha)$ is given by
$$T(c_1)=\sum_{\substack{c_2 \in \mathrm{crit}(f_\beta)\\ \mu(c_2)=\mu(c_1)+1}}
\#\Big( \bigsqcup_{r_0 \in \mathcal{S}} \mathcal{N}(f_s^{r_0},g_s^{r_0};c_1,c_2)\Big)c_2.$$
Abbreviate
$$\phi^{\beta \alpha}_0:=\phi^{\beta \alpha}_{f_s^0,g_s^0}, \qquad
\phi^{\beta \alpha}_1:=\phi^{\beta \alpha}_{f_s^1,g_s^1}.$$
Using again that the number of boundary points modulo two of a one-dimensional manifold with boundary vanishes the geometric statement (\ref{homhom}) translates into the algebraic statement
$$\phi^{\beta \alpha}_0+\phi^{\beta \alpha}_1+\partial^\beta T+T\partial^\alpha=0$$
which modulo two can be expressed equivalently as
\begin{equation}\label{equivalence}
\phi^{\beta \alpha}_0=\phi^{\beta \alpha}_1+\partial^\beta T+T\partial^\alpha.
\end{equation}
In terms of homological algebra this means that $\phi^{\beta \alpha}_0$ is \emph{chain homotopy equivalent} to
$\phi^{\beta \alpha}_1$. Chain homotopy equivalent chain maps induce the same map on homology. Namely if
$$\Phi^{\beta \alpha}_0, \Phi^{\beta \alpha}_1 \colon HM_*(f_\alpha,g_\alpha) \to HM_*(f_\beta,g_\beta)$$
are the by $\phi^{\beta \alpha}_0$, respectively $\phi^{\beta \alpha}_0$, induced maps on homology we claim that
$$\Phi^{\beta \alpha}_0=\Phi^{\beta \alpha}_1.$$
To see how this follows from (\ref{equivalence}) pick a class
$[\xi] \in HM_*(f_\alpha,g_\alpha)$. Its representative necessarily lies in the kernel of $\partial^\alpha$, i.e.,
$$\partial^\alpha \xi=0.$$
Hence we compute using (\ref{equivalence})
\begin{eqnarray*}
\Phi^{\beta \alpha}_0[\xi]&=&[\phi^{\beta \alpha}_0 \xi]=[\phi^{\beta \alpha}_1 \xi+\partial^\beta T\xi]
=[\phi^{\beta \alpha}_1 \xi]=\Phi^{\beta \alpha}_1[\xi].
\end{eqnarray*}
This shows that $\Phi^{\beta \alpha}_1$ coincides with $\Phi^{\beta \alpha}_0$ and the proposition is proven.
\hfill $\square$

\subsection{Functoriality}

In this paragraph we study the composition of the canonical maps between the Morse homologies for different Morse-Smale pairs

\begin{prop}\label{functo}
Suppose that $(f_\alpha,g_\alpha)$, $(f_\beta,g_\beta)$, and $(f_\gamma.g_\gamma)$ are three Morse-Smale pairs. Then
$$\Phi^{\gamma \beta} \circ \Phi^{\beta \alpha}=\Phi^{\gamma \alpha} \colon HM_*(f_\alpha,g_\alpha) \to
HM_*(f_\gamma,g_\gamma).$$
\end{prop}
\textbf{Proof: } The proposition follows by gluing gradient flow lines. Assume that
$(f^{\beta \alpha}_s,g^{\beta \alpha}_s)$ and $(f^{\gamma \beta}_s,g^{\gamma \beta}_s)$ are generic homotopies. Suppose further that $c_1 \in \mathrm{crit}(f^\alpha)$, $c_2 \in \mathrm{crit}(f^\beta)$, and $c_3 \in \mathrm{crit}(f^\gamma)$
such that
$$\mu(c_1)=\mu(c_2)=\mu(c_3).$$
For large $R$ define
$$(f^{\gamma \beta}\#_R f^{\beta \alpha})_s=\left\{\begin{array}{cc}
f^{\beta \alpha}_{s+R} & s \leq 0\\
f^{\gamma \beta}_{s-R} & s \geq 0
\end{array}\right.$$
and 
$$(g^{\gamma \beta}\#_R f^{\beta \alpha})_s=\left\{\begin{array}{cc}
g^{\beta \alpha}_{s+R} & s \leq 0\\
g^{\gamma \beta}_{s-R} & s \geq 0.
\end{array}\right.$$ 
By gluing trajectories together we find $R_0>0$ such that for $R \geq R_0$ there is a bijection of moduli spaces
\begin{eqnarray*}
& &\mathcal{N}\big((f^{\gamma \beta}\#_R f^{\beta \alpha})_s,(g^{\gamma \beta}\#_R g^{\beta \alpha})_s;c_1,c_3\big) \cong\\
& &\mathcal{N}\big(f^{\beta \alpha}_s, g^{\beta \alpha}_s;c_1,c_2\big) \times
\mathcal{N}\big(f^{\gamma \beta}_s, g^{\gamma \beta}_s;c_2,c_3\big).
\end{eqnarray*}
This implies that
$$\phi^{\gamma \beta} \circ \phi^{\beta \alpha}=\phi^{\gamma \alpha} \colon CM_*(f_\alpha) \to CM_*(f_\gamma)$$
and therefore a fortiori
$$\Phi^{\gamma \beta} \circ \Phi^{\beta \alpha}=\Phi^{\gamma \alpha} \colon HM_*(f_\alpha, g_\alpha) \to 
HM_*(f_\gamma,g_\gamma).$$
The proposition is proved. \hfill $\square$

\subsection{The trivial homotopy}

In this paragraph we explain that the canonical map from the Morse homology of a fixed Morse-Smale pair to itself is the identity by exploring the trivial homotopy. 
\begin{prop}\label{identi}
Suppose that $(f_\alpha,g_\alpha)$ is a Morse-Smale pair. Then
$$\Phi^{\alpha \alpha}=\mathbb{1} \colon HM_*(f_\alpha,g_\alpha) \to HM_*(f_\alpha,g_\alpha).$$
\end{prop}
\textbf{Proof: } We choose as homotopy from $(f_\alpha,g_\alpha)$ to $(f_\alpha,g_\alpha)$ the constant homotopy
$(f_\alpha,g_\alpha)$. Since $(f_\alpha,g_\alpha)$ is a Morse-Smale pair the constant homotopy is actually generic. 
Suppose that $c_1$ and $c_2$ are two critical points of $f_\alpha$ satisfying
$$\mu(c_1)=\mu(c_2).$$
Note that for the constant homotopy we have
$$\mathcal{N}(f_\alpha, g_\alpha;c_1,c_2)=\widetilde{\mathcal{M}}(f_\alpha,g_\alpha;c_1,c_2),$$
the moduli space of parametrized gradient flow lines from $c_1$ to $c_2$. Note that 
$$\mathrm{dim}\widetilde{\mathcal{M}}(f_\alpha,g_\alpha;c_1,c_2)=\mu(c_1)-\mu(c_2)=0.$$
However, if $c_1 \neq c_2$, then there is a free $\mathbb{R}$-action on the moduli space by timeshift and the quotient consisting of unparametrized gradient flow lines has negative dimension. This implies
$$\widetilde{\mathcal{M}}(f_\alpha,g_\alpha;c_1,c_2)=\emptyset, \qquad c_1 \neq c_2.$$
On the other hand if $c_1=c_2=:c$, then the moduli space has precisely one element, namely the critical point itself interpreted as a constant gradient flow line. Therefore
$$\#\widetilde{\mathcal{M}}(f_\alpha,g_\alpha;c,c)=1.$$
This two facts imply that $\Phi^{\alpha \alpha}$ is the identity and the proposition is proven. \hfill $\square$ 

\subsection{Proof that the canonical maps are isomorphisms}

In this paragraph we proof Theorem~\ref{caniso}, namely that the canonical maps are isomorphisms.
\\ \\
\textbf{Proof of Theorem~\ref{caniso}: } Suppose that $(f_\alpha,g_\alpha)$ and $(f_\beta,g_\beta)$ are Morse-Smale pairs.
From Proposition~\ref{functo} and Proposition~\ref{identi} we obtain
$$\Phi^{\beta \alpha} \circ \Phi^{\alpha \beta}=\Phi^{\beta \beta}=\mathbb{1} \colon HM_*(f_\beta,g_\beta) \to
HM_*(f_\beta,g_\beta).$$
Interchanging the roles of $\alpha$ and $\beta$ we also have
$$\Phi^{\alpha \beta} \circ \Phi^{\beta \alpha}=\mathbb{1} \colon HM_*(f_\alpha,g_\alpha) \to
HM_*(f_\alpha,g_\alpha).$$
This proves that
$$\Phi^{\alpha \beta} \colon HM_*(f_\alpha,g_\alpha) \to HM_*(f_\beta,g_\beta)$$
is an isomorphism with inverse
$$(\Phi^{\alpha \beta})^{-1}=\Phi^{\beta \alpha}$$
and the theorem follows. \hfill $\square$

\section{Spectral numbers}

Spectral numbers associate via a minimax process a critical value to a class in Morse homology. An interesting aspect of
spectral numbers is their continuity with respect to the Morse function. This allows to extend spectral numbers even to functions which are not Morse. Spectral numbers are as well a crucial ingredient to incorporate techniques from persistent homology into symplectic geometry, see \cite{polterovich-rosen-samvelyan-zhang}.

\subsection{Spectral numbers for Morse functions}

Suppose that $M$ is a closed manifold and $f \colon M \to \mathbb{R}$ is a Morse function. Assume that 
$\xi \in CM_*(f)$. Then $\xi$ can be written as a formal sum of critical points of $f$
$$\xi=\sum_{c \in \mathrm{crit}(f)} a_c c$$
where the coefficient $a_c$ lie in the field $\mathbb{Z}_2$. Assume that $\xi$ is not the zero vector. We define
$$\sigma_f(\xi):=\max\big\{f(c): c \in \mathrm{crit}(f),\,\,a_c \neq 0\big\}.$$
If $g$ is a Morse-Smale metric with respect to $f$ we have Morse homology $HM_*(f,g)$ so that for $\alpha \neq 0
\in HM_*(f,g)$ we can set
$$\sigma_f(\alpha):=\min\big\{\sigma(\xi):[\xi]=\alpha\big\}$$
and refer to $\sigma_f(\alpha)$ as the \emph{spectral number} associated to the homology class $\alpha$. At the moment it is not clear that the spectral number is independent of the Morse-Smale metric. However, when we prove in the following paragraph Lipschitz continuity of the spectral number with respect to $f$, we get as a byproduct that the spectral number is independent of the choice of the metric. 

\subsection{Lipschitz continuity} 

As we explained in Chapter~\ref{invariance} there is a canonical isomorphism between Morse homologies associated to different Morse-Smale pairs. This allows us to fix a homology class in Morse homology up to canonical isomorphism and explore how the spectral number depends on the Morse function. Hence suppose that $f_-,f_+ \colon M \to \mathbb{R}$ are two Morse functions.
In order to construct the canonical isomorphisms between the two Morse homologies we have to interpolate between the two Morse functions. In order to get good estimates between the two spectral numbers we choose a special kind of interpolation.
Namely we pick a smooth monotone cutoff function $\beta \in C^\infty(\mathbb{R},[0,1])$ such that for $T>0$ it holds that
$$\beta(s)=\left\{\begin{array}{cc}
0 & s \leq -T\\
1 & s\geq T.
\end{array}\right.$$ 
Using the cutoff function $\beta$ we interpolate between $f_-$ and $f_+$ via the time dependent function
$$f_s=\beta(s) f_++(1-\beta(s))f_-, \quad s \in \mathbb{R}.$$
We further choose a Morse-Smale metric $g_-$ for $f_-$ and a Morse-Smale metric $g_+$ for $f_+$ and interpolate between
$g_-$ and $g_+$ by a smooth family of Riemannian metrics $g_s$. For the following estimate which improves Lemma~\ref{timeenest} the special kind of interpolation between $f_-$ and $f_+$ is crucial but it does not matter how we precisely interpolate between $g_-$ and $g_+$. 
\begin{lemma}\label{convactest}
Suppose that $x$ is a solution of the gradient flow equation for the time dependent gradient $\nabla_{g_s}f_s$ such that
$$\lim_{s \to \pm \infty}x(s)=x_\pm \in \mathrm{crit}(f_\pm).$$
Then
$$f_+(x_+) \leq f_-(x_-)+\max(f_+-f_-).$$
\end{lemma}
\textbf{Proof: } We differentiate and take advantage of the gradient flow equation
\begin{eqnarray*}
\frac{d}{ds}f_s(x(s))&=&\frac{d}{ds}\Big(\beta(s) f_++(1-\beta(s))f_-\Big)(x(s))\\
&=&\beta'(s)f_+(x(s))-\beta'(s)f_-(x(s))+df_s(x(s))\partial_s x(s)\\
&=&\beta'(s)(f_+-f_-)(x(s))-df_s(x(s))\nabla_{g_s} f_s(x(s))\\
&=&\beta'(s)(f_+-f_-)(x(s))-g_s\big(\nabla_{g_s}f_s(x(s)),\nabla_{g_s} f_s(x(s))\big)\\
&\geq&\beta'(s)(f_+-f_-)(x(s)).
\end{eqnarray*}
Integrating this inequality we obtain
\begin{eqnarray*}
f_+(x_+)-f_-(x_-)&=&\int_{-\infty}^\infty \frac{d}{ds}f_s(x(s))ds\\
&\leq& \int_{-\infty}^\infty \beta'(s)(f_+-f_-)(x(s))ds\\
&\leq&\int_{-\infty}^\infty\beta'(s)\max(f_+-f_-)ds\\
&=&\max(f_+-f_-)\int_{-\infty}^\infty \beta(s)ds\\
&=&\max(f_+-f_-).
\end{eqnarray*}
This implies the lemma. \hfill $\square$
\\ \\
In the following Corollary we use the canonical isomorphism 
$$\Phi \colon HM_*(f_-,g_-) \to HM_*(f_+,g_+)$$
as explained in Chapter~\ref{invariance}.
\begin{cor}
Suppose that $\alpha \neq 0 \in HM_*(f_-,g_-)$, then
$$\sigma_{f_+}(\Phi (\alpha)) \leq \sigma_{f_-}(\alpha)+\max(f_+-f_-).$$
\end{cor}\label{spe1}
\textbf{Proof: } Suppose that $\xi \in CM_*(f_-,g_-)$ is a representative of the homology class $\alpha$, i.e.,
$$[\xi]=\alpha.$$
We can write $\xi$ as formal sum
$$\xi=\sum_{c \in \mathrm{crit}(f_-)} a_c c$$
with coefficients $a_c \in \mathbb{Z}_2$.
Abbreviate by
$$\phi=\phi_{f_s, g_s} \colon CM_*(f_-) \to CM_*(f_+)$$
the chain map associated to the interpolation $(f_s,g_s)$ which induces the canonical isomorphism $\Phi$ on homology. 
By definition of $\phi$ we have
\begin{eqnarray*}
\phi(\xi)&=&\sum_{c \in \mathrm{crit}(f_-)} a_c \phi(c)\\
&=&\sum_{c \in \mathrm{crit}(f_-)} a_c\Big(\sum_{c' \in \mathrm{crit}(f_+)}\#_2\mathcal{N}(f_s,g_s;c,c')c'\Big)\\
&=&\sum_{c' \in \mathrm{crit}(f_+)}\Big(\sum_{c \in \mathrm{crit}(f_-)} a_c \#_2 \mathcal{N}(f_s,g_s;c,c')\Big)c'.
\end{eqnarray*}
For $c' \in \mathrm{crit}(f_+)$ abbreviate
$$b_{c'}:=\sum_{c \in \mathrm{crit}(f_-)} a_c \#_2 \mathcal{N}(f_s,g_s;c,c')$$
so that we can write the above formula more compactly as
$$\phi(x)=\sum_{c' \in \mathrm{crit}(f_+)} b_{c'} c'.$$
Suppose that
$$b_{c'} \neq 0$$
for a critical point $c'$ of $f_+$. Then there has to exist a critical point $c$ of $f_-$ such that
\begin{equation}\label{ac}
a_c \neq 0
\end{equation}
and in addition
$$\#_2\mathcal{N}(f_s,g_s;c,c') \neq 0.$$
In particular, there has to exist a gradient flow line of $\nabla_{g_s}f_s$ from $c$ to $c'$. 
Hence by Lemma~\ref{convactest} the estimate
$$f_+(c') \leq f_-(c)+\max(f_+-f_-)$$
holds true. Together with (\ref{ac}) this implies
\begin{eqnarray*}
\sigma_{f_+}(\phi(\xi))&=&\max\big\{f_+(c'): c' \in \mathrm{crit}(f_+),\,\,b_{c'} \neq 0\big\}\\
&\leq&\max\big\{f_-(c): c \in \mathrm{crit}(f_-),\,\,a_c \neq 0\big\}+\max(f_+-f_-)\\
&=&\sigma_{f_-}(\xi)+\max(f_+-f_-).
\end{eqnarray*}
Therefore
\begin{eqnarray*}
\sigma_{f_+}(\Phi(\alpha)) &\leq&\min_{\substack{\xi \in CM_*(f_-)\\ [\xi]=\alpha}}\sigma_{f_+}(\phi(\xi))\\
&\leq&\min_{\substack{\xi \in CM_*(f_-)\\ [\xi]=\alpha}}\sigma_{f_-}(\xi)+\max(f_+-f_-)\\
&=&\sigma_{f_-}(\alpha)+\max(f_+-f_-).
\end{eqnarray*}
This finishes the proof of the Corollary. \hfill $\square$
\\ \\
By taking of advantage of the inverse canonical isomorphism
$$\Psi=\Phi^{-1} \colon HM_*(f_+,g_+) \to HM_*(f_-,g_-)$$
we can obtain a lower bound on $\sigma_{f_+}(\Phi(\alpha))$ as well.
\begin{cor}\label{spe2}
Under the assumption of Corollary~\ref{spe2} we have
$$\sigma_{f_-}(\alpha)+\min(f_+-f_-) \leq \sigma_{f_+}(\Phi(\alpha)).$$
\end{cor}
\textbf{Proof: } Applying Corollary~\ref{spe1} the the inverse isomorphism $\Psi$ we obtain the estimate
\begin{eqnarray*}
\sigma_{f_-}(\alpha)&=&\sigma_{f_-}(\Psi \circ \Phi(\alpha))\\
&\leq& \sigma_{f_+}(\Phi(\alpha))+\max(f_--f_+)
\end{eqnarray*}
and therefore
\begin{eqnarray*}
\sigma_{f_+}(\Phi(\alpha))&\leq&\sigma_{f_-}(\alpha)-\max(f_--f_+)\\
&=&\sigma_{f_-}(\alpha)+\min(f_+-f_-).
\end{eqnarray*}
This proves the Corollary. \hfill $\square$
\\ \\
Combining Corollary~\ref{spe1} and Corollary~\ref{spe2} we proved the following proposition.
\begin{prop}\label{spe3}
Suppose that $\alpha \neq 0 \in HM_*(f_-,g_-)$, then
$$\big|\sigma_{f_-}(\alpha)-\sigma_{f_+}(\Phi(\alpha))\big| \leq ||f_+-f_-||_{C^0}=\max|f_+-f_-|.$$
\end{prop}
Since Morse homology $HM_*(f,g)$ depends as well on the choice of a Morse-Smale metric the spectral number a priori depends as well on $g$, so that we should write more precisely $\sigma_{f,g}(\alpha)$. However, an immediate consequence of Proposition~\ref{spe3} is that spectral numbers are independent of the choice of the Morse-Smale metric. 
\begin{cor}
Spectral numbers do not depend on the choice of the Morse-Smale metric.
\end{cor}
\textbf{Proof: } Suppose that $f$ is a Morse function and $g_-$ and $g_+$ are two Morse-Smale metrics for $f$. Then in view of Proposition~\ref{spe3} we have for $\alpha \neq 0 \in HM_*(f,g_-)$
$$\big|\sigma_{f,g_-}(\alpha)-\sigma_{f,g_+}(\Phi(\alpha))\big|\leq ||f-f||_{C^0}=0$$
implying that
$$\sigma_{f,g-}(\alpha)=\sigma_{f,g_+}(\Phi(\alpha)).$$
This finishes the proof of the Corollary. \hfill $\square$
\\ \\
In the following we identify Morse homologies for different Morse-Smale pairs via the canonical isomorphisms so that we
can speak of $HM_*$ without reference to the Morse-Smale pair. Denote by $C_{\mathrm{Morse}}$ the space of all Morse functions on $M$. Given $\alpha \neq 0 \in HM_*$ we can interpret spectral numbers as a map
$$\rho_\alpha \colon C_{\mathrm{Morse}} \to \mathbb{R}, \quad f \mapsto \sigma_f(\alpha).$$
Proposition~\ref{spe3} tells us that the map $\rho_\alpha$ is Lipschitz continuous with Lipschitz constant equal to one with respect to the $C^0$-topology on $C_{\mathrm{Morse}}$. Our goal is to use this fact to extend spectral invariants continuously to all smooth functions on $M$ not necessarily Morse. In order to achieve this goal we first show in the following paragraph that Morse functions are dense in the space of smooth functions. 

\subsection{Denseness of Morse functions}

Suppose that $M$ is a closed manifold. We denote by 
$$C^\infty=C^\infty(M,\mathbb{R})$$
the space of smooth functions on $M$ endowed with the $C^\infty$-topology and bv
$$C_{\mathrm{Morse}} \subset C^\infty$$
the subspace of Morse functions. The main result of this paragraph is the following proposition.
\begin{prop}\label{morsedense}
$C_{\mathrm{Morse}}$ is an open and dense subset in $C^\infty$.
\end{prop}
We prove the proposition by an analogous but much easier scheme as we proved denseness of Morse-Smale metrics in Chapter~\ref{transversality}. For a prove of the proposition only requiring the finite dimensional version of Sard's theorem, see
\cite[Chapter 6]{milnor}.
\\ \\
\textbf{Proof of Proposition~\ref{morsedense}: }Given $f \in C^\infty$ define a section
$$s_f \colon M \to T^*M, \quad x \mapsto df(x).$$
Zeros of this section are precisely critical points. If $x \in s_f^{-1}(0)$, then the vertical differential at $x$
$$Ds_f(x) \colon T_x M \to T^*_xM$$ is given by
\begin{equation}\label{vertdiff}
Ds_f(x)v=d^2 f(x)(v,\cdot), \quad v \in T_xM.
\end{equation}
As a linear map between finite dimensional vector spaces it is obviously Fredholm and since the two vector spaces have the same dimension its index is zero. Hence $Ds_f(x)$ is surjective if and only if it is injective. Moreover, by (\ref{vertdiff})
it is injective if and only if $x$ is a Morse critical point. Therefore
$$s_f \pitchfork 0\quad \Longleftrightarrow \quad f \in C_{\mathrm{Morse}}.$$ 
For $k \geq 2$ abbreviate
$$C^k=C^k(M,\mathbb{R})$$
the Banach space of $C^k$-functions on $M$ endowed with the $C^k$-topology. Define a section
$$S \colon C^k \times M \to T^*M, \quad (f,x) \mapsto df(x).$$
Suppose that 
$$(f,x) \in S^{-1}(0),$$
i.e. $x$ is a critical point of $f$. The vertical differential of $S$ at $(f,x)$ is given by
$$DS(f,x) \colon C^k \times T_x M \to T^*_x M, \quad (\widehat{f},\widehat{x}) \mapsto d\widehat{f}(x)+d^2 f(x)(\widehat{x},\cdot).$$
Given a cotangent vector $w \in T_x^*M$ we can always find a function $\widehat{f} \in C^k$ such that
$$d\widehat{f}(x)=w.$$
Therefore the vertical differential $DS(f,x)$ is surjective and the universal moduli space
$$S^{-1}(0) \subset C^k \times M$$
consisting of critical points of arbitrary functions is a Banach manifold. We consider the map
$$\Pi \colon S^{-1}(0) \to C^k, \quad (f,x) \mapsto f.$$
By Sard's theorem there exists a subset of second category
$$C^k_{\mathrm{reg}} \subset C^k$$
consisting of regular values of $\Pi$. If $f \in C^k_{\mathrm{reg}}$, then $s_f$ is transverse to the zero section and
hence $f$ is a Morse function of class $C^k$. Since $M$ is closed the set of critical points is compact and therefore
the set $C^k_{\mathrm{reg}}$ is open as well. It follows that
$$C_{\mathrm{Morse}}=C^k_{\mathrm{reg}} \cap C^\infty$$
is dense in $C^k$ for every $k \geq 2$. In particular, $C_{\mathrm{Morse}}$ is dense in $C^\infty$. Using once more that the set of critical points on a closed manifold is compact it is open as well. This finishes the proof of the proposition. 
\hfill $\square$

\subsection{Extension of spectral numbers to smooth functions}

Given $\alpha \neq 0 \in HM_*$ we constructed so far a map
$$\rho_\alpha \colon C_{\mathrm{Morse}} \to \mathbb{R}, \quad f \mapsto \sigma_f(\alpha)$$
which by Proposition~\ref{spe3} is Lipschitz continuous with Lipschitz constant one with respect to the $C^0$-topology
on $C_{\mathrm{Morse}}$, i.e., if $f$ and $f'$ are Morse functions on $M$, then
\begin{equation}\label{lipschitz}
\big|\rho_\alpha(f)-\rho_\alpha(f')\big| \leq ||f-f'||_{C^0}.
\end{equation}
We explain how (\ref{lipschitz}) gives rise to a continuous extension of $\rho_\alpha$ to smooth functions which are not necessarily Morse. Suppose that $f \in C^\infty$. By Proposition~\ref{morsedense} there exists a sequence 
$f_\nu \in C_{\mathrm{Morse}}$ for $\nu \in \mathbb{N}$ such that 
$$f_\nu \xrightarrow
 {
  \substack{C^\infty}
 }
 f.
$$
In particular,
$$f_\nu \xrightarrow
 {
  \substack{C^0}
 }
 f.
$$
This implies that $f_\nu$ is a Cauchy sequence and therefore there exists $\nu_0 \in \mathbb{N}$ such that for every
$\nu,\mu \geq \nu_0$ it holds that
$$||f_\nu-f_\mu||_{C^0} \leq 1.$$
Hence we obtain for $\nu \geq \nu_0$ by the Lipschitz continuity (\ref{lipschitz})
$$|\rho_\alpha(f_\nu)| \leq |\rho_\alpha(f_{\nu_0})|+||f_\nu-f_{\nu_0}||_{C^0} \leq |\rho_\alpha(f_{\nu_0})|+1.$$
In particular, $|\rho_\alpha(f_\nu)|$ is uniformly bounded so that there exists a subsequence $\nu_j$ such that
$$\lim_{j \to \infty}\rho_\alpha(f_{\nu_j}) \in \mathbb{R}$$
exists. We set
$$\rho_\alpha(f):=\lim_{j \to \infty}\rho_\alpha(f_{\nu_j}).$$
\begin{lemma}
$\rho_\alpha(f)$ is well defined, i.e., independent of the choice of the sequence $f_{\nu_j}$. 
\end{lemma}
\textbf{Proof: } Assume that $(f_\nu)_\nu$ and $(f'_\nu)_\nu$ are two sequences meeting the following requirements
\begin{description}
 \item[(i)] $f_\nu, f'_\nu \in C_{\mathrm{Morse}}$ for every $\nu \in \mathbb{N}$,
 \item[(ii)] $f_\nu \xrightarrow
 {
  \substack{C^0}
 }
 f
$ and $f'_\nu \xrightarrow
 {
  \substack{C^0}
 }
 f,
$
 \item[(iii)] $\lim_{\nu \to \infty} \rho_\alpha(f_\nu)$ exists as well as $\lim_{\nu \to \infty} \rho_\alpha(f'_\nu)$.
\end{description}
Under this assumptions we show that
$$\lim_{\nu \to \infty} \rho_\alpha(f_\nu)=\lim_{\nu \to \infty} \rho_\alpha(f'_\nu).$$
To see that pick $\epsilon>0$. By assumption (ii) there exists $\nu_0$ such that
$$||f-f_\nu||_{C^0} \leq \frac{\epsilon}{6},\quad ||f-f'_\nu||_{C^0} \leq \frac{\epsilon}{6},\qquad \nu \geq \nu_0.$$
Therefore we have for $\nu,\mu \geq \nu_0$
$$||f_\nu-f'_\mu||_{C^0}\leq||f-f_\nu||_{C^0}+||f-f'_\mu||_{C^0} \leq \frac{\epsilon}{3}$$
so that we obtain from the Lipschitz continuity (\ref{lipschitz}) 
\begin{equation}\label{ext1}
|\rho_\alpha(f_\nu)-\rho_\alpha(f'_\mu)| \leq ||f_\nu-f'_\mu||_{C^0} \leq \frac{\epsilon}{3}.
\end{equation}
Abbreviate
$$a:=\lim_{\nu \to \infty} \rho_\alpha(f_\nu),\qquad b:=\lim_{\nu \to \infty} \rho_\alpha(f'_\nu).$$
Choose $\nu \geq \nu_0$ such that
\begin{equation}\label{ext2}
|\rho_\alpha(f_\nu)-a| \leq \frac{\epsilon}{3}
\end{equation}
and choose $\mu \geq \nu_0$ such that
\begin{equation}\label{ext3}
|\rho_\alpha(f'_\mu)-b| \leq \frac{\epsilon}{3}.
\end{equation}
Using (\ref{ext1}),(\ref{ext2}), and (\ref{ext3}) we estimate
\begin{eqnarray*}
|a-b| &\leq&|a-\rho_\alpha(f_\nu)|+|\rho_\alpha(f_\nu)-\rho_\alpha(f'_\mu)|+|\rho_\alpha(f'_\mu)-b| \leq \epsilon.
\end{eqnarray*}
Since $\epsilon$ was arbitrary it follows that
$$a=b$$
and the lemma is proved. \hfill $\square$
\\ \\
In view of the lemma the function $\rho_\alpha$ has a well defined extension to smooth functions. By abuse of notation we still denote the extension by the same letter so that we now map a map
$$\rho_\alpha \colon C^\infty \to \mathbb{R}.$$
\subsection{Lipschitz continuity of the extension}

By (\ref{lipschitz}) we know that before the extension the function $\rho_\alpha$ was Lipschitz continuous. We show in this paragraph that the estimate (\ref{lipschitz}) continuous to hold for the extended function.
\begin{lemma}\label{liplem}
Suppose that $f,f' \in C^\infty$. Then
$$\big|\rho_\alpha(f)-\rho_\alpha(f')\big| \leq ||f-f'||_{C^0}.$$
\end{lemma}
\textbf{Proof: }If $f$ and $f'$ are Morse functions this is precisely (\ref{lipschitz}). Now assume that $f$ and $f'$ are arbitrary functions. Choose sequence $(f_\nu)_\nu$ and $(f'_\nu)_\nu$ such that 
\begin{description}
 \item[(i)] $f_\nu, f'_\nu \in C_{\mathrm{Morse}}$ for every $\nu \in \mathbb{N}$,
 \item[(ii)] $f_\nu \xrightarrow
 {
  \substack{C^0}
 }
 f
$ and $f'_\nu \xrightarrow
 {
  \substack{C^0}
 }
 f.
$
\end{description}
Pick $\epsilon>0$. Then there exists $\nu \in \mathbb{N}$ such that
$$|\rho_\alpha(f)-\rho_\alpha(f_\nu)|\leq \frac{\epsilon}{4}, \quad
|\rho_\alpha(f')-\rho_\alpha(f'_\nu)|\leq \frac{\epsilon}{4}$$
as well as
$$||f-f_\nu||_{C^0}\leq \frac{\epsilon}{4},
\quad ||f'-f'_\nu||_{C^0}\leq \frac{\epsilon}{4}.$$
We estimate using the Lipschitz estimate for Morse functions
\begin{eqnarray*}
|\rho_\alpha(f)-\rho_\alpha(f')|&\leq&|\rho_\alpha(f)-\rho_\alpha(f_\nu)|+|\rho_\alpha(f_\nu)-\rho_\alpha(f'_\nu)|
+|\rho_\alpha(f'_\nu)-\rho_\alpha(f')|\\
&\leq&\frac{\epsilon}{2}+||f_\nu-f'_\nu||_{C^0}\\
&\leq&\frac{\epsilon}{2}+||f_\nu-f||_{C^0}+||f-f'||_{C^0}+||f'-f'_\nu||_{C^0}\\
&\leq&\epsilon+||f-f'||_{C^0}.
\end{eqnarray*}
Since $\epsilon>0$ was arbitrary we conclude that
$$\big|\rho_\alpha(f)-\rho_\alpha(f')\big| \leq ||f-f'||_{C^0}$$
and the lemma is proved. \hfill $\square$

\subsection{Spectrality}

Suppose that $f \in C^\infty$. Abbreviate
$$\mathscr{S}_f:=\big\{r \in \mathbb{R}: r=f(x),\,\,x \in \mathrm{crit}(f)\big\}$$
the set of critical values or \emph{spectrum} of $f$. We next show that $\rho_\alpha(f)$ lies in the spectrum of $f$.
\begin{lemma}\label{spectral}
Assume that $\alpha \neq 0 \in HM_*$ and $f \in C^\infty$. Then
$$\rho_\alpha(f) \in \mathscr{S}_f.$$
\end{lemma}
\textbf{Proof: } If $f$ is Morse, this is clear by construction of the spectral number. For the general case we again approximate $f$ by Morse functions. Namely we choose a sequence $(f_\nu)_\nu$ such that
\begin{description}
 \item[(i)] $f_\nu\in C_{\mathrm{Morse}}$ for every $\nu \in \mathbb{N}$,
 \item[(ii)] $f_\nu \xrightarrow
 {
  \substack{C^\infty}
 }
 f.
$
\end{description}
By spectrality for Morse functions we conclude that there exist $x_\nu \in \mathrm{crit}(f_\nu)$ such that
$$f_\nu(x_\nu)=\rho_\alpha(f_\nu).$$
Since $M$ is closed there exists a subsequence $\nu_j$ and $x \in M$ such that
$$\lim_{j \to \infty}x_{\nu_j}=x.$$
Since the points $x_{\nu_j}$ are critical points of $f_{\nu_j}$ they are solutions of the equation
$$df_{\nu_j}(x_{\nu_j})=0$$
and since $f_{\nu_j}$ converges to $f$ in the $C^\infty$-topology we conclude that
$$df(x)=0,$$
i.e.,
$$x \in \mathrm{crit}(f).$$
Therefore
\begin{eqnarray*}
f(x)&=&\lim_{j \to \infty}f_{\nu_j}(x_{\nu_j})=\lim_{j \to \infty}\rho_\alpha(f_{\nu_j})=\rho_\alpha(f)
\end{eqnarray*}
and the lemma follows. \hfill $\square$
\\ \\
If we abbreviate
$$\mathscr{S}=\big\{(f,r): f \in C^\infty,\,\,r \in \mathscr{S}_f\big\} \subset C^\infty \times \mathbb{R}$$
we have a natural projection
$$\pi \colon \mathscr{S} \to C^\infty, \quad (f,r) \mapsto f$$
so that we can think of $\mathscr{S}$ as a bundle over $C^\infty$, which we refer to as the \emph{spectral bundle}.
The fibre over a point $f \in C^\infty$ is then given by
$$\pi^{-1}(f)=\mathscr{S}_f.$$
Using this notion we can interpret by Lemma~\ref{spectral} the functions $\rho_\alpha$ for $\alpha \neq 0 \in HM_*$ as
sections
$$\rho_\alpha \colon C^\infty \to \mathscr{S}$$
which are by Lemma~\ref{liplem} Lipschitz continuous with Lipschitz constant one with respect to the $C^0$-topology.

\subsection{The homological spectral bundle and the action gap}

If $f \in C^\infty$ we abbreviate
$$\mathscr{S}_f^h:=\big\{\rho_\alpha(f): \alpha \neq 0 \in HM_*\big\} \subset \mathscr{S}_f$$
and set
$$\mathscr{S}^h=\big\{(f,r): f \in C^\infty,\,\,r \in \mathscr{S}^h_f\big\} \subset \mathscr{S}$$
which we refer to as the \emph{homological spectral bundle}. We can think of the homological spectral bundle as a subbundle of the spectral bundle. Note that for each $f \in C^\infty$ we have an obvious string of inequalities
\begin{equation}\label{critest}
\#\mathrm{crit}(f) \geq \#\mathscr{S}_f \geq \#\mathscr{S}_f^h.
\end{equation}
Therefore an estimate on the cardinality of points in the fiber of the homological spectral bundle gives a lower bound on
the number of critical points of the corresponding function which does not need to be Morse. For $f \in C^\infty$ we introduce the action gap as follows 
$$\Delta(f):=\min\big\{|r-r'|: r \neq r' \in \mathscr{S}^h_f\big\}.$$
We have the following theorem.
\begin{thm}
Suppose that $f$ and $f'$ are smooth functions on $M$ satisfying
$$||f-f'||_{C^0}<\frac{\Delta(f)}{2}.$$
Then 
$$\#\mathrm{crit}(f') \geq \# \mathscr{S}^h_f.$$
\end{thm}
\textbf{Proof: } Since the difference of $f$ and $f'$ in the $C^0$-norm is less than half the action gap of $f$ it follows from the Lipschitz continuity of the spectral numbers in Lemma~\ref{liplem} that
$$\#\mathscr{S}_{f'}^h \geq \#\mathscr{S}_f^h.$$
Hence in view of (\ref{critest}) we obtain the estimate
$$\#\mathrm{crit}(f') \geq \#\mathscr{S}_{f'}^h \geq \#\mathscr{S}_f^h.$$
This proves the Theorem. \hfill $\square$


\begin{thebibliography}{99}
\bibitem{abraham-marsden-ratiu} R.\,Abraham, J.\,Marsden, T.\,Ratiu, \emph{Manifolds, tensor analysis, and applications,}
2nd ed., Applied Mathematical Sciences \textbf{75}, Springer (1988).
\bibitem{albers-frauenfelder} P.\,Albers, U.\,Frauenfelder, \emph{Exponential decay for sc-gradient flow lines}, J.\,Fixed Point Theory Appl. \textbf{13}, no.\,2, (2013), 571--586.
\bibitem{audin-damian} M.\,Audin, M.\,Damian, \emph{Morse theory and Floer homology}, Translated from the 2010 French original by Reinie Ern\'e, Universitext, Springer, London, EDP, Sciences, Les Ulis (2014).
\bibitem{banyaga-hurtubise} A.\,Banyaga, D.\,Hurtubise, \emph{Lectures on Morse homology}, Kluwer Texts in Mathematical Sciences \textbf{29}, Kluwer Academic Publishers Group, Dordrecht (2004).
\bibitem{floer1} A.\,Floer, \emph{Morse theory for Lagrangian intersections}, J.\,Diff.\,Geom. \textbf{28} (1988), 513--547.
\bibitem{floer2} A.\,Floer, \emph{The unregularized gradient flow of the symplectic action}, Comm.\,Pure Appl.\,Math.
 \textbf{41} (1988), 775--813.
\bibitem{floer3} A.\,Floer, \emph{Witten's complex and infinite dimensional Morse theory}, J.\,Diff.\,Geom. \textbf{30} (1989), 575--611. 
\bibitem{floer-hofer} A.\,Floer, H.\,Hofer, \emph{Coherent orientations for periodic orbit problems in symplectic geometry},
 Math.\,Z. \textbf{212} (1993), 13--38.
\bibitem{gaio-salamon} A.\,Gaio, D.\,Salamon, \emph{Gromov-Witten invariants of symplectic quotients and adiabatic limits},
 J.\,Symplectic Geom. \textbf{3} (2005), 55--159.
\bibitem{mcduff-salamon} D.\,McDuff, D.\,Salamon, \emph{J-holomorphic Curves and Symplectic
Topology} 2nd edition, Amer.\,Math.\,Soc., Providence, RI (2012).
 \bibitem{milnor} J.\,Milnor, \emph{Morse theory. Based on lecture notes by M.\,Spivak and R.\,Wells}, Annals of Mathematics Studies \textbf{51}, Princeton University Press, Princeton, N.J. (1963).
 \bibitem{milnor2} J.\,Milnor, \emph{Topology from the Differential Viewpoint}, The University Press of Virginia (1965).
 \bibitem{nicolaescu} L.\,Nicolaescu, \emph{An invitation to Morse theory}, Second edition. Universitext. Springer, New York
  (2011).
 \bibitem{palis-demelo}J.\,Palis, W.\,de Melo, \emph{Geometric theory of dynamical systems}, Springer-Verlag, 
 New York-Berlin (1982).
 \bibitem{polterovich-rosen-samvelyan-zhang} L.\,Polterovich, D.\,Rosen, K.\,Samvelyan, J.\,Zhang, \emph{Topological 
  Persistence in Geometry and Analyis}, arXiv:1904.04044
\bibitem{robbin-salamon} J.\,Robbin, D.\,Salamon, \emph{The spectral flow and the Maslov index}, Bull.\,London Math.\,
 Soc. \textbf{27}, no.\,1 (1995), 1--33.
 \bibitem{salamon} D.\,Salamon, \emph{Lectures on Floer homology}, Symplectic geometry and topology (Park City, UT, 1997),
  147--229, IAS/Park City Math.\,Ser.,\,\textbf{7}, Amer.\,Math.\,Soc., Providence, RI (1999).
 \bibitem{schwarz} M.\,Schwarz, \emph{Morse homology}, Progress in Mathematics \textbf{111}, Birkh\"auser Verlag, Basel
 (1993).
 \bibitem{schwarz2} M.\,Schwarz, \emph{Equivalences for Morse homology}, in: Geometry and Topology in Dynamics (Winston-
  Salem, NC, 1998/San Antonio, TX, 1999), Contemporary Mathematics, vol.\,\textbf{246} (1999), American Mathematical  
  Society, Providence, RI, 1999), 197--216. 
 \bibitem{simcevic} T.\,Simcevic, \emph{A Hardy Space Approach to Lagrangian Floer Gluing}, arXiv:1410.5998.
 \bibitem{sylvester} J.\,Sylvester, \emph{A demonstration of the theorem that every homogeneous quadratic polynomial is
 reducible by real orthogonal substitutions to the form of a sum of positive and negative squares}, Philosophical Magazine
 4th Series \textbf{4}(23) (1852), 138--142.
 \bibitem{weber} J.\,Weber, \emph{The Morse-Witten complex via dynamical systems}, Expo.\,Math. \textbf{24}, no.\,2 (2006),
  127--159.
 \bibitem{zehnder} E.\,Zehnder, \emph{Lectures on dynamical systems. Hamiltonian vector fields and symplectic capacities},
 EMS Textbooks in Mathematics, European Mathematical Society (EMS), Z\"urich (2010).
 \bibitem{ziltener} F.\,Ziltener, \emph{The invariant symplectic action and decay for vortices}, J.\,Symplectic Geom.
 \textbf{7} (2009), 357--376.
\end{thebibliography}
\end{document}